\documentclass[12pt, a4paper, reqno]{amsart}

\usepackage{epsfig}
\usepackage{graphicx}
\usepackage{setspace}
\usepackage{mathrsfs}
\usepackage{cancel}
\usepackage{faktor}
\usepackage{mhequ}
\usepackage{stmaryrd}
\usepackage{amssymb}
\usepackage{version}
\usepackage[T1]{fontenc}
\usepackage{amsmath}
\usepackage{amsthm}

\newcommand{\cqfd}{\; \;\underset{\Box}{}}
\newcommand{\R}{\textrm{I\hspace{-.19em}R}}

\newcommand{\e}{\mathrm{e}}

\newcommand{\w}{\widetilde}
\renewcommand{\t}{\tilde}
\newcommand{\wg}{\tilde{g}}
\newcommand{\nw}{\tilde{\nabla}}
\newcommand{\contrainte}{\mathcal{C}}
\renewcommand{\d}{\partial}
\newcommand{\Az}{\mathring{\A}}
\newcommand{\Bz}{\mathring{B}}
\newcommand{\Ps}{P^{*}}
\newcommand{\Pgpi}{P_{(g,\pi)}}
\newcommand{\Pgpis}{P_{(g,\pi)}^{*}}
\newcommand{\Psw}{\widetilde{P}^{*}}
\newcommand{\Pgpisw}{P_{(\t{g},\t{\pi})}^{*}}
\newcommand{\Sz}{\mathring{S}}
\newcommand{\Uz}{\mathring{U}}
\newcommand{\phibf}{\mathbf{\Phi}}
\newcommand{\phibfs}{\mathbf{\Phi}^{*}}
\newcommand{\phisun}{\mathbf{\Phi}^{*}_{1}}
\newcommand{\phisdeux}{\mathbf{\Phi}^{*}_{2}}
\newcommand{\phiwsun}{\widetilde{\mathbf{\Phi}}^{*}_{1}}
\newcommand{\phiwsdeux}{\widetilde{\mathbf{\Phi}}^{*}_{2}}
\newcommand{\phio}{\mathbf{\Phi_{0}}}

\newcommand{\phiog}{\mathbf{\Phi_{0}}(g)}
\newcommand{\phiogpi}{\mathbf{\Phi_{0}}(g, \pi)}

\newcommand{\phigpi}{\mathbf{\Phi}(g, \pi)}
\newcommand{\phigpiws}{\mathbf{\Phi}(\t{g},\t{\pi})^{*}}
\newcommand{\phiogs}{\mathbf{\Phi_{0}}(g)^{*}}

\newcommand{\phiogpis}{\mathbf{\Phi_{0}}(g,\pi)^{*}}
\newcommand{\phiigpis}{\mathbf{\Phi_{i}}(g,\pi)^{*}}

\newcommand{\phigpis}{\mathbf{\Phi}(g,\pi)^{*}}
\newcommand{\phii}{\mathbf{\Phi_{i}}}

\newcommand{\phiigpi}{\mathbf{\Phi_{i}}(g,\pi)}
\newcommand{\nat}{\mathbb{N}}
\newcommand{\G}{\mathcal{G}}
\newcommand{\Gplus}{\mathcal{G}^{+}}
\renewcommand{\L}{\mathcal{L}}
\renewcommand{\S}{\mathcal{S}}
\newcommand{\T}{\mathcal{T}}
\newcommand{\K}{\mathcal{K}}
\newcommand{\F}{\mathcal{F}}

\newcommand{\precig}[1]{\hspace{0.05cm}^{(#1)}\!}
\newcommand{\preciv}[1]{\hspace{0.05cm}^{#1}\!}
\renewcommand{\epsilon}{\varepsilon}
\newcommand{\Epsilon}{\t{\varepsilon}}

\newcommand{\cde}[3]{\Gamma_{#1 \; #3}^{\hspace{0.1cm} #2}}
\newcommand{\wcde}[3]{\widetilde{\Gamma}_{#1 \; #3}^{\hspace{0.1cm} #2}}
\newcommand{\ape}[1]{A_{#1}}
\newcommand{\ade}[3]{A_{#1 \; #3}^{\hspace{0.1cm} #2}}
\newcommand{\wade}[3]{\widetilde{A}_{#1 \; #3}^{\hspace{0.1cm} #2}}

\newcommand{\cdez}[3]{\mathring{\Gamma}_{#1 \; #3}^{\hspace{0.1cm} #2}}

\newcommand{\infeg}{\leqslant}
\newcommand{\A}{\mathcal{A}}
\newcommand{\supeg}{\geqslant}
\newcommand{\equi}{\Leftrightarrow}
\newcommand{\cinf}{\mathscr{C}^{\infty}}

\newcommand{\demi}{\frac{1}{2}}
\newcommand{\tdemi}{\tfrac{1}{2}}

\newcommand{\Norm}[2]{||{#1}||_{#2}}

\newcommand{\bNorm}[2]{\Big\| {#1} \Big\|_{#2}}
\newcommand{\BNorm}[2]{\bigg\| {#1} \bigg\|_{#2}}
\newcommand{\norm}[2]{|{#1}|_{#2}}
\newcommand{\normgz}[1]{|{#1}|_{\gz}}
\newcommand{\leb}[2]{{L}^{#1}_{#2}}
\newcommand{\lebd}{L^{2}_{\delta}}
\newcommand{\sob}[3]{W^{#1,#2}_{#3}}
\newcommand{\sobd}[1]{W^{#1,2}_{\delta}}
\newcommand{\m}{\mathcal{M}}
\newcommand{\N}{\mathcal{N}}

\newcommand{\ie}{\emph{i.e.} \,}
\newcommand{\gz}{\mathring{g}}
\newcommand{\hz}{\mathring{h}}
\newcommand{\Kz}{\mathring{K}}
\newcommand{\piz}{\mathring{\pi}}
\newcommand{\n}{\nabla}
\newcommand{\nz}{\mathring{\nabla}}
\newcommand{\Tz}{\mathring{T}}
\newcommand{\wn}{\widetilde{\nabla}}
\newcommand{\wlapla}{\widetilde{\Delta}}
\newcommand{\laplaz}{\mathring{\Delta}}

\newcommand{\laplag}{\Delta_{g}}

\newcommand{\laplahz}{\Delta_{\mathring{h}}}
\newcommand{\pshz}{\, \textbf{.}_{\mathring{h}}}

\newcommand{\quotient}[2]{#1_{\! \diagup_{\! #2}}}

\newcommand{\Riem}{\operatorname{Riem}}
\newcommand{\Ric}{\operatorname{Ric}}
\renewcommand{\div}{\operatorname{div}}
\newcommand{\tr}{\operatorname{tr}}
\newcommand{\Coker}{\operatorname{Coker}}
\renewcommand{\Im}{\operatorname{Im}}

\def\restriction#1#2{\mathchoice
                  {\setbox1\hbox{${\displaystyle #1}_{\scriptstyle #2}$}
                  \restrictionaux{#1}{#2}}
                  {\setbox1\hbox{${\textstyle #1}_{\scriptstyle #2}$}
                  \restrictionaux{#1}{#2}}
                  {\setbox1\hbox{${\scriptstyle #1}_{\scriptscriptstyle #2}$}
                  \restrictionaux{#1}{#2}}
                  {\setbox1\hbox{${\scriptscriptstyle #1}_{\scriptscriptstyle #2}$}
                  \restrictionaux{#1}{#2}}}
    \def\restrictionaux#1#2{{#1\,\smash{\vrule height .8\ht1 depth .85\dp1}}_{\,#2}}

\newtheorem{lemperso}{Lemma}
\newtheorem{thmperso}{Theorem}
\newtheorem{corperso}{Corollary}
\newtheorem{propperso}{Proposition}
\newtheorem{defperso}{Definition}

\title[Hilbert manifold structure for AH initial data]{Hilbert manifold structure for asymptotically hyperbolic relativistic initial data}

\author[E.  Delay]{Erwann
Delay} \address{Erwann Delay, Avignon Université, Laboratoire de Mathématiques d'Avignon (EA 2151)
F-84916 Avignon}
\email{Erwann.Delay@univ-avignon.fr}
\urladdr{http://www.math.univ-avignon.fr}

\author[J. Fougeirol]{Jérémie Fougeirol}
\address{Jérémie Fougeirol\\
Avignon Université, Laboratoire de Mathématiques d'Avignon (EA 2151)
F-84916 Avignon}
\address{Institut Montpelliérain Alexander Grotendieck\\
CNRS and Montpellier Université\\
Place Eugène Bataillon\\
34095 Montpellier Cedex 5, France.}
\email{jeremie.fougeirol@free.fr}
\date{December 12,  2016}
\usepackage{xcolor}
\usepackage{cite}

\newcommand{\corec}[1]{{ #1}}

\begin{document}

\begin{abstract}
We provide a Hilbert manifold structure à la Bartnik for the space of asymptotically hyperbolic initial data for the vacuum constraint equations. The adaptation led us to prove
new weighted Poincaré and Korn type inequalities for AH manifolds with inner boundary
and weakly regular metric.
\end{abstract}

\maketitle

\tableofcontents


%\bf 2010 MSC : 53C21, 53A45, 53A30,  58J05, 35J61

\section{Introduction}

This work follows on from a paper of $2005$ \cite{Bartnik2005} in which \corec{R.} Bartnik described a Hilbert manifold structure for the space of asymptotically flat solutions of the Einstein equations (see also \cite{McCormick2014},\cite{McCormick2015}, \cite{RaiSaraykar2016}). The work is done with rather weak regularity assumptions concerning the metric involved (curvature constant only modulo weighted $\leb{2}{}$ terms) and so can be related to the context of the bounded $L^{2}$ curvature conjecture of KRS \cite{SergiuKlainerman2013}. Actually \corec{R.} Bartnik showed in \cite{Bartnik1986} that these assumptions on the regularity are the weakest possible to define the ADM mass of the manifold, explaining why we chose the same regularity conditions. In an undergoing work \corec{\cite{JFthese}}, we are showing that these very assumptions allow us to properly define the mass of an asymptotically hyperbolic manifold compatible with previous definitions of the mass (see \cite{ChruscielHerzlich2003} , \cite{Herzlich2005} or \cite{DGS2013}) and with the Hilbert manifold structure exposed here.
In order to overcome difficulties arising in the asymptotically hyperbolic case, we had to create a Hessian-type operator $\Tz$ and a differential operator of order two, called $\Uz$, built up with the first  derivatives of the Killing operator $\Sz$. In particular, we obtain Poincaré and Korn-type estimates of second order on an asymptotically hyperbolic manifold with boundary. These estimates are the key to prove triviality of the adjoint's kernel.\\

{\bf Acknowledgement :} J\'er\'emie Fougeirol thanks Marc Herzlich for some helpful remarks.
\section{Notations and conventions}
%\addcontentsline{toc}{chapter}{Notations and conventions}
%$\d$ désigne l'opérateur de dérivation partielle usuelle.\\
Let $(\m,g)$ be a Riemannian manifold. We define  $T^{r}_{m}(\m)$
to be the  bundle of tensor covariant of  rank $m$ and contravariant of  rank $r$ .
For all $u \in T^{r}_{m}(\m)$, $\norm{u}{g}$ will denote the norm of $u$ with respect to the metric $g$ and notation $\norm{u}{g,x}$ allows us to precise the point of the manifold we consider. $d \mu (g)$ is the Riemannian measure determined by $g$. $\Riem \, g , \Ric \, g$ and $R(g)$ are respectively the Riemann tensor, the Ricci tensor and the scalar curvature of the metric $g$. For a Riemannian metric $g$ with connection $\n$, we set the following notations concerning the Hessian and Laplacian of a function $u$.
\begin{eqnarray*}
\n^{2}_{ij}u &=& \n_{i}\n_{j}u\\
\Delta u &=& \corec{\tr}_{g}\n^{2}u = g^{ij} \n^{2}_{ij}u
\end{eqnarray*}
w.r.t is the abbreviation of "with respect to".\\
The work presented in this article is done on a $n$-dimensional manifold as often as possible and, otherwise explicitely stated, the results presented here are valid in any dimension $n$. Nevertheless, Sobolev inequalities strongly constrain the dimension to $n=3$ in several proofs ; this is clearly specified when this is the case. We chose to leave $n$ and specify $n=3$ in the concerned results rather than replace $n$ by its specific value since it helps to understand where the dimension plays a role.\\
Concerning constants in norm inequalities, the constant \corec{$c$} will design in general a constant depending on the background metric $\gz$ and the decaying rate $\delta$ and its expression may change from line to line in a proof. The nature of the dependence of the constant $C$ will be systematically specified because it will depend on other parameters.

\subsection{Conformally compact manifold.}$\mbox{ }$\\
Let $(\m, \gz)$ be a $\cinf$ $n$-dimensional complete non compact Riemannian manifold.
The manifold $(\m, \gz)$ is conformally compact if there exists a Riemannian metric $\hz$ such that $\gz = \rho^{-2} \hz$ , where $(\m , \hz)$ is a $\cinf$ compact Riemannian manifold with boundary $\corec{\d_{\infty}\m}$ and $\rho$ is a function on $\bar{\m}$, called defining function on $\corec{\d_{\infty}\m}$:
\begin{itemize}
\item[\textbullet] $\rho \in \cinf(\bar{\m})$
\item[\textbullet] $\rho \supeg 0$ on $\bar{\m}$
\item[\textbullet] $\corec{\d_{\infty}\m} = \lbrace x \in \m : \rho(x)=0\rbrace$
\item[\textbullet] $\restriction{d\rho}{\corec{\d_{\infty}\m}}$ \corec{never vanishes}.
\end{itemize}
We call asymptotically hyperbolic metric on $\m$ every Riemannian metric $\gz$ verifying:
\begin{itemize}
\item[\textbullet] $\gz$ is conformally compact.
\item[\textbullet] $\norm{d\rho}{\hz, \corec{\d_{\infty}\m}}^{2} = 1.$
\end{itemize}

\corec{In some inequalities, we will allow $\m$ to have  an inner boundary $\d \m$. In that case,
$\bar{\m}$ will the disjoint union of $M$, $\d \m$ and $\d_\infty\m$.}

The existence of these two metrics on $\m$ can be source of confusion thereafter, so we have to fix some notations.
$\nz, \laplaz, \cdez{i j}{k} \; , \norm{u}{\gz}$ (resp. $\hspace{0.05cm}^{(\hz)}\n, \laplahz, \hspace{0.05cm}^{(\hz)}\cde{i j}{k} \; , \norm{u}{\hz})$ will respectively denote the connection, Laplace operator, Christoffel symbols and tensor norm w.r.t $\gz$ (resp. $\hz$).
There are some correspondances between these quantities.\\
For volume measures, in dimension $n$, clearly $\; \; d\mu (\gz) = \rho^{-n} d\mu (\hz).$\\
For Christoffel symbols, $\; \;\cdez{i j}{k} = \hspace{0.05cm}^{(\hz)}\cde{i}{k}{j} \; - \frac{1}{\rho} \left( \delta^{k}_{j} \d_{i} \rho + \delta^{k}_{i} \d_{j} \rho - \hz^{kl} \hz_{ij} \d_{l} \rho \right). $

In particular, for the Hessian of $\rho$ (e.g. equation (C.12) from \cite{Chrusciel2004}) %gives
\begin{equation}\label{hessrhoah}
\nz^{2}_{ij} \rho = \precig{\hz}\nz^{2}_{ij} \rho + \rho \bigg(2 \frac{\d_{i}\rho}{\rho} \frac{\d_{j}\rho}{\rho} - \gz_{ij}\norm{d\rho}{\hz}^{2}\bigg).
\end{equation}
Taking the trace of (\ref{hessrhoah}), we obtain the expression of the Laplace operator of $\rho$
\begin{equation}\label{laplarhoah}
\laplaz \rho = \rho^{2} \laplahz \rho - (n-2) \rho \norm{d\rho}{\hz}^{2}.
\end{equation}
More generally,
\begin{equation}\label{laplaah}
\forall u \in \cinf_{c}(\m) \; \; , \; \; \laplaz u = \rho^{2} \left(\laplahz u - \frac{n-2}{\rho} \, d\rho \pshz du\right) ,
\end{equation}
where $\pshz$ is the scalar product w.r.t the metric $\hz$.\\
For tensor norms, for all $u \in T^{r}_{m}(\m) \; , \; \norm{u}{\gz} =  \rho^{m-r}\norm{u}{\hz}.$\\

\begin{defperso}{A cut-off function:}\label{cutoff}
\\$(\m, \gz)$ is a conformally compact manifold with defining function $\rho$.\\
Let $\chi: \R \rightarrow \R$ be a smooth cut-off function such that:
\begin{itemize}
\item[\textbullet] $\chi (\R) \subset \lbrack 0,1 \rbrack$
\item[\textbullet] supp $\chi \subset (- \infty, 2 \rbrack$
\item[\textbullet] $\restriction{\chi}{(- \infty, 1 \rbrack}=1$
\end{itemize}
then for $R$ large enough, we can define a cut-off function on $\m$ by
$$\chi_{R}(x) = \chi(-\ln(\rho(x))/R).$$
Setting $\Omega_{R} = \lbrace x \in \m : \rho(x) \,\corec{>}\, e^{-2R}\rbrace$ , then $\chi_{R}$ verifies
\begin{eqnarray*}
\chi_{R} &=& \begin{cases}
1 \; \corec{\text{on}} \; \; \Omega_{R \diagup 2} \\
0 \; \corec{\text{on}} \; \; \m \setminus \Omega_{R}
\end{cases}
\end{eqnarray*}
In other words, $\chi_{R}$ is a cut-off function near the boundary \corec{at infinity} of $\m$.\\
%De plus, pour $S$ assez grand and $R \gg S$, on peut définir une autre fonction de troncature sur $\m$ par $\chi_{S,R}(x) = (1-\chi_{S}(x))\chi_{R}(x)$. Le  support de $\chi_{S,R}, \; \Omega_{S,R} = \lbrace e^{-2R} < \rho < e^{-S} \rbrace$ est donc une couronne.
\end{defperso}

\subsection{Sobolev and Hölder weighted spaces.}$\mbox{ }$\\
Thanks to the background metric $\gz$ , we can define the following norm
\begin{eqnarray*}
\forall \; 1\infeg p < \infty \;\; \corec{\corec{\text{and}}} \; \; \delta\in \R &:& \; \; \Norm{u}{p,\delta} = \left( \int_{\m} \norm{u}{\gz}^{p} \; \rho^{p \delta} d\mu (\gz) \right)^{1 \slash p}\\
\text{For} \; p= \infty \; \; \corec{\corec{\text{and}}} \; \;  \delta\in \R &:& \; \; \Norm{u}{\infty,\delta} = \underset{\m}{sup} \left( \rho^{\delta} \norm{u}{\gz} \right)
\end{eqnarray*}

The Lebesgue weighted space $\leb{p}{\delta}$ is now defined as the space of measurable functions of \corec{$L^{p}_{\text{loc}}$} whose norm mentioned above is finite. The Sobolev weighted space $\sob{k}{p}{\delta}$  is then the space of measurable functions of $W^{k,p}_{\text{loc}}$ whose following norm is finite:

$$\Norm{u}{k,p,\delta} = \sum_{|\alpha| \infeg k} \Norm{\nz^{\alpha}u}{p,\delta} ,$$
where $\alpha$ is a multi-index of size $n$ and $\nz^{\alpha}u = \nz_{i_{1}}^{\alpha_{1}} \ldots \nz_{i_{n}}^{\alpha_{n}}u$ ,
$$\alpha = (\alpha_{1}, \ldots , \alpha_{n}) \; \; \corec{\corec{\text{and}}} \; \; |\alpha| = \sum_{i = 1}^{n} \alpha_{i} .$$
NB: for $\delta=0$, we get back to norms of the classic Lebesgue and Sobolev spaces.\\

 $\sob{k}{p}{\delta}(T^{r}_{m} \m)$ 
%(resp. $\sob{k}{p}{\delta}(T^{r}_{m} \m , \hz, \rho)$ )
  will refer to Sobolev spaces of sections of the $(r,m)$-tensor bundle over $\m$.
 %relativement à $\gz$ (resp. $\hz$). Quand il n'y a pas lieu de confondre, on omettra de préciser la métrique afin d'alléger les notations
We will speak of $\sob{k}{p}{\delta}$- norm to indicate the norm of a \corec{tensor field} $u \in \sob{k}{p}{\delta}(T^{r}_{m} \m )$. For a domain $U\subset \m$, $\Norm{u}{p, \delta ; U}$ will be the restriction to $U$ of the $\sob{k}{p}{\delta}$- norm of $u$. The Hölder weighted space $C^{s, \alpha}_{\delta} (\m,g)$ , with $0< \alpha<1$ is endowed with the norm
$$\Norm{u}{C^{s, \alpha}_{\delta}} = \underset{|k| \infeg s}{\text{max}} \Norm{\nz^{k}u}{C^{0, \alpha}_{\delta}}$$
with $$\Norm{u}{C^{0, \alpha}_{\delta}} = \underset{x \in \m}{\text{sup}} \rho^{\delta} \norm{u}{\gz}  + \underset{x \in \m}{\text{sup}} \rho^{\delta} \left( \underset{d_{\gz}(x,y) \infeg 1}{\text{sup}} \frac{\norm{\w{u}(x) - \w{u}(y)}{\w{g}}}{d_{\gz}(x,y)^{\alpha}} \right) ,$$
where $\w{u}$ and $\w{g}$ represent tensors $u$ and $g$ in an appropriate orthonormal basis.

\subsection{Elliptic operators.}$\mbox{ }$\\
Here we recall classic results on elliptic operators that we may find in \cite{Andersson1993} for example.
Let $B_{1}$ and $B_{2}$ be two tensor bundles over a conformally compact manifold $(\m,\corec{\gz})$ with defining function $\rho$ and $A: \cinf(B_{1}) \rightarrow \cinf(B_{2})$ be a partial differential linear operator of order $m$ define\corec{d} by
\begin{equation}\label{defAop}
A= \sum_{|\alpha| \infeg m} a_{\alpha}\nz^{\alpha}
\end{equation}
Set $s \in \nat$. We say \corec{that the} operator $A$ of the form (\ref{defAop}) has symbol in $\mathcal{OP}^{m}_{s}$ if
$$a_{\alpha} \in C^{s_{\alpha}}_{-|\alpha|}L(B_{1}, B_{2}) , \; \; \; \text{with} \; \; s_{\alpha} = \text{max} (s,|\alpha|-m+1) .$$
We say \corec{that} $A$ is an elliptic operator if
\corec{\begin{itemize}
\item[\textbullet] For all $\alpha$ such that $|\alpha|=m$ , for all
$\xi^{\alpha} = \xi_{1}^{\alpha_{1}} \ldots \xi_{n}^{\alpha_{n}} \ne 0$,\\
 $a_{\alpha} \xi^{\alpha} : B_{1} \rightarrow B_{2}$ is a tensor bundles isomorphism.
\item[\textbullet] For all $\xi^{\alpha}$, there exists two constants $c_{1}$ and $c_{2}$ such that
$$\Norm{a_{\alpha} \xi^{\alpha}}{\corec{\gz}} < c_{1} \; \norm{\xi^{\alpha}}{\corec{\gz}} \hspace{.8cm} \text{and} \hspace{.8cm} \Norm{(a_{\alpha} \xi^{\alpha})^{-1}}{\corec{\gz}} < c_{2} \; \norm{\xi^{\alpha}}{\corec{\gz}} .$$
\end{itemize}
}

\begin{lemperso}\label{lemestellip}%(Lemme 2.4 %\cite{Andersson1993})
Set $s_{0} \in \nat$. For \corec{every} elliptic operator $A$ with symbol in $\mathcal{OP}^{m}_{s_{0}} \;$, there exists  a positive constant $c = c \, (\gz , \delta)$ such that the following inequality is valid for all $s \infeg s_{0}$ :
\begin{equation}\label{sbelemme}
\Norm{u}{2,s+m,\delta} \infeg c \, \left(\Norm{Au}{2,s,\delta} + \Norm{u}{2,0,\delta}\right).
\end{equation}
\end{lemperso}

\begin{thmperso}\label{estellip}%(prop 2.6 %\cite{Andersson1993})
For all $R >0$, let $\Omega_{R}$ be \corec{as} in Definition \ref{cutoff}. Given $\delta \in \R$ and $A$ an elliptic operator with symbol in $\mathcal{OP}^{m}_{0}$.
Suppose there exists $R$ large enough so that $A$ verifies
\begin{equation}
\forall u \in \cinf_{c}(\m \setminus \corec{\overline\Omega_{R}}) \; \; , \; \; \Norm{u}{2,\delta} \infeg C \, \Norm{Au}{2,\delta},
\end{equation}
where $C$ depends on $A$.\\
Then we can choose $R$ large enough so that the following inequality
\begin{equation}\label{sbe}
\Norm{u}{m,2,\delta} \infeg C \, \left(\Norm{Au}{2,\delta} + \Norm{u}{2,\delta;\Omega_{R}}\right)
\end{equation}
is valid for all $u \in \cinf_{c}(\m)$. In particular, $A: \sobd{m} \rightarrow \lebd$ is semi-Fredholm, $\ie A$ has finite dimensional kernel and closed range.
\end{thmperso}

\section{Preliminary analysis}

In this section, we introduce some useful inequalities (see  \cite{Andersson1993} , th 2.3 for example):

\begin{propperso}{Weighted Hölder inequalities (in any dimension):}
\begin{list}{$\bullet$}{}
\item
Set $\delta \in \R$ so that $\delta = \delta_{1} + \delta_{2}$ , let $p,q,r \in \nat$ be such that
$$ 1\infeg p\infeg q \infeg r \infeg \infty \;\; \text{and} \;\; \frac{1}{p} = \frac{1}{q} + \frac{1}{r} \, ,$$ then
\begin{equation}\label{h1}
\Norm{uv}{p ,\delta} \infeg \Norm{u}{q ,\delta_{1}} \; \Norm{v}{r ,\delta_{2}}.
\end{equation}
\item Set $\delta \in \R$ and $\lambda \in \lbrack 0,1 \rbrack$ , let $p,q,r \in \nat$ be such that
$$ 1\infeg p\infeg q \infeg r \infeg \infty \;\; \corec{\corec{\text{and}}} \;\; \frac{1}{p} = \frac{\lambda}{q} + \frac{1- \lambda}{r} \, ,$$ then
\begin{equation}\label{h2}
\Norm{u}{p ,\delta} \infeg \Norm{u}{q ,\delta}^{\lambda} \; \Norm{u}{r ,\delta}^{1- \lambda}.
\end{equation}
\end{list}
\end{propperso}

\begin{thmperso}{}
Let $(\m,g)$ be a conformally compact $n$-dimensional manifold.\\[.2cm]
Weighted Sobolev inclusion:\\
For all $ 1\infeg p\infeg q< \infty $ , for all $k\supeg k'$, \; if $\delta \infeg \delta ' - \frac{n}{p} \frac{(q-p)}{q}$, \; then $\sob{k}{q}{\delta} \subset \sob{k'}{p}{\delta '}$ ,\\
\corec{and} there exists a positive constant $c = c \, (\gz , \delta , \delta ' , k, k',n, p,q)$ such that
$$\Norm{u}{k', p ,\delta'} \infeg c \, \Norm{u}{k, q ,\delta}.$$
If $ q= \infty$, for all $ 1\infeg p\infeg \infty $; for all $ k \supeg k'$ , \; if $\delta \infeg \delta ' - \frac{n}{p}$, \; then $\sob{k}{\infty}{\delta} \subset \sob{k'}{p}{\delta '}$ ,\\
\corec{and} there exists a positive constant $c = c \, (\gz , \delta , \delta ' , k, k',n,p)$ such that
$$\Norm{u}{k', p ,\delta'} \infeg c \, \Norm{u}{k, \infty ,\delta}.$$
Weighted Hölder inclusion:\\
Set $k \in \nat$, $n=3$ and $p=2$, \; then for all $ 0 < \alpha \infeg \demi$, for all $\delta \infeg \delta'$ , there exists a positive constant $c = c \, (\gz , \delta , \delta ', k, k' , \alpha)$ such that
\begin{equation}\label{inchold}
\forall u \in \sobd{2}(\m) \; \; , \; \;\Norm{u}{C^{0, \alpha}_{\delta '}} \infeg c \, \Norm{u}{2, 2 ,\delta}.
\end{equation}
Weighted Sobolev inequalities:\\
Set $1\infeg p < \infty$ and let $k,j$ be integers. \corec{ In each of the following cases, there exists a positive constant $c = c \,(\gz , \delta , k,n,j,p,q)$ such that
for all $u \in \sob{j+k}{p}{\delta}(\m)$,}
\begin{list}{$\bullet$}{}
\item If $pk<n, \; \Norm{u}{j, q ,\delta} \infeg c \, \Norm{u}{j+k, p ,\delta} \; , \; \; \forall \, p\infeg q \infeg \frac{np}{n-kp}$.\\
\item If $pk=n, \; \Norm{u}{j, q ,\delta} \infeg c \, \Norm{u}{j+k, p ,\delta} \; , \; \; \forall \, p\infeg q < \infty$.\\
\item If $pk>n, \; \Norm{u}{j, q ,\delta} \infeg c \, \Norm{u}{j+k, p ,\delta} \; , \; \; \forall \, p\infeg q \infeg \infty$.
\end{list}
\medskip
Ehrling inequality:\\
For all $\epsilon >0$ , for all integers $j,k$ such that $0 <j<k$ , there exists a positive constant $C(\epsilon)$ such that
\begin{equation}\label{Ehrling}
\forall u \in \sob{k}{p}{\delta} \; \; , \; \; \Norm{u}{j, p ,\delta} \infeg \epsilon  \Norm{u}{k, p ,\delta} + C(\epsilon)  \Norm{u}{p ,\delta}.
\end{equation}
Rellich Theorem:\\
For all $k>k'$ and $\delta < \delta'$ , the inclusion $\sobd{k} \subset \sob{k'}{2}{\delta'}$ is compact.
\end{thmperso}
A consequence of the Sobolev inequalities (cf. \cite{Chrusciel2004} for example) is
\begin{propperso}{}
In dimension $n=3$, for all $k > \tfrac{3}{2}$ , 
\begin{equation}\label{tauxdec}
u \in \sobd{k} \Rightarrow u = o(\rho^{- \delta}).
\end{equation}
\end{propperso}

\begin{lemperso}{}
In dimension 3, for all $\delta \in \R$, there exists a positive constant $c = c \, (\gz, \delta_{1} , \delta_{2})$ such that
\begin{equation}\label{19g}
\Norm{uv}{2 ,\delta} \infeg c \; \Norm{u}{1,2,\delta_{1}} \; \Norm{v}{1,2 ,\delta_{2}} \; \; , \; \; \text{with} \; \delta_{1}+\delta_{2} = \delta.
\end{equation}
\end{lemperso}
Remark: In the particular case where $\delta_{1}, \delta_{2}$ and $\delta$ are non positive, then $\delta_{1}$ and $\delta_{2}$ are both \corec{greater than} $\delta$. The weighted Sobolev inclusion leads to
\begin{equation}\label{19gbis}
\Norm{uv}{2 ,\delta} \infeg c \; \Norm{u}{1,2,\delta} \; \Norm{v}{1,2 ,\delta}. % \; \; , \; \; \forall \; \delta \infeg \delta_{1},\delta_{2} \infeg 0.
\end{equation}

\begin{lemperso}{}
For each of the following inequalities, in dimension $3$, there exists a positive constant $c = c \, (\gz, \delta)$ such that
\begin{eqnarray}
\forall u \in \sob{2}{2}{\delta}(\m) \; \; , \; \; \Norm{u}{\infty ,\delta} &\infeg& \epsilon \; \Norm{u}{2,2,\delta} +c \epsilon^{-3}\Norm{u}{1,2 ,\delta}. \label{uinfd}\\
\forall u \in \sob{1}{2}{\delta}(\m) \; \; , \; \; \Norm{u}{3 ,\delta} &\infeg& \epsilon \; \Norm{u}{1,2,\delta} + c \epsilon^{-1}\Norm{u}{2 ,\delta}. \label{u3d}
\end{eqnarray}
\end{lemperso}

\begin{lemperso}{}
In all dimension $n$, we define for $R$ real large enough, $$E_{R}:= \m \setminus \Omega_{R} = \lbrace \rho \infeg e^{-2R} \rbrace .$$ Then $\forall u \in \leb{p}{\delta}(E_{R}) , \; \forall \; \delta \in \R$,
\begin{equation}\label{49bis}
\Norm{u}{p,\eta;E_{R}} \infeg e^{2R (\delta - \eta)} \Norm{u}{p,\delta;E_{R}} \; \; \; \; ,  \forall \, \delta \infeg \eta \, , \, \forall \, 1 \infeg p \infeg \infty.
\end{equation}
\end{lemperso}

\section{The Hessian type operator $\Tz$}

%Dans le but d'utiliser le théorème des fonctions implicites, on s'intéresse au noyau de $D \phiogs$ and on commence par montrer sa coercivité dans la Proposition \ref{prop3.3}. 
Here we give some preliminary results concerning the operator $\Tz$ defined for a function $N$ by: 
\begin{equation}\label{defTz}
\Tz =\Tz(N) := \nz^{2}N - N \gz.
\end{equation}
\begin{lemperso}{}\label{lemtz}
For all $\delta \in \R$, there exists a positive constant $c>0$ depending on $\gz$ such that for all $N \in \sob{2}{2}{-\delta}(\m)$,
\begin{equation}
\Norm{\Tz(N)}{2, -\delta} \supeg \Norm{\nz^{2}N}{2,-\delta} - c\, \Norm{N}{2,-\delta}.
\end{equation}
\end{lemperso}
\textbf{Proof}: The result stems from the definition of $\Tz$ and the Triangle inequality. $\cqfd$ \\

\begin{lemperso}{}\label{lemn1inftz}
For all $\delta \in \rbrack -(n+1) \slash 2 \, , 0 \rbrack$, there exists a positive constant $c = c\,(\gz , \delta)$  such that for all $N \in \sob{2}{2}{-\delta}(\m)$,
\begin{equation}\label{lemn1inftzbis}
\Norm{N}{1,2,-\delta} \infeg c \, \Norm{\Tz(N)}{2, -\delta}.
\end{equation}
\end{lemperso}
\textbf{Proof}: By density, we can suppose $N \in \cinf_{c}(\m)$.
We use the {proof} of Proposition \ref{n2inftzinf} which establishes (\ref{lemn1inftzbis}) if the support of $N$ is in a neighborhood of the boundary at infinity. Here we ignore the $\delta$ restriction due to positivity of the interior boundary term since $N$ is compactly supported. We obtain the result near the boundary for $\delta \in \rbrack -(n+1) \slash 2 \, ; 0 \rbrack$ and conclude with kernel triviality of $\Tz$ for \corec{$-\delta < \frac{n+1}{2}$} (see \cite{ChruscielDelayExotic}) thanks to a {proof} similar to the one of Theorem \ref{estellip}. $\cqfd$ \\[.2cm]

Combination of Lemmas \ref{lemtz} and \ref{lemn1inftz} gives

\begin{propperso}{}\label{propn2infTz}
For all $\delta \in \rbrack -(n+1) \slash 2 \, , 0 \rbrack$ , there exists a positive constant $c = c \, (\gz, \delta)$ such that
\begin{equation}\label{n2infctz}
\Norm{N}{2,2,-\delta} \infeg c \, \Norm{\Tz(N)}{2,-\delta}.
\end{equation}
\end{propperso}
%La proposition suivante est l'analogue en asymptotiquement hyperbolique %de la Proposition 3.3 de \cite{Bartnik2005}.

We will need the next three lemmas which give general equalities on  $(\m, \gz)$, a $n$-dimensional asymptotically hyperbolic Riemannian manifold with $\gz = \rho^{-2} \hz$.
\corec{Here we allow an possible inner boundary $\d\m$}. From now on, $d\sigma(\gz)$ will be the measure induced by $\gz$ on $\d\m$ and $\eta$ is the exterior unit normal to $\d\m$.
\corec{The term $o(1)$ will  tend to zero when approaching $\d_\infty\m$.}
\begin{lemperso}{}\label{lem2}
Let $(M, \gz)$ be a $n$-dimensional asymptotically hyperbolic manifold and\\ $N \in \cinf(\m)$ with a compact support \corec{on $\m$}. $\forall \delta \in \R$ ,
\begin{equation*}
\int_{\m} 2N \langle dN , \frac{d\rho}{\rho} \rangle_{\gz} \, \rho^{2\delta} \, d\mu(\gz) = -\int_{\m} \lbrack 2 \delta +1 - n +o(1) \rbrack N^{2} \rho^{2\delta} \, d\mu(\gz) \, + \int_{\d\m} N^{2}\langle \frac{d\rho}{\rho} , \eta \rangle_{\gz} \, \rho^{2\delta} \, d\sigma(\gz) \, .
\end{equation*}
\end{lemperso}
\textbf{Proof}: Integration by parts gives
\begin{eqnarray}\label{2bis}
\int_{\m} \nz_{i}(-N^{2}\nz^{i}(\rho^{-1})\rho^{2\delta+1}) \, d\mu(\gz)& =& \int_{\m} 2N \langle \,  dN , \frac{d\rho}{\rho} \rangle_{\gz} \, \rho^{2\delta} \, d\mu(\gz) - \int_{\m} N^{2} \laplaz (\rho^{-1})\rho^{2\delta+1} \, d\mu(\gz) \nonumber\\
&&+ \int_{\m} (2 \delta +1) N^{2} \norm{d\rho}{\hz}^{2}\rho^{2\delta} \, d\mu(\gz).
\end{eqnarray}
Let us compute (see \cite{Chrusciel2004}, (D.4) for instance)
\begin{equation}\label{hesszrho-1}
\nz^{2}(\rho^{-1}) = \rho^{-1} \norm{d\rho}{\hz}^{2} \gz - \rho^{-2} \precig{\hz} \n^{2}\rho.
\end{equation}
Taking the $\gz$-trace,
\begin{equation*}
\laplaz (\rho^{-1}) = n \rho^{-1} \norm{d\rho}{\hz}^{2} - \laplahz \rho.
\end{equation*}
The metric $\hz$ being defined (and so bounded) \corec{until $\d_\infty\m$} and $\rho$ being a smooth function on \corec{$\bar\m$} , $\laplahz \rho$ is a smooth function bounded on \corec{$\bar\m$} and so, we can write $\laplahz \rho =  \, O(1)=o(\rho^{-1})$ near \corec{$\d_\infty\m$}. We obtain
\begin{equation}\label{laplazrho-1}
\laplaz (\rho^{-1}) = \rho^{-1} \big(n\norm{d\rho}{\hz}^{2} + o(1)\big).
\end{equation}
According to $ \norm{d\rho}{\hz}^{2} = 1 + o(1)$ near the boundary \corec{at infinity} on an asymptotically hyperbolic manifold, (\ref{2bis}) become
\begin{eqnarray}\label{2ter}
\int_{\m} \nz_{i}(-N^{2}\nz^{i}(\rho^{-1})\rho^{2\delta+1}) \, d\mu(\gz)& =& \int_{\m} 2N \langle dN , \frac{d\rho}{\rho} \rangle_{\gz} \, \rho^{2\delta} \, d\mu(\gz) \nonumber\\
&& + \int_{\m} \lbrack 2 \delta +1 - n +o(1) \rbrack N^{2} \rho^{2\delta} \, d\mu(\gz).
\end{eqnarray}

From the Divergence theorem,
\begin{equation}\label{tdflem2}
\int_{\m} \nz_{i}(-N^{2}\nz^{i}(\rho^{-1})\rho^{2\delta+1}) \, d\mu(\gz) = \int_{\d\m} N^{2}\langle \frac{d\rho}{\rho} , \eta \rangle_{\gz} \, \rho^{2\delta} \, d\sigma(\gz).
\end{equation}
We end the proof replacing the left-hand side of (\ref{2ter}) by its expression (\ref{tdflem2})$.\cqfd$\\

%The two following lemmas are obtained thanks to integration by parts, the Divergence theorem and Lemma \ref{lem2}: \erw{donner des preuves courtes des lemmes : dire ce qu'on integre comme divergence}
\begin{lemperso}{}\label{lem4}
Let $(M, \gz)$ be a $n$-dimensional asymptotically hyperbolic manifold and\\ $N \in \cinf(\m)$ with a compact support \corec{ on $\m$}. $\forall \delta \in \R$ ,
\begin{eqnarray}\label{4}
-2\int_{\m} 
\begin{comment}(\nz^{2}N - N\gz)\end{comment}
 \Tz  (dN , \frac{d\rho}{\rho}) \rho^{2\delta} \, d\mu(\gz) &=& \int_{\m} \lbrace 2 \delta +1 - n +o(1) \rbrace \norm{dN}{\gz}^{2} \rho^{2\delta} \, d\mu(\gz) \nonumber\\
&& - \int_{\m} \lbrack 2 \delta +1 - n +o(1) \rbrack N^{2} \rho^{2\delta} \, d\mu(\gz) \nonumber\\
&& + \int_{\d\m} (N^{2}-\norm{dN}{\gz}^{2})\langle \frac{d\rho}{\rho} , \eta \rangle_{\gz} \, \rho^{2\delta} \, d\sigma(\gz) ,
\end{eqnarray}
\end{lemperso}

\textbf{Proof}: We integrate by parts the term 
%\begin{eqnarray}\label{4bis}
%\int_{\m} 
$\nz_{i}(\norm{dN}{\gz}^{2} \, \nz^{i}(\rho^{-1}) \, \rho^{2\delta+1})$ and the result follows on from the Divergence theorem and Lemma \ref{lem2}. $\cqfd$

\begin{lemperso}{}\label{lem6}
Let $(M, \gz)$ be a $n$-dimensional asymptotically hyperbolic manifold and\\ $N \in \cinf(\m)$ with a compact support \corec{ on $\m$}. $\forall \delta \in \R$ ,
\begin{eqnarray}\label{6}
-\int_{\m} 
\begin{comment}(\laplaz N - n N)\end{comment}
\corec{\tr}_{\gz}\Tz \, N \, \rho^{2\delta} \, d\mu(\gz)& =& - \int_{\m} \lbrack \,  \delta (2 \delta +1 - n) - n + o(1) \rbrack N^{2} \rho^{2\delta} \, d\mu(\gz) \nonumber\\
&& +\int_{\m} \norm{dN}{\gz}^{2}\rho^{2\delta} \, d\mu(\gz) - \int_{\d\m} N \langle dN , \eta \rangle_{\gz} \, \rho^{2\delta} \, d\sigma(\gz) \nonumber\\
&& +\int_{\d\m} \delta N^{2} \langle dN , \eta \rangle_{\gz} \, \rho^{2\delta} \, d\sigma(\gz) ,
\end{eqnarray}
\end{lemperso}
%\begin{comment}
\textbf{Proof}: We integrate by parts the term 
%\begin{eqnarray}\label{4bis}
%\int_{\m} 
$\nz_{i}(N\nz^{i}N \rho^{2\delta})$ and the result follows on from the Divergence theorem and Lemma \ref{lem2}. $\cqfd$ \\[.3cm]

The next proposition stems from the two previous lemmas and will play an important role to prove the adjoint kernel triviality:
\begin{propperso}{}\label{n2inftzinf}
For every $ \epsilon >0$ , for all $\delta \in \rbrack -(n+1) \slash 2 , -1 \lbrack$, there exists $ R_{\epsilon,\delta}>0 $ such that for all $R>R_{\epsilon,\delta}$, there exists a positive constant $c$ such that
\begin{equation}\label{n2infctzbord}
\forall N \in \cinf_{c}(E_{R})\;,\;\;\Norm{N}{2,2,-\delta; E_{R}} \infeg c \, \Norm{\Tz}{2,-\delta; E_{R}}.
\end{equation}
\end{propperso}
\textbf{Proof}: From $\Tz$ and $\corec{\tr}_{\gz}\Tz$ expressions and  $\nz_{n}N$ (resp. $\nz_{T}N$) being the component of $d N$ normal (resp. tangential) to $\corec{\d_{\infty}\m}$, $\nz_{n}N := \langle dN , \eta \rangle_{\gz} \; \; \; \text{and} \; \; \; \norm{dN}{\gz}^{2} = \norm{ \nz_{n}N}{\gz}^{2}+ \norm{\nz_{T}N}{\gz}^{2}$.\\ 
$(\ref{6}) - \tdemi (\ref{4})$ give
\begin{eqnarray}\label{7}
&&\int_{\m} \Tz (dN , \frac{d\rho}{\rho}) \rho^{2\delta} \, d\mu(\gz) -\int_{\m}N \corec{\tr}_{\gz}\Tz \rho^{2\delta} \, d\mu(\gz) \nonumber\\
 &=& \int_{\m} \lbrace \tfrac{n+1}{2} - \delta +o(1) \rbrace \norm{dN}{\gz}^{2} \rho^{2\delta} \, d\mu(\gz) + \int_{\m} \lbrack -2 \delta^{2} + n \delta + \tfrac{n+1}{2} +o(1) \rbrack N^{2} \rho^{2\delta} \, d\mu(\gz) \nonumber\\
&& + \int_{\d\m} \lbrace (\delta - \tdemi)N^{2} + \tdemi \norm{dN}{\gz}^{2} \rbrace \langle \frac{d\rho}{\rho} , \eta \rangle_{\gz} \rho^{2\delta} \, d\sigma(\gz) - \int_{\d\m} N \nz_{n}N \, \rho^{2\delta} \, d\sigma(\gz).
\end{eqnarray}
Application on $E_{R}$:
$E_{R}$ possesses two disjoint boundary components. A boundary at infinity, noted \corec{$\d E_{\infty}=\d_\infty \m$}, and an inner  boundary $\d \Omega_{R} = \{\rho = \e^{-2R} \}$. Since $N \in \cinf_{c}(E_{R})$, $N$ is surely null near $\d E_{\infty}$ but not necessarily on $\d \Omega_{R}$ and that's the reason why boundary terms in (\ref{7}) will only concern $\d \Omega_{R}$. If $\eta_{R}$  is the normal to $\d \Omega_{R}$ exterior to $E_{R}$ and considering that when $R \rightarrow +\infty \; , \; \eta_{R} - \tfrac{d\rho}{\rho} \rightarrow 0$ , so that $\langle \tfrac{d\rho}{\rho} , \eta_{R} \rangle_{\gz} = \tfrac{\norm{d\rho}{\gz}^{2}}{\rho^{2}} + o(1) = 1 + o(1)$.\\
\begin{eqnarray*}
&&\int_{E_{R}} \Tz (dN , \frac{d\rho}{\rho}) \rho^{2\delta} \, d\mu(\gz) -\int_{E_{R}}N\, \corec{\tr}_{\gz}\Tz \rho^{2\delta} \, d\mu(\gz) \nonumber\\
 &=& \int_{E_{R}} \lbrace \tfrac{n+1}{2} - \delta +o(1) \rbrace \norm{dN}{\gz}^{2} \rho^{2\delta} \, d\mu(\gz) + \int_{E_{R}} \lbrack -2 \delta^{2} + n \delta + \tfrac{n+1}{2} +o(1) \rbrack N^{2} \rho^{2\delta} \, d\mu(\gz) \nonumber\\
&& + \int_{\d E_{R}} \lbrace (\delta - \tdemi + o(1))N^{2} + \tdemi \norm{dN}{\gz}^{2} - N \nz_{n}N + o(1) \rbrace \rho^{2\delta} \, d\sigma(\gz).
\end{eqnarray*}
According to the following inequalities
$$
\begin{cases}
 \norm{\corec{\tr}_{\gz}\Tz}{\gz}^{2} \infeg n \, \norm{\Tz}{\gz}^{2}.\\
\Tz (dN , \frac{d\rho}{\rho}) \infeg \tfrac{a}{2}\norm{\Tz}{\gz}^{2} + \frac{1}{2a} \norm{dN}{\gz}^{2} \norm{d\rho}{\hz}^{2} \; \; \; \; \; , \forall a>0 \, , \; a >>1. \\
- N \,\corec{\tr}_{\gz}\Tz \infeg \tfrac{b}{2}\,\norm{\corec{\tr}_{\gz}\Tz}{\gz}^{2} + \frac{1}{2b} N^{2} \norm{d\rho}{\hz}^{2} \; \; \; \; \; , \forall b>0 \, , \; b >>1.
\end{cases}
$$
$\forall \epsilon >0 , \exists R_{\epsilon}>0 $ such that $ \forall R>R_{\epsilon}$ ,
\begin{eqnarray*}
\Big(\frac{a + bn}{2}\Big)\int_{E_{R}} \norm{\Tz}{\gz}^{2} \rho^{2\delta} \, d\mu(\gz) &\supeg& \int_{E_{R}} \lbrace \tfrac{n+1}{2} - \delta - \epsilon \rbrace \norm{dN}{\gz}^{2} \rho^{2\delta} \, d\mu(\gz)\\
&& + \int_{E_{R}} \lbrack -2 \delta^{2} + n \delta + \tfrac{n+1}{2} - \epsilon \rbrack N^{2} \rho^{2\delta} \, d\mu(\gz) \nonumber\\
&& + \int_{\d E_{R}} (\delta - 1 - \epsilon)N^{2} \rho^{2\delta} \, d\sigma(\gz)\\
&& + \int_{\d E_{R}} \left\lbrace \tdemi \norm{\nz_{T}N}{\gz}^{2} + \tdemi (N - \nz_{n}N)^{2} \right\rbrace \rho^{2\delta} \, d\sigma(\gz).
\end{eqnarray*}
\begin{itemize}
\item[\textbullet] The $N^{2}$ interior term
%, nous avons un trinôme du second degré en $\delta$ dont le discriminant est
%$$\Delta = n^{2} + 4(n+1) = (n+2)^{2}>0$$
%et dont les racines sont
%$$\delta_{\pm} = \frac{n \pm (n+2)}{4} = \Big \lbrace -\demi \, ; \frac{n+1}{2} \Big \rbrace.$$
%Ce terme
 is non negative if $\delta \in \rbrack -\demi \, ; \frac{n+1}{2} \lbrack$.\\
\item[\textbullet] The $\norm{dN}{\gz}^{2}$ interior term is non negative if $\delta < (n+1) \slash 2$.\\
\item[\textbullet] The boundary term is non negative if $\delta > 1$.
\end{itemize}
Moreover, a quick calculation shows that on the interval $\lbrack 0 \, ; \frac{n+1}{2} \lbrack$ ,
$$\tfrac{n+1}{2} - \delta \infeg -2 \delta^{2} + n \delta + \tfrac{n+1}{2}.$$
Consequently, for $\delta \in \rbrack 1 \, ; \frac{n+1}{2} \lbrack$, \\
$\forall \epsilon >0 , \exists R_{\epsilon}>0 $ such that $ \forall R>R_{\epsilon}$ ,
\begin{eqnarray*}
\Big(\frac{a + bn}{2}\Big)\int_{E_{R}} \norm{\Tz}{\gz}^{2} \rho^{2\delta} \, d\mu(\gz) &\supeg& \int_{E_{R}} \lbrace \tfrac{n+1}{2} - \delta - \epsilon \rbrace \Big(N^{2} + \norm{dN}{\gz}^{2}\Big) \rho^{2\delta} \, d\mu(\gz) .
\end{eqnarray*}
Namely
$$ \Norm{N}{1,2,\delta; E_{R}} \infeg c \, \Norm{\Tz}{2,\delta; E_{R}}. $$
Combining this inequality with Lemma \ref{lemtz} valid in particular on $E_{R}$, we end the proof of (\ref{n2infctzbord}).$\cqfd$

\section{\corec{The} Killing operator $\Sz$}

Let $\Sz$ be the Killing operator defined on 1-forms by
\begin{equation}
\Sz(Y)_{ij} = \demi (\nz_{i}Y_{j} + \nz_{j}Y_{i}) = \nz_{(i}Y_{j)}.
\end{equation}
The trace of this operator is
\begin{equation}
\corec{\tr}_{\gz} \Sz(Y) = \gz^{ij}\Sz(Y)_{ij} = \nz^{i}Y_{i} =: \corec{\div }Y.
\end{equation}

The next three lemmas are respectively versions of Lemma D.1  and Propositions D.2 and D.3 from \cite{Chrusciel2004} with \corec{inner} boundary:
\begin{lemperso}{}\label{lemd1}
Let $V$ be a vector field and $Y$ a $1-$form both compactly supported on $\m$. Then,
\begin{eqnarray*}
&&\int_{\m} (\mathring{S}(Y) + \tdemi \corec{\tr}_{\gz}(\mathring{S}(Y))\gz)(Y,V) \, d\mu(\gz)\\
&&=  - \demi  \int_{\m} \big\{ \nz V (Y,Y) + \tdemi\corec{ \corec{\div }}V \norm{Y}{\gz}^{2} \big\} \, d\mu(\gz)
 + \demi \int_{\d\m} \langle Y , V \rangle_{\gz} \; \langle Y , \eta \rangle_{\gz} \, d\sigma(\gz)\\
&& + \frac{1}{4} \int_{\d\m} \norm{Y}{\gz}^{2} \; \langle V , \eta \rangle_{\gz} \, d\sigma(\gz) ,
\end{eqnarray*}
\end{lemperso}

\begin{lemperso}{}\label{propd2}
Let $u$ be a function, $V$ a vector field and $Y$ a $1-$form all compactly supported on $\m$. Then,
\begin{eqnarray*}
&&\int_{\m} \e^{2u}(\mathring{S}(Y) + \tdemi \corec{\tr}_{\gz}(\mathring{S}(Y))\gz)(Y,V) \, d\mu(\gz)\\
&& =  - \demi  \int_{\m} \e^{2u} \Big\{ \nz V (Y,Y) + \tdemi\corec{ \corec{\div }}V \norm{Y}{\gz}^{2} \Big\} \, d\mu(\gz) + \demi \int_{\d\m} \e^{2u} \langle Y , V \rangle_{\gz} \; \langle Y , \eta \rangle_{\gz} \, d\sigma(\gz)\\
&&- \demi \int_{\m} \e^{2u} \Big\{2 \langle du , Y \rangle_{\gz} \; \langle Y , V \rangle_{\gz} +  \langle du , V \rangle_{\gz} \norm{Y}{\gz}^{2}\Big\} \, d\mu(\gz) +\frac{1}{4} \int_{\d\m} \e^{2u} \norm{Y}{\gz}^{2} \; \langle V , \eta \rangle_{\gz} \, d\sigma(\gz) ,
\end{eqnarray*}
\end{lemperso}

\begin{lemperso}{}\label{lempd3}
Let $Y$ be a $1-$form compactly supported $\m$ and $u,v \in \cinf(\m)$ two functions defined in a \corec{neighborhood of the support of $Y$}. Then,
\begin{eqnarray*}\label{pc2}
&&-2 \int_{\m} v \, \e^{2u} \mathring{S}(Y) (dv,dv) \langle dv , Y \rangle_{\gz}  \, d\mu(\gz) \\
&&= \int_{\m}  \e^{2u} \langle dv , Y \rangle_{\gz} \Big\{ \, \langle dv , Y \rangle_{\gz} \left[ \norm{dv}{\gz}^{2} + v \laplaz v + 2 v \langle dv , du \rangle_{\gz} \right] + 2 v \, \nz^{2}v (Y, dv) \, \Big\} \, d\mu(\gz) \nonumber\\
&& - \int_{\d\m} v \, \e^{2u} \langle dv , Y \rangle^{2}_{\gz} \langle dv , \eta \rangle_{\gz} \, d\sigma(\gz) ,
\end{eqnarray*}
\end{lemperso}

Respective versions of Corollaries D.4 and D.5 from \cite{Chrusciel2004} with \corec{possible inner} boundary come out from these previous lemmas:
\begin{corperso}{}\label{cord10}
Let $(M, \gz)$ be an asymptotically hyperbolic manifold and $Y$ a $1-$form\\ compactly supported on \corec{$\m$}. $\forall \delta \in \R$ ,
\begin{eqnarray}
&&2\int_{\m} \rho^{2\delta}(\mathring{S}(Y) + \tdemi \corec{\tr}_{\gz}(\mathring{S}(Y))\gz)(Y,\frac{d\rho}{\rho}) \, d\mu(\gz) \nonumber\\
&& = \int_{\m} \rho^{2\delta} \Big\{ (\tfrac{n+1}{2}- \delta + o(1)) \norm{Y}{\gz}^{2} -(2 \delta + 1) \langle \frac{d\rho}{\rho} , Y \rangle^{2}_{\gz} \Big\} \, d\mu(\gz) \nonumber\\
&& + \frac{1}{2} \int_{\d\m} \rho^{2\delta} \norm{Y}{\gz}^{2} \; \langle \frac{d\rho}{\rho} , \eta \rangle_{\gz} \, d\sigma(\gz) + \int_{\d\m} \rho^{2\delta} \langle Y , \frac{d\rho}{\rho} \rangle_{\gz} \; \langle Y , \eta \rangle_{\gz} \, d\sigma(\gz) , \label{eqcord10}
\end{eqnarray}
\end{corperso}
\textbf{Proof}: We apply lemma \ref{propd2} with
$$
\begin{cases}
V= d(\rho^{-1}) = - \rho^{-2} d\rho\\
u=(\delta + \tdemi) \ln \rho
\end{cases}
\; \; \text{and} \; \; 
\begin{cases}
\nz V= \nz^{2}(\rho^{-1}) = \rho^{-1} \norm{d\rho}{\hz}^{2} \gz - \rho^{-2} \precig{\hz} \n^{2}\rho  \; \; \; \mbox{from} \; \; (\ref{hesszrho-1})\\
\corec{\div }V = \laplaz (\rho^{-1}) = \rho^{-1} (n\norm{d\rho}{\hz}^{2} + o(1)) \; \; \; \mbox{from} \; \; (\ref{laplazrho-1})
\end{cases}
.$$
$$du=(\delta + \tdemi) \frac{d\rho}{\rho} \; \; \text{and} \; \; \e^{2u}= \rho^{2 \delta+1}.$$
The metric $\hz$ being defined (and so bounded) \corec{until $\d_{\infty}\m$} and $\rho$ being a smooth function on \corec{$\bar\m$} , $\precig{\hz} \n^{2}\rho$ is a smooth function bounded on \corec{$\bar\m$} and so, we can write $\precig{\hz} \n^{2}\rho =  \, o(\rho^{-1})$ near \corec{$\d_{\infty}\m$}. Thus,
$$
\begin{cases}
\nz V(Y,Y)  = \rho^{-1} (\norm{d\rho}{\hz}^{2} \norm{Y}{\gz}^{2} + o(1) \norm{Y}{\gz}^{2})\\
\corec{\div }V \norm{Y}{\gz}^{2} = \rho^{-1} (n \norm{d\rho}{\hz}^{2} + o(1))\norm{Y}{\gz}^{2}
\end{cases}
.$$
We obtain
\begin{eqnarray*}
&&\int_{\m} \rho^{2\delta}(\mathring{S}(Y) + \tdemi \corec{\tr}_{\gz}(\mathring{S}(Y))\gz)(Y,\frac{d\rho}{\rho}) \, d\mu(\gz)\\
&& = \demi  \int_{\m} \rho^{2\delta} \Big\{ \left[1 - (\delta + \tdemi) + \tfrac{n}{2} \right] \norm{d\rho}{\hz}^{2}\norm{Y}{\gz}^{2} -(2 \delta + 1) \langle \frac{d\rho}{\rho} , Y \rangle^{2}_{\gz} + o(1) \norm{Y}{\gz}^{2} \Big\} \, d\mu(\gz) \\
&& + \frac{1}{4} \int_{\d\m} \rho^{2\delta} \norm{Y}{\gz}^{2} \; \langle \frac{d\rho}{\rho} , \eta \rangle_{\gz} \, d\sigma(\gz) + \demi \int_{\d\m} \rho^{2\delta} \langle Y , \frac{d\rho}{\rho} \rangle_{\gz} \; \langle Y , \eta \rangle_{\gz} \, d\sigma(\gz).
\end{eqnarray*}
We end the proof with $ \norm{d\rho}{\hz}^{2} = 1 + o(1)$ on an asymptotically hyperbolic manifold. $\cqfd$

\begin{corperso}{}\label{cord11}
Let $(M, \gz)$ be an asymptotically hyperbolic manifold and $Y$ a $1-$form\\ compactly supported \corec{on $\m$}. Then,
\begin{eqnarray}
&& 2 \int_{\m} \rho^{2\delta}\mathring{S}(Y)(\frac{d\rho}{\rho},\frac{d\rho}{\rho}) \; \langle \frac{d\rho}{\rho} , Y \rangle_{\gz} \, d\mu(\gz) \nonumber\\
&& = \int_{\m} \rho^{2\delta}  \left(n-1- 2\delta + o(1) \right) \; \langle \frac{d\rho}{\rho} , Y \rangle^{2}_{\gz} \, d\mu(\gz) + \int_{\d\m} \rho^{2\delta} \; \langle \frac{d\rho}{\rho} , Y \rangle^{2}_{\gz} \langle \frac{d\rho}{\rho} , \eta \rangle_{\gz} \, d\sigma(\gz) ,\nonumber \\
&& \label{eqcord11}
\end{eqnarray}
\end{corperso}
\textbf{Proof}: We apply lemma \ref{lempd3} with
$$
\begin{cases}
v= \rho^{-1}\\
u=(\delta + 2) \ln \rho
\end{cases}
, \; \;
\begin{cases}
dv= d(\rho^{-1}) = - \rho^{-1}\frac{d\rho}{\rho} \\
du=(\delta + 2) \frac{d\rho}{\rho} 
\end{cases}
\; \; \text{and} \; \;
\begin{cases}
\e^{2u}= \rho^{2 \delta+4}\\
\norm{dv}{\gz}^{2} = \rho^{-2} \norm{d\rho}{\hz}^{2}
\end{cases}
$$
together with
$$
\begin{cases}
\nz^{2} v= \nz^{2}(\rho^{-1}) = \rho^{-1} \left( \norm{d\rho}{\hz}^{2} \gz - \rho^{-1} \precig{\hz} \n^{2}\rho \right) \; \; \mbox{from} \; \; (\ref{hesszrho-1})\\
\laplaz v = \laplaz (\rho^{-1}) = \rho^{-1} \left( n \norm{d\rho}{\hz}^{2} + o(1) \right) \; \; \mbox{from} \; \; (\ref{laplazrho-1})
\end{cases}
.$$
Since $ \norm{d\rho}{\hz}^{2} = 1 + o(1)$ on an asymptotically hyperbolic manifold, we end up with
\begin{eqnarray*}
&2 \int_{\m} \rho^{2\delta}\mathring{S}(Y)(\frac{d\rho}{\rho},\frac{d\rho}{\rho}) \; \langle \frac{d\rho}{\rho} , Y \rangle_{\gz} \, d\mu(\gz)
& = \int_{\m} \rho^{2\delta}  \Big\{ \, n + 3 - 2 (\delta +2)+ o(1) \Big\} \langle \frac{d\rho}{\rho}  , Y \rangle^{2}_{\gz} \, d\mu(\gz)\\
&& + \int_{\d\m} \rho^{2\delta} \; \langle \frac{d\rho}{\rho} , Y \rangle^{2}_{\gz} \langle \frac{d\rho}{\rho} , \eta \rangle_{\gz} \, d\sigma(\gz). \cqfd \nonumber
\end{eqnarray*}

The following lemma establishes a Korn-type inequality for the Killing operator $\Sz$:
\begin{lemperso}{}\label{lemyinfcs}
\corec{Assume $\m$ has no inner boundary.} Then for all $\delta > -(n+1) \slash 2$ and $ \delta \ne -(n-1) \slash 2$, there exists a positive constant $c = c \, (\gz,\delta)$ such that for all 1-form $Y \in \sob{1}{2}{-\delta}(T^{*}\m)$,
\begin{equation*}
\Norm{Y}{1,2,-\delta} \infeg c \, \Norm{\mathring{S}(Y)}{2, -\delta}.
\end{equation*}
\end{lemperso}
\textbf{Proof}: Here again we base ourselves on the {proof} of Theorem \ref{estellip} but for the operator $\mathring{S}$.
Lemma $2.8$ from \cite{Chrusciel2004} (for $\gz$ with $N=0$) replaces Lemma \ref{lemestellip}, in order to get a Korn-type inequality
\begin{eqnarray*}
\Norm{Y}{1,2,-\delta} &\infeg& c \, (\Norm{\mathring{S}(Y) - \corec{\tr}_{\gz}(S(Y)) \gz}{2, -\delta} + \Norm{Y}{2,-\delta})\\
&\infeg& c \, (\Norm{\mathring{S}(Y)}{2, -\delta} + \Norm{Y}{2,-\delta}) \, ,
\end{eqnarray*}
where $c = c \, (\gz,\delta)$ is a positive constant.\\
We now use Proposition $D12$ from \cite{Chrusciel2004}: Let $(\m, \gz)$ be a conformally compact manifold with $\gz = \rho^{-2} \hz$ . For all $\delta \ne -(n+1) \slash 2$ and $ \delta \ne -(n-1) \slash 2$, there exists two constants $c_{\delta} >0$ and $\rho_{\epsilon, \delta} >0$ such that for all differentiable vector field $Y$ compactly supported in $\lbrace \rho < \rho_{\epsilon, \delta} \rbrace$ :
\begin{equation}\label{inegsy}
\Norm{Y}{2,-\delta} \infeg c_{\delta} \, \Norm{\mathring{S}(Y)}{2, -\delta}.
\end{equation}
Defining $R$ such that $\rho_{\epsilon, \delta} = e^{-2 R} $, we set $\Omega_{R}$ as in Definition \ref{cutoff} and we have (\ref{inegsy}) on $\m \setminus \Omega_{R}$, as in hypothesis of Theorem \ref{estellip}. The rest of the {proof} is analogous to the one of Theorem \ref{estellip} and we can choose $R$ large enough so that for all $Y \in \cinf_{c}(T^{*}\m)$ ,
\begin{equation*}
\Norm{Y}{1,2,-\delta} \infeg c \, \left(\Norm{\mathring{S}(Y)}{2,-\delta} + \Norm{Y}{2,-\delta;\Omega_{R}}\right).
\end{equation*}
Hence, the operator $\mathring{S} : \sob{1}{2}{-\delta}(T^{*}\m) \rightarrow \leb{2}{-\delta}(T^{*}\m)$ has a finite dimensional kernel. So we can write $\sob{1}{2}{-\delta}(T^{*}\m) = \text{ker} \mathring{S} \oplus (\text{ker} \mathring{S})^{\perp}$. From then, there exists a positive constant $c$ such that for all $Y \in (\text{ker} \mathring{S})^{\perp}$,
\begin{equation*}
\Norm{Y}{1,2,-\delta} \infeg c \, \Norm{\mathring{S}(Y)}{2,-\delta}.
\end{equation*}
It remains to show that \, ker $\mathring{S} \cap \sob{1}{2}{-\delta} = \lbrace 0 \rbrace$ , $\forall \delta > -(n+1) \slash 2$.\\
For that matter, we use the coordinate system $(x^{1} = \rho, x^{2}, \ldots , x^{n}) = (\rho,\theta)$ on a neighborhood of the boundary $\lbrack 0, \epsilon \rbrack \times \corec{\d_{\infty}\m} $ that we may find in \cite{Chrusciel2004}. From the expression of the defining function $\rho$, the metric $\gz$ can be written:
$$\gz = \rho^{-2} \hz = \rho^{-2} (d\rho^{2}+ \hat{g}(\rho)) \; \; \; \text{with} \; \; \hat{g}(\rho)(\d_{\rho},.) = 0. $$
The same conventions and notations as in the paper are in order: the index $\rho$ will be the radial coordinate one whereas indices relative to tangential coordinates will be designated by latin capital indices. Finally lower case latin indices designate any components. Christoffel symbols, in this coordinate system, are given in \cite{Chrusciel2004}, just as the equation
\begin{equation}\label{eqs} 
\nz_{i}X_{j} + \nz_{j}X_{i} = 0
\end{equation}
 which becomes the system:
\begin{comment}
\begin{equation}
\left\{
\begin{array}{l}
\d_{\rho} X_{\rho} + \rho^{-1} X_{\rho} = 0\\
\d_{\rho} X_{A} + \d_{A} X_{\rho} + 2 \rho^{-1} X_{A} - \hat{g}^{CD}(\rho)\hat{g}_{DA}'(\rho) X_{C} = 0\\
\d_{A} X_{B} + \d_{B} X_{A} - 2 \hat{\Gamma}_{A \;B}^{\hspace{0.1cm} C}(\rho) X_{C} +( \hat{g}_{AB}'(\rho) -2 \rho^{-1}  \hat{g}_{AB}(\rho)) X_{\rho} = 0
\end{array}
\right.
\end{equation}
\end{comment}
\begin{align}
&\d_{\rho} X_{\rho} + \rho^{-1} X_{\rho} = 0\label{polyhomo1}\\
&\d_{\rho} X_{A} + \d_{A} X_{\rho} + 2 \rho^{-1} X_{A} - \hat{g}^{CD}(\rho)\hat{g}_{DA}'(\rho) X_{C} = 0 \label{polyhomo2}\\
&\d_{A} X_{B} + \d_{B} X_{A} - 2 \hat{\Gamma}_{A \;B}^{\hspace{0.1cm} C}(\rho) X_{C} +( \hat{g}_{AB}'(\rho) -2 \rho^{-1}  \hat{g}_{AB}(\rho)) X_{\rho} = 0\nonumber
\end{align}
where $f':= \d_{\rho} f$.\\
Solving equation (\ref{polyhomo1}) give 
$$X_{\rho} = \rho^{-1} K(\theta) .$$
As the metric $\gz$ is polyhomogeneous , $\hat{g}$ can be written as a development of powers of $\rho$ and $\ln \rho$, and the first terms only contain powers of $\rho$. We set
$$ \hat{g}^{CD} \hat{g}_{DA}' := T^{C}_{\hspace{0.1cm} A} ,$$
where $T$ is and order two tensor whose development is
$$T^{C}_{\hspace{0.1cm} A}(\rho,\theta) = \rho \preciv{1} T^{C}_{\hspace{0.1cm} A}(\theta) + o(1).$$
Andersson and Chru\'{s}ciel have shown in \cite{Andersson1996} that the solution $X$ of 
$$(\Delta - Ric \, \gz )X = 0$$
 is also polyhomogeneous 
$$X_{A}(\rho,\theta) = \rho^{s} Z_{A}(\theta) + o(\rho^{s}).$$
Replacing in equation (\ref{polyhomo2}), we find $s=-2$ and $X_{A} = \rho^{-2} Z_{A}(\theta) + o(\rho^{-2})$.\\
We obtain the form of the solution of (\ref{eqs}) near the boundary
$$ X = \rho^{-1} K(\theta)\, d\rho + \lbrack \rho^{-2} Z_{A}(\theta) + o(\rho^{-2}) \rbrack \, dx^{A}.$$
Moreover, $X \in \sob{1}{2}{-\delta}$ , leading to $Z=0$.\\
Indeed, suppose $Z \ne 0$. The $\gz$-norm of $X$ is
$$\norm{X}{\gz}^{2} \simeq \rho^{-2}\norm{K}{\gz}^{2} + \rho^{-4}\norm{Z}{\gz}^{2} , $$
with $\norm{K}{\gz}^{2} \simeq \norm{Z}{\gz}^{2} \simeq \rho^{2}\, O(1)$.
$$\norm{X}{\gz}^{2} \simeq \rho^{-2} \, O(1). $$
Consequently
\begin{eqnarray*}
\Norm{X}{2,-\delta} < \infty &\equi& \int_{0}^{\epsilon} \norm{X}{\gz}^{2} \rho^{-2 \delta} \rho^{-n} d\rho < \infty\\
&\equi& \int_{0}^{\epsilon} \rho^{-2} \rho^{-2 \delta} \rho^{-n} d\rho < \infty\\
&\equi& -2 \delta - n - 2 >-1\\
&\equi& \delta  <-(n+1) \slash 2.
\end{eqnarray*}
Given that $\delta > -(n+1) \slash 2$ , then $X \notin \sob{1}{2}{-\delta}$. Hence, $Z$ is necessarily null near the boundary. Analysing the $X$ development coefficients, we realize $Z=0$ lead to $\norm{X}{\gz} = O(\rho^{\infty})$.
We conclude the {proof} using the unique continuation theorem from \cite{Mazzeo1991}.$\cqfd$ \\[.2cm]

\section{The contraint operator $\phibf$}

Let $\m$ be a $n$-dimensional connected non compact oriented manifold. We consider $\m$ as a \corec{spacelike  hypersurface} of a $(n+1)$-dimensional Lorentzian manifold $(\N, \gamma)$, from now on refered to as spacetime. We will distinguish the two manifolds by different indices: Latin indices will take values from $1$ to $n$ and are spatial indices whereas Greek indices will take values from $0$ to $n$ and are spacetime indices. $K$ is the second fundamental form of $\m$ in $\N$ defined by
\begin{equation}\label{defk}
K(X,Y) = \gamma (X, \precig{\gamma}\n_{Y} \vec n) ,
\end{equation}
where $\precig{\gamma}\n$ is the spacetime connection on $T\N$ , $X,Y \in T\m$ and $\vec n$ is the future-directed unit normal to $\m$ in $\N$. It is convenient to consider the conjugate momentum $\pi$ as a reparametrisation of $K$
\begin{equation}\label{defpi}
\pi^{ij} = \w{\pi}^{ij} \sqrt{g} \; \; \; \text{with} \; \; \w{\pi}^{ij} = K^{ij} - \corec{\tr}_{g}K g^{ij}.
\end{equation}
where $\sqrt{g}$ is the volume measure of the metric $g$.
\begin{equation*}\label{defsqrtg}
\sqrt{g} = \frac{\sqrt{det(g)}}{\sqrt{det(\gz)}} \, d\mu(\gz).
\end{equation*}
Remark: $\n (\sqrt{g}) = 0$ since the covariant derivative of the volume form is null.\\
$\w{\pi}$ is a section of the bundle $S^{2}T\m$ whereas $\pi$ is a section of the bundle $\w{S}=S^{2}T\m \otimes \Lambda^{3}T^{*}\m$.
We consider a smooth and polyhomogeneous asymptotically hyperbolic metric $\gz$ asymptotically hyperbolic on $\m$ as model. We ask $\gz$ to satisfy the following integrability condition
\begin{equation}
Riem \, \gz_{ijkl} - \gz_{il} \gz_{jk} + \gz_{ik} \gz_{jl}\in \lebd \label{integriem}.
\end{equation}
In particular,
\begin{equation}
Ric \, \gz_{jl} + (n-1) \gz_{jl} \in \lebd \label{integriccg}.
\end{equation}
\begin{comment}
We define the volume measure of $\gz$ with the Levi-Civita tensor $T_{LC}$ 
$$d\mu(\gz) =  T_{LC}.$$
In a direct basis, in dimension 3, components of the Levi-Civita tensor are given by
$$T_{abc} = \sqrt{det(\gz)} \, \Epsilon_{abc} \, ,$$ 
where $\Epsilon$ is the Levi-Civita symbol of order $3$.\\
A direct computation yields:
\begin{lemperso}{}\label{LV0}
The covariant derivative of the Levi-Civita tensor is null.
\end{lemperso}
%\begin{comment}
Preuve:
\begin{eqnarray*}
\nz_{i} T_{abc} &=& \d_{i}T_{abc} - \cdez{a}{k}{i} T_{k bc} - \cdez{b}{k}{i} T_{a k c} - \cdez{c}{k}{i} T_{abk}\\
&=& \mathbf{\Epsilon}_{abc} \, \d_{i}\sqrt{det(\gz)} - \sqrt{det(\gz)} \, \lbrack \cdez{a}{k}{i} \, \Epsilon_{k bc} - \cdez{b}{k}{i} \, \Epsilon_{a k c} - \cdez{c}{k}{i} \, \Epsilon_{abk}\rbrack\\
&=& \Epsilon_{abc} \, \d_{i}\sqrt{det(\gz)} - \sqrt{det(\gz)} \cdez{k}{k}{i} \, \Epsilon_{abc}\\
&=& \Epsilon_{abc} \lbrack \d_{i}\sqrt{det(\gz)} - \cdez{k}{k}{i}\sqrt{det(\gz)}\rbrack\\
&=& 0 \, ,
\end{eqnarray*}
using the expression of the contraction of the Christoffel symbol
$$\cdez{k}{k}{i} = \frac{\d_{i}\sqrt{det(\gz)}}{\sqrt{det(\gz)}} \cqfd$$
\end{comment}
For any sufficiently regular Riemannian metric $g$ on $\m$, we define the constraint operator $\phibf=(\phio , \phii) = \phibf (g,\pi)$ as follows:
\begin{eqnarray}
\phiogpi &:=& \left(R(g) - 2 \Lambda - \norm{K}{g}^{2} + (\corec{\tr}_{g}K)^{2}\right) \sqrt{g} \nonumber\\
&=& \left(R(g) - 2 \Lambda\right) \sqrt{g} - \left(\norm{\pi}{g}^{2} - \tfrac{1}{n-1} (\corec{\tr}_{g}\pi)^{2}\right) \slash \sqrt{g}.\\
\phiigpi &:=& 2 (\n^{j}K_{ij} - \n_{i}(\corec{\tr}_{g}K) )\sqrt{g} \nonumber\\
&=& 2 g_{ij} \n_{k} \pi^{jk} =  2 g_{ij} \n_{k} \w{\pi}^{jk} \sqrt{g}. \label{ECCG}
\end{eqnarray}

%A consequence of Lemma \ref{LV0} we used in (\ref{ECCG}) is that
%\begin{equation}\label{nulsqrtg}
%\n (\sqrt{g}) = 0 .
%\end{equation}
We set
\begin{equation}\label{defKz}
\Kz = \tau \gz \, ,
\end{equation}
where $\tau$ is a real parameter.\\
The cosmological constant $\Lambda$ is normalized here in dimension $n$ by
\begin{equation}\label{cosmocg}
2 \Lambda = n (n-1)(\tau^{2}-1),
\end{equation}
so that $\phibf(\gz,\Kz)=0$ at infinity.
Taking the $\gz$-trace of (\ref{integriccg}) considering (\ref{cosmocg}), we end up with the integrability condition
\begin{equation}\label{integscal}
R(\gz) - 2\Lambda + n(n-1)\tau^{2} \in \lebd.
\end{equation}

The conjugate momentum $\piz$ is then
\begin{equation}\label{defpiz}
\piz^{ij} =  (\Kz^{ij} - \corec{\tr}_{\gz}\Kz \gz^{ij}) \, d\mu(\gz) = \tau (1-n)\gz^{ij}\, d\mu(\gz) .
\end{equation}
Remark : $\nz \piz = \nz \Kz = 0$.\\
If the spacetime satisfies Einstein's equations, the normalisation chosen insures that the constraint operator and the energy-momentum tensor are related by $$ \phibf_{\alpha} = 16 \pi G T_{\vec n\alpha} \sqrt{g} \, ,$$
where $G$ is Newton's gravitationnal constant. $\xi = (N,X^{i})$ is the \emph{lapse-shift} associated to the spacetime foliation. 
We study the constraint operator $\phibf$ for Riemannian metrics of the form $g= \gz + h$. $g$ is asymptotic to $\gz \, , \; \ie \norm{g- \gz}{\gz} = \norm{h}{\gz} \underset{\infty}{\longrightarrow} 0$.\\

$\S:= S^{2} T^{*}\m$ is the bundle of symmetric bilinear forms on $\m$. $\w{\S}:= S^{2}T \m \otimes \Lambda^{3} T^{*}\m$ is the bundle of symmetric $2$-tensors-valued densities ($3$-forms) on $\m$. $\T := T\N$ is the spacetime tangent bundle. The following spaces will be of particular interest in the sequel:
\begin{equation*}
\begin{array}{ccl}
\G& := &\sob{2}{2}{\delta} (\S).\\
\K& := &\lbrace \pi : \pi -\piz \in \sob{1}{2}{\delta} (\w{\S}) \rbrace.\\
\G^{+}& := &\lbrace g : g-\gz \in \G , g>0\rbrace.\\
\G^{+}_{\lambda}& := &\lbrace g \in \G^{+} : \lambda \gz < g < \lambda^{-1} \gz\rbrace \; , \; 0 < \lambda < 1.\\
\L^{*}& := &\leb{2}{\delta}(\T^{*} \otimes \Lambda^{3} T^{*}\m) \, \, \text{is the dual space of} \, \, \L := \leb{2}{- \delta}(\T).
\end{array}
\end{equation*}
From $(\ref{inchold})$, tensors in $\G$ are Hölder-continuous
% of order $\demi$
and thus, matrices inequalities in spaces $\G^{+}$ and $\G^{+}_{\lambda} $ are satisfied pointwise. In particular, for all metric $g \in \G^{+}_{\lambda} $, metrics $g$ and $\gz$ are equivalent in the following sense
\begin{equation}\label{holdcontcg}
\lambda \gz_{ij}(x) \, v^{i}v^{j} < g_{ij}(x) \, v^{i}v^{j} < \lambda^{-1} \gz_{ij}(x) \, v^{i}v^{j} \; , \; \; \; \forall x \in \m \, , \forall v \in T\m
\end{equation}
So $\norm{g}{\gz} \simeq c \, \norm{\gz}{\gz} \simeq c \, \norm{g}{g} \simeq \sqrt{n}$.\\
$\F = \Gplus \times \K$ will be the phase space of the contraint operator $\phibf$. We will use $(g,\pi)$ as well as $(g,K)$ to express coordinates on $\F$.\\[.5cm]
Let $\cdez{}{}{}$ and $\nz$ (resp. $\cde{}{}{}$ and $\n$) be the Christoffel symbols and the Levi-Civita connection for $\gz$ (resp. $g$). We define
\begin{equation}\label{defA}
\ade{i}{k}{j} = \cde{i}{k}{j}-\cdez{i}{k}{j}
\end{equation}
Remark: $A$ is a symmetric tensor, arising from symmetry of the Christoffel symbols. We easily show 
\begin{equation}\label{A}
\ade{i}{k}{j} = g^{kl} \ape{ilj} = \demi g^{kl} (\nz_{i} g_{jl} + \nz_{j} g_{il} - \nz_{l} g_{ij})
\end{equation}

The scalar curvature of $g$ can be express with $\nz$ and $\ade{i}{j}{k}$ (see eq. (21) of \cite{Bartnik2005}):
\begin{lemperso}{}
\begin{eqnarray}\label{courbscalcov}
R(g) &=& g^{jk} Ric\, \gz_{jk} + g^{jk}(\nz_{i} \ade{j}{i}{k} - \nz_{j} \ade{i}{i}{k} + \ade{j}{l}{k}\ade{i}{i}{l} - \ade{j}{i}{l}\ade{k}{l}{i}) \nonumber\\
&=& g^{jk} Ric\, \gz_{jk} + Q(g^{-1}, \nz g) + g^{ik} g^{jl}(\nz^{2}_{ij} g_{kl} - \nz^{2}_{ik} g_{jl})
\end{eqnarray}
where Q is a sum of quadratic terms in $g^{-1}, \nz g$.
\end{lemperso}
This result relies on the following fact:
%\erw{mettre des ref pour les lemmes}
\begin{lemperso}{}\label{lemriccicov}
\begin{equation}\label{riccicov}
Ric \, g_{jk} - Ric \, \gz_{jk} = \nz_{i} \ade{j}{i}{k} - \nz_{j} \ade{i}{i}{k} + \ade{j}{\mu}{k}\ade{i}{i}{\mu} - \ade{j}{i}{\mu}\ade{k}{\mu}{i}
\end{equation}
\end{lemperso}
\textbf{Proof}:
\begin{eqnarray*}
Ric \, g_{jk} - Ric \, \gz_{jk} &=& \d_{i}\ade{j}{i}{k} - \d_{j}\ade{i}{i}{k} + \lbrack \cde{l}{i}{i}\cde{j}{l}{k} - \cdez{l}{i}{i}\cdez{j}{l}{k} \rbrack - \lbrack \cde{k}{l}{i}\cde{j}{i}{l} - \cdez{k}{l}{i}\cdez{j}{i}{l} \rbrack\\
&=& \d_{i}\ade{j}{i}{k} - \d_{j}\ade{i}{i}{k} + \lbrack \ade{l}{i}{i}\ade{j}{l}{k} + \cdez{l}{i}{i}\ade{j}{l}{k} + \cdez{j}{l}{k}\ade{l}{i}{i} \rbrack\\
&& - \lbrack \ade{k}{l}{i}\ade{j}{i}{l} + \cdez{k}{l}{i}\ade{j}{i}{l} + \cdez{j}{i}{l}\ade{k}{l}{i} \rbrack
\end{eqnarray*}
We end the proof adding and substracting $\cdez{j}{l}{i}\ade{l}{i}{k}$.$\cqfd$\\[.3cm]
\begin{comment}
\begin{eqnarray*}
Ric \, g_{jk} - Ric \, \gz_{jk} &=& \lbrack \d_{i}\ade{j}{i}{k} + \cdez{l}{i}{i}\ade{j}{l}{k} \mathbf{- \cdez{j}{l}{i}\ade{l}{i}{k}} - \cdez{k}{l}{i}\ade{j}{i}{l} \rbrack \\
&&- \lbrack \d_{j}\ade{i}{i}{k} + \cdez{l}{i}{j}\ade{i}{l}{k} - \mathbf{\cdez{j}{l}{i}\ade{l}{i}{k}} - \cdez{k}{l}{j}\ade{i}{i}{l} \rbrack \\
&& + \ade{l}{i}{i}\ade{j}{l}{k} - \ade{k}{l}{i}\ade{j}{i}{l} \; \; \; \; \;  \text{en ajoutant and retranchant le terme en gras}\\
&=& \nz_{i} \ade{j}{i}{k} - \nz_{j} \ade{i}{i}{k} + \ade{j}{l}{k}\ade{i}{i}{l} - \ade{j}{i}{l}\ade{k}{l}{i}. \cqfd
\end{eqnarray*}
\end{comment}

%La proposition suivante, analogue de la Proposition 3.1 de \cite{Bartnik2005} en asymptotiquement hyperbolique, 
Here we show $\phibf$ is a well-defined mapping between the Hilbert spaces $\F$ and $\L^{*}$.
\begin{propperso}{}\label{phidef} 
Set $(g,\pi) \in \G^{+}_{\lambda}\times \K$ , with $ 0 < \lambda < 1$. Then in dimension $n=3$, for all $\delta \infeg 0 $, there exists a positive constant $c = c \, (\lambda, \gz,\delta)$ such that
\begin{eqnarray}\label{16g}
\Norm{\phiogpi}{2,\delta} &\infeg& c \, \big(1 + \Norm{g-\gz}{2,2,\delta}^{2} + \Norm{\pi- \piz}{1,2,\delta}^{2}\big)\\
\Norm{\phiigpi}{2,\delta} &\infeg& c \, \Big(\Norm{\nz (\pi-\piz)}{2,\delta} +\Norm{\nz g}{1,2,\delta} (1+ \Norm{\pi- \piz}{1,2,\delta})\Big)
\end{eqnarray}
\end{propperso}
\textbf{Proof}: From $R(g)$ expression $(\ref{courbscalcov})$, we can write
\begin{eqnarray*}
\phiogpi &=& (R(g) - 2 \Lambda) \sqrt{g} - (\norm{\pi}{g}^{2} - \tfrac{1}{n-1} (\corec{\tr}_{g}\pi)^{2}\rbrack \slash \sqrt{g}\\
&=& \lbrack R(g) \mathbf{- R(\gz) + R(\gz)} - 2 \Lambda \mathbf{+ n(n-1) \tau^{2} - n(n-1) \tau^{2}} \rbrack \sqrt{g}\\
&& - \lbrack \norm{\pi- \piz}{g}^{2}- \norm{\piz}{g}^{2} + 2(\pi - \piz)_{ij}\piz^{ij} + 2 \, \norm{\piz}{g}^{2} \rbrack \slash \sqrt{g}\\
&& + \tfrac{1}{n-1} \lbrack (\corec{\tr}_{g}(\pi- \piz))^{2} + (\corec{\tr}_{g}\piz)^{2} + 2 \, \corec{\tr}_{g}(\pi- \piz) \, \corec{\tr}_{g}\piz \rbrack \slash \sqrt{g}\\
&=& \lbrack R(g) \mathbf{- R(\gz) + R(\gz)} - 2 \Lambda \mathbf{+ n(n-1) \tau^{2} - n(n-1) \tau^{2}} \rbrack \sqrt{g}\\
&& - \lbrack \norm{\pi- \piz}{g}^{2} + 2(\pi - \piz)_{ij}\piz^{ij} + \norm{\piz}{g}^{2} - \tfrac{1}{n-1}(\corec{\tr}_{\gz}\piz)^{2} \rbrack \slash \sqrt{g} \\
&& + \tfrac{1}{n-1} \lbrack (\corec{\tr}_{g}(\pi- \piz))^{2} + ((g- \gz)_{ij}\piz^{ij})^{2} + 2 (g- \gz)_{ij}\piz^{ij}\corec{\tr}_{\gz}\piz + 2 \, \corec{\tr}_{g}(\pi- \piz) \, \corec{\tr}_{g}\piz \rbrack \slash \sqrt{g}\\
&=& \lbrack R(g)- R(\gz) + R(\gz) - 2 \Lambda + n(n-1) \tau^{2}\rbrack \sqrt{g} - \lbrack \norm{\pi- \piz}{g}^{2} + 2(\pi - \piz)_{ij}\piz^{ij} \rbrack \slash \sqrt{g} \\
&& + \tfrac{1}{n-1} \lbrack (\corec{\tr}_{g}(\pi- \piz))^{2} + ((g- \gz)_{ij}\piz^{ij})^{2} + 2 (g- \gz)_{ij}\piz^{ij}\corec{\tr}_{\gz}\piz + 2 \, \corec{\tr}_{g}(\pi- \piz) \, \corec{\tr}_{g}\piz \rbrack \slash \sqrt{g}.
\end{eqnarray*}
Since $g \in \G^{+}_{\lambda}$, we can use (\ref{holdcontcg}) and from Cauchy-Schwarz inequality and $2ab \infeg a^{2} + b^{2}$
\begin{eqnarray*}
\norm{\phiogpi}{\gz} &\infeg& \lbrack \norm{R(g) - R(\gz)}{\gz} + \norm{R(\gz) - 2 \Lambda + n(n-1) \tau^{2}}{\gz} \rbrack \sqrt{g}\\
&& + c \, \lbrack 1 + \norm{\pi- \piz}{g}^{2}  + \norm{g- \gz}{g}^{2} \rbrack \slash \sqrt{g}.
\end{eqnarray*}
From (\ref{riccicov}), $Ric \, g - Ric \, \gz \simeq \nz A + A^{2} \simeq (\nz g) ^{2} + g \nz^{2}g + g^{-2} (\nz g) ^{2}$.\\
Using (\ref{19g}),
\begin{eqnarray*}
\Norm{Ric \, g - Ric \, \gz}{2, \delta} &\infeg& (c \, \Norm{\nz g}{1,2,\delta}^{2} + c \, \Norm{\nz ^{2} g}{2,\delta})\\
&\infeg& c \, \Norm{g - \gz}{2,2,-\delta}.
\end{eqnarray*}
In particular, we have the following integrability conditions
\begin{eqnarray}
Ric \, g - Ric \, \gz \in \lebd . \label{integdifric}\\
R(g) - R(\gz) \in \lebd . \label{integdifscal}
\end{eqnarray}
Thanks to (\ref{integdifscal}),(\ref{integscal}) and (\ref{19gbis}),
\begin{eqnarray*}
\Norm{\phiogpi}{2,\delta}&\infeg& c \left( 1 + \Norm{(\pi- \piz)^{2}}{2, \delta} +  \Norm{(g- \gz)^{2}}{2, \delta} \right)\\
&\infeg& c \left( 1 + \Norm{\pi- \piz}{1,2, \delta}^{2} +  \Norm{g- \gz}{2,2, \delta}^{2} \right).
\end{eqnarray*}
Hence, $\phiogpi \in \L^{*}$.\\
Regarding $\phiigpi$ , using (\ref{defA}),  
\begin{equation*}
\phiigpi = 2 g_{ij}(\nz_{k} (\pi- \piz)^{jk} + \ade{k}{j}{l} \, (\pi-\piz)^{kl}) + \ade{k}{j}{l} \, \piz^{kl}).
\end{equation*}
Considering (\ref{A}), $\phiigpi$ is of the form
\begin{equation}\label{phiidef}
\phiigpi \simeq g (\nz (\pi- \piz) + g^{-1} \nz g \, (\pi- \piz) +  g^{-1} \nz g \, \piz).
\end{equation}
\begin{eqnarray*}
\Norm{\phiigpi}{2,\delta} &\infeg& c \, (\Norm{\nz (\pi- \piz)}{2,\delta} + \Norm{\nz g \, (\pi- \piz)}{2,\delta} + \Norm{\nz g \, \piz}{2,\delta})\\
 &\infeg& c \, (\Norm{\nz (\pi- \piz)}{2,\delta} + \Norm{\nz g}{1,2,\delta} \Norm{\pi- \piz}{1,2,\delta} + \Norm{\nz g}{2,\delta}\Norm{\piz}{\infty,0}) \\
 &\infeg& c \, \big(\Norm{\nz (\pi- \piz)}{2,\delta} + \Norm{\nz g}{1,2,\delta} (1+ \Norm{\pi- \piz}{1,2,\delta})\big). \cqfd
\end{eqnarray*}
%La proposition suivante est l'analogue du Corollaire 3.2 de \cite{Bartnik2005} en asymptotiquement hyperbolique.
\begin{propperso}{}
Let $(g,\pi) \in \F$. Then in dimension $n=3$, for all $\delta \infeg 0 $,\\ $\phibf : \F \rightarrow \L^{*}$ is a smooth map between Hilbert spaces.
\end{propperso}
\textbf{Proof}: We recall the proof of \cite{Bartnik2005} for completeness. From Proposition \ref{phidef}, \\
$\Norm{\phigpi}{\L^{*}} \infeg c (1 + \Norm{g-\gz}{\G}^{2} + \Norm{\pi-\piz}{\K}^{2})$ , $\ie \phibf$ is locally bounded on $\F$. The polynomial structure of the constraint operator allows us to show $\phibf$ is smooth, $\ie$ indefinitely differentiable in a Fréchet sense. From the expression (\ref{courbscalcov}) of scalar curvature and given (\ref{phiidef}) , $\phibf$ can be expressed as
$$\phigpi = F (g, g^{-1}, \sqrt{g} , 1 \slash \sqrt{g} , \nz g, \nz^{2} g, \pi , \nz \pi) \, , $$
where $F= F(a_{1}, \ldots ,a_{8})$ is a polynomial function quadratic in $a_{5}$ and $a_{7}$ and linear in the remaining parameters. The map $g \mapsto (g, g^{-1}, \sqrt{g} , 1 \slash \sqrt{g})$ is analytic on the space of positive definite matrices and the maps $g \mapsto \nz g \, , \, g \mapsto \nz^{2} g$ and $\pi \mapsto \nz \pi$ are bounded linear, thus smooth, from $\F$ to $\L^{*}$, which are Hilbert spaces. A result from Hille \cite{Hille1957}  on locally bounded polynomial functionals shows $\phibf$ admit continuous Fréchet-derivatives of all orders.$\cqfd$\\

The set $\contrainte = \lbrace (g, \pi) \in \G^{+} \times \K : \phigpi = 0 \rbrace := \phibf^{-1}(\lbrace 0 \rbrace) \subset \F$ is the set of initial data for the vacuum Einstein's equations. To prove that $\contrainte$ is a submanifold of $\F$, we show that $0$ is a regular value of $\phibf$ , so we \corec{are interested} in the surjectivity of the differential of $\phibf$. 

\section{Expressions of the linearization of $\phibf$ and its adjoint}
%\erw{phrase and ref}
In this section we recall the expression of the linearization of $\phibf$ and its adjoint that we may find in \cite{Bartnik2005} or \cite{Fischer1979} for example.
\begin{propperso}
\begin{eqnarray}
D \phiogpi .(h,p) &=& (\n^{i} \n^{j} h_{ij} - \laplag \corec{\tr}_{g}h) \sqrt{g} - h_{ij} \lbrack R^{ij} - \tfrac{1}{2}  (R(g) - 2 \Lambda) g^{ij} \rbrack \sqrt{g}\nonumber\\
&& + h_{ij} \big(\tfrac{2}{n-1} \corec{\tr}_{g}\pi \pi^{ij} - 2 \pi^{i}_{\, k} \pi^{kj} + \tfrac{1}{2} \norm{\pi}{g}^{2} \, g^{ij} - \tfrac{1}{2(n-1)} (\corec{\tr}_{g}\pi)^{2} g^{ij} \big) \slash \sqrt{g} \nonumber\\
&& + p^{ij} (\tfrac{2}{n-1}\corec{\tr}_{g}\pi g_{ij} - 2 \pi_{ij}) \slash \sqrt{g}. \label {def0}\\[.3cm]
D \phiigpi .(h,p) &=& \pi^{jk}(2 \n_{k}h_{ij}- \n_{i} h_{jk}) + 2 h_{ij} \n_{k}\pi^{jk} + 2 g_{ik} \n_{j} p^{jk}. \label{defi}
\end{eqnarray}
\end{propperso}

Using notations of \cite{Bartnik2005} , 
\begin{eqnarray*}
\delta_{g}\delta_{g} h &=& \n^{i}\n^{j} h_{ij}.\\
E^{ij} &=& R^{ij} - \tfrac{1}{2}  (R(g) - 2 \Lambda) g^{ij}.\\
\Pi^{ij} &=& \big(\tfrac{2}{n-1} \corec{\tr}_{g}\pi \pi^{ij} - 2 \pi^{i}_{\, k} \pi^{kj} + \tfrac{1}{2} \norm{\pi}{g}^{2} \, g^{ij} - \tfrac{1}{2(n-1)} (\corec{\tr}_{g}\pi)^{2} g^{ij} \big) \slash (\sqrt{g})^{2}.
\end{eqnarray*}
%\begin{comment}
We can express $D \phibf$ in the matricial following form
\begin{equation}
D\phigpi.(h,p) =
\left[
\begin{array}{cc}
\sqrt{g}(\delta_{g}\delta_{g} - \laplag \corec{\tr}_{g} + \Pi - E) & -2 K\\
\hat{\pi}\n + 2 \delta_{g} \pi & 2 \delta_{g}
\end{array}
\right]
\left[
\begin{array}{c}
 h\\
p
\end{array}
\right] \, ,
\end{equation}
with $\; \hat{\pi} \n h = \hat{\pi}_{i}^{jkl} \n_{j} h_{kl} = (\pi^{jk} \delta_{i}^{l} + \pi^{jl} \delta_{i}^{k} - \pi^{kl} \delta_{i}^{j})\n_{j} h_{kl}$.\\[.5cm]
%\end{comment}
\begin{comment}
\begin{eqnarray*}
D \phiogpi .(h,p) &=& (\n^{i} \n^{j} h_{ij} - \laplag \corec{\tr}_{g}h) \sqrt{g} - h_{ij} E^{ij} \sqrt{g}\\
&& + h_{ij} S^{ij} \sqrt{g} + p^{ij} (\corec{\tr}_{g}\pi g_{ij} - 2 \pi_{ij}) \slash \sqrt{g}\\[.3cm]
D \phiigpi .(h,p) &=& \pi^{jk}(2 \n_{j}h_{ik}- \n_{i} h_{jk}) + 2 h_{ij} \n_{k}\pi^{jk} + 2 g_{ik} \n_{j} p^{jk}
\end{eqnarray*}
\end{comment}
To prove surjectivity of the differential of $\phibf$ , we investigate injectivity of the adjoint operator. Integrating by parts and ignoring boundary terms leads (cf. \cite{Fischer1979} for example) to the expression of the formal $L^{2}(d \mu (\gz))$-adjoint of $D \phigpi$ :
\begin{equation*}
\int_{\m} D \phigpi .(h,p) \, (N,X^{i}) = \int_{\m} (h,p) \bullet D \phigpis (N,X^{i}).
\end{equation*}

\begin{propperso}{}
\begin{eqnarray*}\label{adjphiocg}
(h,p) \bullet D \phiogpis N &=& h_{ij} \lbrack \n^{i} \n^{j} N - g^{ij}\laplag N - \lbrack R^{ij} - \tfrac{1}{2}  (R(g) - 2 \Lambda) g^{ij} \rbrack N \rbrack \sqrt{g}\\
&&+ N h_{ij} \big(\tfrac{2}{n-1} \corec{\tr}_{g}\pi \pi^{ij} - 2 \pi^{i}_{\, k} \pi^{kj} + \tfrac{1}{2} \norm{\pi}{g}^{2} \, g^{ij} - \tfrac{1}{2(n-1)} (\corec{\tr}_{g}\pi)^{2} g^{ij} \big) \slash \sqrt{g}\\
&& +N p^{ij} (\tfrac{2}{n-1} \corec{\tr}_{g}\pi g_{ij} - 2 \pi_{ij}) \slash \sqrt{g}.\\
(h,p) \bullet D \phiigpis X^{i} &=& h_{ij}(X^{k} \n_{k} \pi^{ij} +   \n_{k} X^{k} \pi^{ij} - 2  \n_{k}  X^{(i} \pi^{j)k}) -2 p^{ij} \n_{(i} X_{j)}.
\end{eqnarray*}
\end{propperso}

Then we can put $D \phibf^*$ in the matricial form
\begin{equation}\label{dphis}
D\phigpis.(N,X) =
\left[
\begin{array}{cc}
\sqrt{g}(\n^{2} - g \laplag + \Pi - E) & \n \pi - \hat{\pi}\n\\
 -2 K & - \L_{g}
\end{array}
\right]
\left[
\begin{array}{c}
 N\\
X
\end{array}
\right] \, ,
\end{equation}
with
\begin{eqnarray*}
(\n \pi - \hat{\pi}\n)X &=& \L_{X}\pi = \n_{X}\pi^{ij} - \hat{\pi}_{l}^{kij} \n_{k} X^{l}.\\
\L_{g}(X)&=& \L_{X}g = 2 \, \n_{(i}X_{j)} = 2 \, S(X).
\end{eqnarray*}
$D\phigpis_{1}.\xi$ and $D\phigpis_{2}.\xi \;$ \corec{will denote} the two components of $D\phigpis$ in (\ref{dphis}).\\
 $\lebd \xi$ (resp. $\sob{1}{2}{\delta} \nz \xi$) is the set of terms of the form $u \, \xi$ (resp. $u \, \nz \xi$) such that $\Norm{u}{2,\delta} \infeg C$ (resp. $\Norm{u}{1,2,\delta} \infeg C$ ), where $C$ is a constant depending on $\gz, \delta$ and $\Norm{(g,\pi)}{\F}$.
\begin{eqnarray}
D\phigpis_{1}.\xi &=& \lbrack \n_{i} \n_{j} N - g_{ij}\laplag N + (\Pi_{ij} - E_{ij}) \rbrack N \rbrack \sqrt{g} + (\n \pi - \hat{\pi}\n)X \nonumber\\
&=& D \phibf (g,0)^{*} \, (N,0) + \Pi_{ij} N \sqrt{g} + (\n \pi - \hat{\pi}\n)X \label{expdphi1dphi0} \\[.3cm]
(\n \pi - \hat{\pi}\n)X&=&X^{k}\n_{k}\pi_{ij} - (\pi^{k}_{\, i} \delta_{lj} + \pi^{k}_{\, j} \delta_{li} - \pi_{ij} \delta_{l}^{k}) \n_{k} X^{l} \nonumber\\
%&=& X \nz \pi + XA \pi - (\pi^{k}_{\, i} \delta_{lj} + \pi^{k}_{\, j} \delta_{li} - \pi_{ij} \delta_{l}^{k}) \nz_{k} X^{l}\\
%&=& X \nz (\pi-\piz) + XA (\pi - \piz) + XA \piz +  (\pi - \piz) \nz X - (\piz^{k}_{\, i} \delta_{lj} + \piz^{k}_{\, j} \delta_{li} - \piz_{ij} \delta_{l}^{k}) \n_{k} X^{l}\\
&=& \lebd X +  \sob{1}{2}{\delta} \nz X +(n-1) \tau (2 \Sz (X) - \gz \,\corec{\tr}_{\gz}\Sz(X)). \nonumber\\[.3cm]
\Pi(g, \pi) N &=& \lebd N + \Pi(\gz, \piz) N \nonumber\\
%&=& \lebd N + \big(\tfrac{2}{n-1} \corec{\tr}_{\gz}\piz \piz^{ij} - 2 \piz^{i}_{\, k} \piz^{kj} + \tfrac{1}{2} \norm{\piz}{\gz}^{2} \, \gz^{ij} - \tfrac{1}{2(n-1)} (\corec{\tr}_{\gz}\piz)^{2} \gz^{ij} \big) N\slash (\sqrt{\gz})^{2}\\
&=& \lebd N + - \tdemi(n-1)(n-4)\tau^{2} \gz N.\nonumber
\end{eqnarray}
So we have the integrability condition
\begin{equation}\label{integPi}
\Pi(g, \pi) + \tdemi(n-1)(n-4)\tau^{2} \gz \in \lebd.
\end{equation}
Taking into account (\ref{integriccg}) and (\ref{integscal}),
\begin{equation}\label{integE}
E + (n-1) \gz -\tdemi n(n-1)\tau^{2} \gz \in \lebd.
\end{equation}
On one hand,
\begin{eqnarray}
D\phigpis_{1}.\xi \slash \sqrt{g} &=& \n^{2} N - g\laplag N + (n-1) \gz +\lbrack \Pi + \tdemi(n-1)(n-4)\tau^{2} \gz \rbrack N \nonumber\\
&& +(n-1) \tau (2 \Sz (X) - \gz \,\corec{\tr}_{\gz}\Sz(X))  - \lbrack E+ (n-1) \gz -\tdemi n(n-1)\tau^{2} \gz \rbrack N \rbrack \nonumber\\
&& - (n-1)(n-2)\tau^{2} \gz N + \lebd \xi +  \sob{1}{2}{\delta} \nz X \nonumber\\
&=& \n^{2} N - g\laplag N + (n-1) g + (n-1) \tau (2 \Sz (X) - \gz\, \corec{\tr}_{\gz}\Sz(X)) \nonumber\\
&& - (n-1)(n-2)\tau^{2} \gz N + \lebd \xi +  \sob{1}{2}{\delta} \nz X. \label{expdphi1}
\end{eqnarray}
On the other hand,
\begin{eqnarray*}
D\phigpis_{2}.\xi &=& -2 KN - 2 S(X)\\
%&=& -2 \Kz N - 2 \Sz (X) -2 (K-\Kz) N + 2 A X\\
&=& - 2 (\Sz (X) + \tau \gz N) + \sob{1}{2}{\delta} \xi.
\end{eqnarray*}
From the definition of the operator $T = \n^{2}N -Ng$ and the expression of $\Sz$ as a function of $D\phigpis_{2}.\xi$, we obtain
\begin{eqnarray}\label{dphis12}
D\phigpis_{1}.\xi \slash \sqrt{g}&=& T- g\,\corec{\tr}_{g}T + (n-1) \tau (2 \Sz (X) - \gz\, \corec{\tr}_{\gz}\Sz(X))\nonumber\\
&& - (n-1)(n-2)\tau^{2} \gz N + \lebd \xi +  \sob{1}{2}{\delta} \nz X. \label{dphis1}\\
D\phigpis_{2}.\xi &=& - 2 (\Sz (X) + \tau \gz N) + \sob{1}{2}{\delta} \xi. \label{dphis2}
\end{eqnarray}

It is useful to restructure $D\phibfs$ into the operator $\Ps$ defined by
\begin{eqnarray}\label{Ps}
\Ps(\xi) = \Pgpis(\xi) &=&
\left[
\begin{array}{c}
g^{1 \slash 4}\left(\n^{i}\n_{j}N - \delta^{i}_{\, j} \laplag N + (\Pi^{i}_{\, j} - E^{i}_{\, j})N\right) + g^{-1\slash 4}\L_{X}\pi^{i}_{\, j}\\
 - 2 g^{-1 \slash 4} \n_{l}( K^{i}_{\, j} N +  S(X)^{i}_{\, j})
\end{array}
\right]
 \nonumber \\
&=& \zeta \circ
\left[
\begin{array}{cc}
1 & 0\\
0 & \n
\end{array}
\right]
\circ D\phigpis \xi \, ,
\end{eqnarray}
where $g^{1 \slash 4} = (det(g) \slash det(\gz))^{1 \slash 4} \, d\mu(\gz)$ is a density of weight $\tdemi$ and
\begin{equation}\label{defzeta}
\zeta = \zeta(g) =
\left[
\begin{array}{cc}
g^{-1 \slash 4} g_{jk} & 0\\
0 & g^{1 \slash 4} g^{ik}
\end{array}
\right].
\end{equation}
Finally, we can put $\Pgpis(\xi)$ into the form
\begin{equation}\label{defPs}
\Pgpis(\xi) =
\left(
\begin{array}{c}
g^{-1 \slash 4} \, D\phigpis_{1}.\xi  \\
 g^{1 \slash 4} \, \n D\phigpis_{2}.\xi 
\end{array}
\right).
\end{equation}
Expression (\ref{Ps}) of $\Ps$ allows us to rewrite the $L^{2}(d \mu (\gz))$-adjoint of $\Ps$ as follows
\begin{equation}\label{defP}
\Pgpi = D\phigpi \circ 
\left[
\begin{array}{cc}
1 & 0\\
0 & - \delta_{g}
\end{array}
\right]
\circ \zeta \, ,
\end{equation}
with $\delta_{g}q = \n^{l}(q_{l}^{ij})$ so that $P(f^{i}_{\, j},q_{li}^{\;j}) = D\phibf(f_{ij},q_{l}^{ij})$ and so the composition $P \Ps$ is well defined.\\
%La proposition suivante est l'analogue en asymptotiquement hyperbolique de la Proposition 3.3 de \cite{Bartnik2005}.

\section{Elliptic estimates relative to the adjoint}
In this section, we gather elliptic estimates satisfied by the adjoint operator $D \phibfs$.
\begin{propperso}{}\label{prop3.3g}
Set $\delta \in \rbrack -(n+1) \slash 2 \, , 0 \rbrack$, with $n=3$,  and $ \delta \ne -(n-1) \slash 2$. There exists a positive constant  
$C = C \, (\gz, \lambda, \delta , \Norm{g}{\F})$ such that the following elliptic estimate is satisfied:
 $\forall \xi \in \sob{2}{2}{-\delta} (\T)$ ,
\begin{equation}\label{35cg}
\Norm{\xi}{2,2,-\delta} \infeg c \, \big(\Norm{D \phigpis_{1}.\xi }{2, -\delta} + \Norm{D \phigpis_{2}.\xi }{1,2, -\delta} \big) + C \, \Norm{\xi}{1,2,-2\delta} \, ,
\end{equation}
\end{propperso}
\textbf{Proof}: Considering expression (\ref{dphis2}) of $\Sz$ as a function of $D\phigpis_{2}.\xi$,
\begin{eqnarray}\label{T-gtrT}
T- g\,\corec{\tr}_{g}T &=& D\phigpis_{1}.\xi \slash \sqrt{g} + (n-1)\tau \big(D\phigpis_{2}.\xi - \tdemi \gz \,\corec{\tr}_{\gz} D\phigpis_{2}.\xi \big) \nonumber\\
&& + \lebd \xi +  \sob{1}{2}{\delta} \nz X
\end{eqnarray}
From Proposition \ref{propn2infTz} and (\ref{h1}),
\begin{eqnarray*}
\Norm{N}{2,2,-\delta} &\infeg& c \, \big( \Norm{D\phigpis_{1}.\xi}{2,-\delta} +  \, (n-1) \tau \, (1+ \tfrac{n^{2}}{4})^{\demi} \Norm{D\phigpis_{2}.\xi}{2,-\delta} \big)\\
&& + C \, (\Norm{\xi}{\infty,-2\delta} + \Norm{\nz \xi}{3,-2\delta}).
\end{eqnarray*}
Using (\ref{uinfd}), (\ref{u3d}) and Sobolev inclusion ($\delta \infeg 0$), there exists a positive constant\\  
$C = C \, (\gz, \lambda, \delta , \Norm{g}{\F})$ such that
\begin{eqnarray}\label{estimeeN2}
\Norm{N}{2,2,-\delta} &\infeg& c \, \big( \Norm{D\phigpis_{1}.\xi}{2,-\delta} +  \, (n-1) \tau \, (1+ \tfrac{n^{2}}{4})^{\demi} \Norm{D\phigpis_{2}.\xi}{2,-\delta} \big) \nonumber\\
&& + \epsilon \; \Norm{\xi}{2,2,-\delta} + C \, \Norm{\xi}{1,2,-2\delta}.
\end{eqnarray}
Besides, for every sufficiently regular 1-form $X$ on $\m$, we have the following identity for the metric $\gz$ (cf equation (29) of \cite{Bartnik2005} for example)
\begin{equation}\label{b29}
\nz^{2}_{kj}X_{i} :=\nz_k\nz_jX_i= Riem \, \gz_{ijkl} X^{l} + \nz_{k} \Sz (X)_{ij} + \nz_{j} \Sz (X)_{ik} - \nz_{i} \Sz (X)_{jk},
\end{equation}
Hence
\begin{eqnarray}\label{b29bis}
\Norm{\nz^{2}X}{2,-\delta} &\infeg& \Norm{Riem \,\gz \; X}{2,-\delta} + c \, \Norm{\nz \Sz (X)}{2,-\delta} \nonumber\\
 %&\infeg& \Norm{Riem \,\gz}{\infty,0} \Norm{X}{2,-\delta} + c \, \Norm{\nz \Sz (X)}{2,-\delta} \nonumber\\
 &\infeg& c \, \Norm{X}{2,-\delta} + c \, \Norm{\nz \Sz (X)}{2,-\delta}.
\end{eqnarray}
A consequence of Lemma \ref{lemyinfcs} is 
\begin{equation}\label{conslemyinfcs}
\Norm{X}{2,-\delta} \infeg \Norm{X}{1,2,-\delta} \infeg c \, \Norm{\Sz(X)}{2,-\delta},
\end{equation}
which imply, with (\ref{b29bis})
\begin{equation}\label{b29ter}
\Norm{\nz^{2}X}{2,-\delta} \infeg c \, \Norm{\Sz(X)}{1,2,-\delta},
\end{equation}
and considering (\ref{conslemyinfcs}),
\begin{equation*}
\Norm{X}{2,2,-\delta} \infeg c_{1} \, \Norm{\Sz (X)}{1,2,-\delta}.
\end{equation*}
From (\ref{dphis2}) and (\ref{uinfd}),
\begin{equation}\label{estsx}
\Norm{\Sz (X)}{1,2,-\delta} \infeg \tfrac{1}{4}\Norm{D\phigpis_{2}.\xi}{1,2,-\delta} + n \tau \, \Norm{N}{1,2,-\delta} + \epsilon \Norm{\xi}{2,2,-\delta} + C \, \Norm{\xi}{1,2,-2\delta}.
\end{equation}
Thus there exists a constant $C$ depending on $\gz, \lambda, \epsilon, \delta$ and $\Norm{(g,\pi)}{\F}$ such that
\begin{equation}\label{estimeeX2}
\Norm{X}{2,2,-\delta} -  n c_{1} \, \tau \, \Norm{N}{1,2,-\delta} \infeg \tfrac{c_{1}}{4} \, \Norm{D\phigpis_{2}.\xi}{1,2,-\delta} + \epsilon \Norm{\xi}{2,2,-\delta} + C \, \Norm{\xi}{1,2,-2\delta}.
\end{equation}
We can choose $\epsilon_{0}<<1$ so that $(\ref{estimeeN2}) + \epsilon_{0} (\ref{estimeeX2})$ combine to yield (\ref{35cg}).$\cqfd$\\[.4cm]
Remark: it may be possible to extend this result to $\delta = -(n-1) \slash 2$ using operator $\Uz$ introduced later on.\\[.3cm]
Combining Proposition \ref{prop3.3g} and Ehrling inequality (\ref{Ehrling}), we get
\begin{corperso}{}\label{corestiadjL2}
Let $\delta \in \rbrack -(n+1) \slash 2 \, , 0 \rbrack$, with $n=3$ and $\delta \ne -(n-1) \slash 2$.\\
Then the following estimate is verified:
 $\forall\xi \in \sob{2}{2}{-\delta} (\T)$ , 
\begin{equation}\label{4.22biscg}
\Norm{\xi}{2,2,-\delta} \infeg c \, \big(\Norm{D \phigpis_{1}.\xi }{2, -\delta} + \Norm{D \phigpis_{2}.\xi }{1,2, -\delta} \big) + C \, \Norm{\xi}{2,-2\delta} \, ,
\end{equation}
where $C$ depends on $\gz, \lambda, \delta$ and $\Norm{(g,\pi)}{\F}$.
\end{corperso}

The next lemma will be very useful during the proof of proposition \ref{proplipcg} since it is the Time-symmetric version.

%La proposition suivante est l'analogue en asymptotiquement hyperbolique %de la Proposition 3.4 de \cite{Bartnik2005}.
\begin{lemperso}{}
In dimension $n=3$, let $\delta \infeg 0$ , then the operator\\ $D \phibf (g,0)^{*} \, (.,0) : \sob{2}{2}{-\delta}(\m) \longrightarrow \leb{2}{-\delta} (\w{\S}) \;$ is bounded and depends on $g$ in a Lipschitz way,
\begin{equation}\label{lip}
\bNorm{\big \lbrack D \phibf (g,0)^{*} - D \phibf (\wg,0)^{*} \big \rbrack \, (N,0)}{2, -\delta} \infeg C \Norm{g - \tilde{g}}{\F} \; \Norm{N}{2,2,-\delta} \, ,
\end{equation}
where constant $C$ depends on $\gz , \delta , \Norm{g}{\F}$ and $\Norm{\tilde{g}}{\F}$.
\end{lemperso}

\textbf{Proof}: Let us recall the expression of $D \phibf (g,0)^{*}$
\begin{equation}\label{adjphio}
D \phibf (g,0)^{*}. (N,0) = \lbrack \n_{i} \n_{j} N - g_{ij}\laplag N - \lbrack R_{ij} - \tfrac{1}{2}  (R(g) - 2 \Lambda) g_{ij} \rbrack N \rbrack \sqrt{g}.
\end{equation}
Let us begin by showing $D \phibf (g,0)^{*} \;$ is bounded.
In order to do this, we introduce the operator $O$ acting on functions
\begin{equation}\label{defO} 
O(N)= \n^{2}N - g \, \laplag N
\end{equation}
and we notice that $O(N)= L (\n^{2}N)$ where $L$ is a linear invertible operator, so
\begin{eqnarray*}
\Norm{O}{2,-\delta} &\infeg& c \, \Norm{\n^{2}N}{2, -\delta} =  c \, \left(\Norm{\nz^{2}N}{2, -\delta} + \Norm{A\, d N}{2, -\delta} \right)\\
& \infeg & C \, \Norm{N}{2,2, -\delta}.
\end{eqnarray*}
Indeed, $A\, d N \simeq g^{-1} \nz g \; d N$. Using Hölder inequality (\ref{h1}) , $(\ref{19g})$ and Sobolev inclusion,
\begin{eqnarray}\label{adnbis}
\Norm{A\, d N}{2, -\delta} & \infeg& \Norm{g^{-1}}{\infty,0} \Norm{\nz g \, d N }{2, -\delta} \nonumber\\
& \infeg& c \, \Norm{\nz g}{1,2, \delta} \Norm{d N}{1,2, -2 \delta} \nonumber\\
& \infeg& C \, \Norm{N}{2,2, -\delta}.
\end{eqnarray}
\begin{eqnarray*}\label{dphiog}
\Norm{ D \phibf (g,0)^{*}. (N,0) \diagup \sqrt{g}}{2,-\delta} &\infeg& \Norm{O}{2, -\delta} +  \Norm{\left(Ric \, g - Ric \, \gz \right) N}{2,-\delta} + \Norm{(n-1) \gz N}{2,-\delta} \nonumber\\
&& + \Norm{\left\lbrack Ric \, \gz + (n-1) \gz \right\rbrack N}{2,-\delta} + \Norm{\tdemi n(n-1)\tau^{2}\,g \, N}{2,-\delta}\\
&& + \tdemi \bNorm{\left\lbrack R(g) - 2 \Lambda + n(n-1)\tau^{2} \right \rbrack\,g \, N}{2,-\delta}.
\end{eqnarray*}
Considering (\ref{integdifric}), (\ref{integriccg}),
\begin{eqnarray*}
\Norm{(Ric \, g - Ric \, \gz) N}{2, -\delta} &\infeg& C \,\Norm{N}{2,2,-\delta}\\
\Norm{\lbrack Ric \, \gz + (n-1) \gz \rbrack N}{2,-\delta}  &\infeg& c \, \Norm{N}{2,2,-\delta}.
\end{eqnarray*}
For the scalar curvature term, using Hölder inequality (\ref{h1}) and Sobolev inclusion together with (\ref{integscal}) and (\ref{integdifscal}),
\begin{eqnarray*}
\Norm{(R(g) - 2 \Lambda + n(n-1)\tau^{2} )\,g\, N}{2, -\delta} &\infeg& \Norm{(R(\gz) - 2 \Lambda + n(n-1)\tau^{2} )\,g\, N}{2, -\delta}\\
&&+ \Norm{(R(g) - R(\gz) )\,g\, N}{2, -\delta}\\
%&\infeg& \Norm{g}{\infty, 0} \bNorm{\left \lbrack R(g) - 2 \Lambda + n(n-1)\tau^{2} \right \rbrack N}{2, -\delta}\\
%&\infeg& c \, \Norm{R(g) - 2 \Lambda + n(n-1)\tau^{2} }{2, \delta} \Norm{N}{\infty, -2\delta}\\
&\infeg& C \, \Norm{N}{2,2,-\delta}.
\end{eqnarray*}
\begin{eqnarray*}
\Norm{(n-1) \gz N}{2, -\delta} &\infeg& (n-1) \Norm{\gz}{\infty, 0} \Norm{N}{2, -\delta}\\
&\infeg& c \, \Norm{N}{2,2,-\delta}.
\end{eqnarray*}
Similarly,
\begin{eqnarray*}
\Norm{\tdemi n(n-1)\tau^{2} g N}{2, -\delta} &\infeg& c \, \Norm{N}{2,2,-\delta}.
\end{eqnarray*}
We end up with
\begin{equation*}%\label{majdphiosqrtg}
\Norm{ D \phibf (g,0)^{*}. (N,0) \diagup \sqrt{g}}{2,-\delta} \infeg C \, \Norm{N}{2,2,-\delta}
\end{equation*}
and finally
\begin{equation}\label{majdphio}
\Norm{D \phibf (g,0)^{*}. (N,0)}{2,-\delta} \infeg C \, \Norm{\sqrt{g}}{\infty,0} \Norm{N}{2,2,-\delta} \infeg C \, \Norm{N}{2,2,-\delta} ,
\end{equation}
where $C$ is a constant depending upon $\gz , \delta$ and $\Norm{g}{\F}$.\\[.2cm]
Proof of (\ref{lip}):
$\wn \; , \wlapla \; , Ric (\wg) \; \text{and} \; R(\wg) \;$ will denote The Levi-Civita connection, the Laplacian, the Ricci tensor and the scalar curvature of the Riemannian metric $\wg$.\\
In order to lighten notations, we set $$D\phiogs N := D \phibf (g,0)^{*} \, (N,0) \; \; \mbox{and} \; \; D \phio (\wg)^{*} N := D \phibf (\wg,0)^{*} \, (N,0)$$
\begin{equation*}
\lbrack D\phiogs - D \phio (\wg)^{*}\rbrack N = (\sqrt{g} - \sqrt{\w{g}}) \, \frac{D\phiogs N}{\sqrt{g}} + \sqrt{\w{g}} \left\lbrack \frac{D\phiogs N}{\sqrt{g}} - \frac{D\phio (\wg)^{*} N}{\sqrt{\wg}} \right\rbrack.
\end{equation*}
\begin{eqnarray}\label{declip}
\bNorm{\lbrack D\phiogs - D\phio (\tilde{g})^{*}\rbrack N}{2, -\delta} 
%&\infeg& \Norm{\sqrt{g} - \sqrt{\w{g}}}{\infty, \delta} \BNorm{\frac{D\phiogs N}{\sqrt{g}}}{2,-2\delta} \nonumber\\
%&&+ \Norm{\sqrt{\w{g}}}{\infty,0} \BNorm{\Big( \frac{D\phiogs N}{\sqrt{g}} - \frac{D\phio (\wg)^{*} N}{\sqrt{\wg}}\Big) }{2,-\delta}\nonumber\\
&\infeg& \Norm{g -\w{g}}{\F} \BNorm{\frac{D\phiogs N}{\sqrt{g}}}{2,-2\delta} \nonumber\\
&&+ c \, \BNorm{\frac{D\phiogs N}{\sqrt{g}} - \frac{D\phio (\wg)^{*} N}{\sqrt{\wg}}}{2,-\delta}.
\end{eqnarray}
\begin{eqnarray*}
\Big( \frac{D\phiogs N}{\sqrt{g}} - \frac{D\phio (\wg)^{*} N}{\sqrt{\wg}}\Big) &=& (\n - \wn) \, d N + g\, \laplag N - \wg \, \wlapla N - \lbrack Ric (g) - Ric (\wg)\rbrack N \\
&&+ \tdemi \left\lbrack (R(g) - 2\Lambda)g - (R(\wg) - 2\Lambda) \wg \right\rbrack N.
\end{eqnarray*}
\begin{eqnarray*}
\BNorm{\frac{D\phiogs N}{\sqrt{g}} - \frac{D\phio (\wg)^{*} N}{\sqrt{\wg}} }{2,-\delta} &\infeg& \Norm{(\n - \wn) \, d N}{2, -\delta} + \Norm{g \, \laplag N - \wg \, \wlapla N}{2,- \delta}\\
&& + \tdemi \bNorm{\left\lbrack (R(g) - 2\Lambda) g - (R(\wg) - 2\Lambda) \wg \right\rbrack N}{2, -\delta}\\
&& - \Norm{\lbrack Ric (g) - Ric (\wg)\rbrack N}{2, -\delta}.
\end{eqnarray*}
\begin{list}{$\bullet$}
\item For the Hessian:
\begin{equation}\label{hess}
\n - \wn = (g^{-1} - \wg^{-1}) \nz g + \wg^{-1} \nz (g - \wg).
\end{equation}
Using Hölder and Sobolev weighted inequalities, Sobolev inclusion ($\delta \infeg 0$) and (\ref{19g}),
\begin{equation}\label{majorhess}
\Norm{(\n - \wn) \, d N}{2, -\delta} 
%& \infeg& \Norm{(g^{-1} - \wg^{-1}) \nz g \; d N}{2, -\delta} + \Norm{\wg^{-1} \nz (g - \wg) \, d N }{2, -\delta}\\
%& \infeg& c \, \Norm{g - \wg}{\infty, \delta} \Norm{\nz g \; d N}{2, -2 \delta} + \Norm{\wg^{-1}}{\infty,0} \Norm{\nz (g - \wg) \, d N }{2, -\delta}\\
%& \infeg& c \, \Norm{g - \wg}{2,2, \delta} \Norm{\nz g \; d N}{2, 0} + c \Norm{\nz (g - \wg) \, d N }{2, 0}.
%\end{eqnarray*}
%En utilisant ,
%\begin{eqnarray}
%\Norm{(\n - \wn) \, d N}{2, -\delta} & \infeg& c \, \Norm{g - \wg}{2,2, \delta} \Norm{\nz g}{1,2, \delta} \Norm{d N}{1,2, -\delta} + c \, \Norm{\nz (g - \wg)}{1,2, \delta} \Norm{d N }{1,2, -\delta} \nonumber\\
 \infeg C \, \Norm{g - \wg}{2,2, \delta} \Norm{N}{2,2, -\delta}.
\end{equation}
\item \item For the Laplacian:
\begin{eqnarray*}
g\, \laplag N - \wg \, \wlapla N &=& g\, \laplag N - \wg \, \laplag N + \wg \, \laplag N - \wg \, \wlapla N\\
&=& (g - \wg ) \laplag N + \wg (\laplag N - \wlapla N)\\
%&=& (g  - \wg ) g^{-1} \n d N + \wg (g^{-1} \n d N - \wg^{-1} \wn d N)\\
&=& (g - \wg ) g^{-1} \n d N + \wg  (g^{-1} - \wg^{-1}) \n d N + \wg \wg^{-1}(\n - \wn) d N.
\end{eqnarray*}
Using Hölder inequality (\ref{h1}), Sobolev inclusion ($\delta \infeg 0$) and Sobolev inequality,
\begin{eqnarray*}
\Norm{g\,\laplag N - \wg\, \wlapla N}{2, -\delta}
% & \infeg& \Norm{(g - \wg ) g^{-1} \n d N}{2, -\delta} + \Norm{\wg (g^{-1} - \wg^{-1}) \n d N }{2, -\delta}\\
%&& + \Norm{\wg \wg^{-1}(\n - \wn) d N}{2, -\delta}\\
%& \infeg& \Norm{g- \wg}{\infty, \delta} \Norm{g^{-1}}{\infty,0}  \Norm{\n d N}{2, -2\delta}\\
%&& + \Norm{g^{-1} - \wg^{-1}}{\infty, \delta} \Norm{\wg}{\infty,0}  \Norm{\n d N}{2, -2\delta}\\
%&&+\Norm{\wg}{\infty,0}\Norm{\wg^{-1}}{\infty,0}\Norm{(\n - \wn) d N}{2, -\delta}\\
& \infeg& c \, \Norm{g - \wg}{2,2, \delta} \Norm{\n d N}{2, -\delta} + c \, \Norm{(\n - \wn) d N}{2, -\delta}.
\end{eqnarray*}
Considering (\ref{adnbis}) and given that $\n \simeq A + \nz$ ,
\begin{equation}\label{hessn}
\Norm{\n d N}{2,-\delta} \infeg C \, \Norm{N}{2,2, -\delta} \,,
\end{equation}
Using $(\ref{majorhess})$ and $(\ref{hessn})$ , we get
\begin{equation*}
\Norm{g\, \laplag N - \wg\, \wlapla N}{2, \delta} \infeg C \, \Norm{g - \wg}{2,2, \delta} \Norm{N}{2,2, -\delta}.
\end{equation*}
\item For the Ricci tensor:\\
We define $\wade{i}{k}{j} = \wcde{i}{k}{j}-\cdez{i}{k}{j}$.\\
Set $T := \wn - \n = \wg^{-1} \nz \wg - g^{-1} \nz g = (g^{-1} - \wg^{-1}) \nz g + \wg^{-1} \nz (g - \wg) $.\\
Using Hölder inequality (\ref{h1}), Sobolev inclusion ($\delta \infeg 0$) and Sobolev inequality,
\begin{eqnarray}\label{normT}
\Norm{T}{1,2, \delta} 
%& \infeg& \Norm{(g^{-1} - \wg^{-1})}{\infty, \delta} \Norm{\nz g }{2, \delta} + \Norm{g^{-1}}{\infty, 0} \Norm{\nz (g - \wg) }{2, \delta} \nonumber \\
%& \infeg& c \, \Norm{g - \wg}{\infty, \delta}  \Norm{\nz g }{2, \delta} + c \, \Norm{\nz (g - \wg) }{2, \delta} \nonumber\\
& \infeg& C \Norm{g - \wg}{2,2, \delta}.
\end{eqnarray}
We can show, adding and substracting $Ric(\gz)$ and using $(\ref{riccicov})$, that
$$\lbrack Ric (g) - Ric (\wg)\rbrack N \simeq (\nz T + \w{A}T + T^{2}) N \, ,$$
which leads to
\begin{equation}\label{majric}
\Norm{\lbrack Ric (g) - Ric (\wg)\rbrack N}{2, -\delta} \infeg \Norm{\nz T N}{2, -\delta} + \Norm{\w{A}T N }{2, -\delta} + \Norm{T^{2} N }{2, -\delta}.
\end{equation}
Using Hölder and Sobolev inequalities along with Sobolev inclusion ($\delta \infeg 0$) and (\ref{normT})
\begin{eqnarray*}
\Norm{\nz T N}{2, -\delta} 
%& \infeg& \Norm{\nz T}{2, 0} \Norm{N }{\infty, -\delta}\\
%& \infeg& \Norm{\nz T}{2, \delta} \Norm{N }{\infty, -\delta}\\
%& \infeg& c \Norm{T}{1,2, \delta} \Norm{N }{2,2, -\delta}\\
&\infeg& C \, \Norm{g - \wg}{2,2, \delta} \Norm{N }{2,2, -\delta}
\end{eqnarray*}
The same method for the term $\w{A}T N $ gives , considering (\ref{19g})
\begin{eqnarray*}
\Norm{\w{A}T N }{2, -\delta}
% &=& \Norm{\wg^{-1} \nz \wg \, T N}{2, -\delta} \infeg \Norm{\wg^{-1}}{\infty,0} \Norm{\nz \wg \, T N}{2, -\delta}\\
%& \infeg& c \, \Norm{\nz \wg}{1,2, \delta}\Norm{T N}{1,2, -2\delta}\\
%& \infeg& c \, \Norm{\nz \wg}{1,2, \delta}\Norm{T N}{1,2, 0}\\
%& \infeg& C \, \Norm{T}{1,2, \delta}\Norm{N}{1,2, -\delta}\\
&\infeg& C \, \Norm{g - \wg}{2,2, \delta} \Norm{N }{2,2, -\delta}
\end{eqnarray*}
Using $(\ref{h1})$ , $(\ref{19gbis})$ , Sobolev inclusion ($\delta \infeg 0$) and Sobolev inequality together with $(\ref{normT})$ ,
\begin{eqnarray*}
\Norm{T^{2} N }{2, -\delta} 
%&\infeg& \Norm{T^{2}}{2, 0} \Norm{N }{\infty, -\delta}\\
%&\infeg& \Norm{T^{2}}{2, \delta} \Norm{N }{\infty, -\delta}\\
%&\infeg& c \, \Norm{T}{1,2, \delta}^{2}  \Norm{N }{2,2, -\delta}\\
&\infeg& C \, \Norm{g - \wg}{2,2, \delta}^{2} \Norm{N }{2,2, -\delta}.
\end{eqnarray*}
Replacing in (\ref{majric}), we obtain
\begin{equation}\label{majric2}
\Norm{\lbrack Ric (g) - Ric (\wg)\rbrack N}{2, -\delta} \infeg C\, \Norm{g - \wg}{2,2, \delta} \Norm{N}{2,2, -\delta}.
\end{equation}
\item For the scalar curvature, 
\begin{eqnarray*}
(R(g)- 2\Lambda) g\,- (R(\wg)- 2\Lambda) \wg &=& (g - \wg)(R(g)- 2\Lambda) + \wg \wg^{-1}(Ric \, g - Ric \, \wg)\\
&&+ \wg (g^{-1} - \wg^{-1}) Ric \, g\\
&=& (g - \wg) \left\lbrack R(g)- 2\Lambda + n(n-1)\tau^{2} \right\rbrack\\
&&-n(n-1)\tau^{2}(g - \wg) + \wg \wg^{-1}(Ric \, g - Ric \, \wg) \\
&&+ \wg (g^{-1} - \wg^{-1}) \big\{ Ric \, g - Ric \, \gz \big\} .
\end{eqnarray*}
Hölder and Sobolev inequalities as well as Sobolev inclusion ($\delta \infeg 0$) and (\ref{majric2}) yield
\begin{equation*}
\Norm{\lbrack (R(g)- 2\Lambda) g\,- (R(\wg)- 2\Lambda) \wg \rbrack N}{2, -\delta} \infeg C \, \Norm{g - \wg}{2,2, \delta}\Norm{N}{2,2, \delta} \, ,
\end{equation*}
given that $\forall u \in \leb{\infty}{0} , \forall v \in \leb{2}{\delta}$ such that $\Norm{v N}{2, -\delta} \infeg C \, \Norm{N}{2,2, -\delta}$ ,
\begin{eqnarray*}
\Norm{(g - \wg ) \, u \, v \, N}{2, -\delta}  &\infeg& \Norm{g - \wg}{\infty,0} \Norm{u}{\infty,0} \Norm{v N}{2, -\delta}\\
&\infeg& C \, \Norm{g - \wg}{2,2, \delta} \Norm{v}{2,\delta} \Norm{N}{2,2, -\delta} ,
\end{eqnarray*}
where $C$ is a positive constant depending on $\gz , \delta$ and $\Norm{g}{\F}$.
\end{list}
Putting the pieces all together in (\ref{declip}) and taking (\ref{majdphio}) into account lead to
\begin{equation*}
\Norm{\lbrack D\phiogs - D\phio (\tilde{g})^{*}\rbrack N}{2, -\delta} \infeg C \, \Norm{g - \wg}{2,2, \delta}\Norm{N}{2,2, \delta}. \cqfd
\end{equation*}

The dependence in $(g,\pi)$ of $P^*$ is controled as follows :
%La proposition suivante est l'analogue en asymptotiquement hyperbolique de la Proposition 3.4 de \cite{Bartnik2005}.
\begin{propperso}{}\label{proplipcg}
Let $\delta \infeg 0$ , then in dimension 3, the operator $P^{*} : \sob{2}{2}{-\delta}(\T) \longrightarrow \leb{2}{-\delta} \;$ is bounded and satisfies
\begin{equation}\label{estimeePs}
\Norm{\xi}{2,2,-\delta} \infeg c \, \Norm{P^{*} \xi }{2, -\delta} + C \, \Norm{\xi}{1,2,-2\delta} \, ,
\end{equation}
where $C$ depends on $\gz, \delta$ and $\Norm{(g,\pi)}{\F}$.\\
Moreover, $\Pgpis$ depends on $(g,\pi) \in \F$ in a Lipschitz way,
\begin{equation}\label{lipcg}
\Norm{(\Pgpis - \Pgpisw) \, \xi}{2, -\delta} \infeg C_{1} \Norm{(g - \t{g}, \pi - \t{\pi})}{\F} \; \Norm{\xi}{2,2,-\delta}\, ,
\end{equation}
where constant $C_{1}$ depends on $\gz , \delta , \Norm{(g,\pi)}{\F}$ and $\Norm{(\t{g},\t{\pi})}{\F}$.
\end{propperso}
\textbf{Proof}: Let us begin by showing $P^{*}$ is bounded, \ie
\begin{equation}\label{majPs}
\Norm{P^{*} \, \xi}{2, -\delta} \infeg C \; \Norm{\xi}{2,2,-\delta}.
\end{equation}
Set
$$
\begin{cases}
\Ps = \Pgpis\\
D\phisun = D\phigpis_{1}\\
D\phisdeux = D\phigpis_{2}
\end{cases}
.$$
From (\ref{defPs}),
\begin{eqnarray}\label{majPsbis}
\Norm{P^{*} \, \xi}{2, -\delta} &\infeg& c \; (\Norm{D\phisun.\xi}{2,-\delta} + \Norm{\n D\phisdeux.\xi}{2,-\delta}) \nonumber \\
 &\infeg& c \; (\Norm{D\phisun.\xi}{2,-\delta} + \Norm{\nz D\phisdeux.\xi}{2,-\delta} + \Norm{A D\phisdeux.\xi}{2,-\delta}).
\end{eqnarray}
From (\ref{dphis1}),(\ref{h1}), (\ref{19g}), Sobolev inequality and inclusion ($\delta \infeg 0$)
\begin{eqnarray}\label{majdphiun}
\Norm{D\phisun.\xi}{2,-\delta} &\infeg& c \, \big(\Norm{T}{2,-\delta} + \Norm{\Sz (X)}{2,-\delta} + \Norm{N}{2,-\delta} \big) + C \, \big(\Norm{\xi}{\infty,-2\delta} + \Norm{\nz X}{1,2,-2\delta} \big) \nonumber \\
&\infeg& c \, \big(\Norm{N}{2,2,-\delta} + \Norm{X}{1,2,-\delta} + \Norm{N}{2,-\delta} \big) + C \, \big(\Norm{\xi}{2,-\delta} + \Norm{\nz X}{1,2,-\delta} \big) \nonumber\\
&\infeg& C \; \Norm{\xi}{2,2,-\delta}.
\end{eqnarray}
From (\ref{dphis2}) along with (\ref{h1}), (\ref{19g}), Sobolev inequality and inclusion ($\delta \infeg 0$)
\begin{equation}\label{majdphideux}
\Norm{D\phisdeux.\xi}{2,-\delta} \infeg c \, \big(\Norm{\Sz (X)}{2,-\delta} + \Norm{N}{2,-\delta} \big) + C \, \Norm{\xi}{1,2,-2\delta}. 
\end{equation}
\begin{eqnarray*}
\Norm{A D\phisdeux.\xi}{2,-\delta} &\infeg&  c \, \big( \Norm{A \Sz (X)}{2,-\delta} + \Norm{A N}{2,-\delta} \big) + \Norm{\xi}{\infty,-2\delta} \\
&\infeg&  C \, \big( \Norm{\Sz (X)}{1,2,-\delta} + \Norm{N}{1,2,-\delta} \big) + \Norm{\xi}{2,2,-2\delta} \\
%&\infeg&  C \, \Norm{\xi}{2,2,-\delta}. \\[.3cm]
\Norm{\nz D\phisdeux.\xi}{2,-\delta} &\infeg& c \, \big(\Norm{\nz \Sz (X)}{2,-\delta} + \Norm{ \nz N}{2,-\delta} \big) + \Norm{\xi}{\infty,-2\delta}\\
&\infeg&  C \, \Norm{\xi}{2,2,-\delta}.
\end{eqnarray*}
Consequently,
\begin{equation}\label{majdphideux12}
\Norm{D\phisdeux.\xi}{2,-\delta} \infeg \Norm{D\phisdeux.\xi}{1,2,-\delta} \infeg C \, \Norm{\xi}{2,2,-\delta}.
\end{equation}
Every term of (\ref{majPsbis}) is controled by $\Norm{\xi}{2,2,-2\delta}$ leading to (\ref{majPs}).\\
Estimate (\ref{estimeePs}) verified by $P^{*}$ directly comes from (\ref{35cg}). We now look into the Lipschitz behaviour of $P^{*}$:\\
Set
$$
\begin{cases}
\Psw = \Pgpisw\\
%D\phiows = D\phiogws\\
D\phiwsun = D\phigpiws_{1}\\
D\phiwsdeux = D\phigpiws_{2}
\end{cases}
.$$
\begin{eqnarray*}
(\Ps - \Psw) \, \xi &=&
\left(
\begin{array}{c}
g^{-1 \slash 4} \, D\phisun.\xi - \t{g}^{-1 \slash 4} \, D\phiwsun .\xi \\
 g^{1 \slash 4} \, \n D\phisdeux.\xi - \t{g}^{1 \slash 4} \, \nw D\phiwsdeux.\xi 
\end{array}
\right)
=:
\left(
\begin{array}{c}
E \\
 F
\end{array}
\right).
\end{eqnarray*}
So
\begin{equation}\label{majlipEF}
\Norm{(\Ps - \Psw) \, \xi}{2, -\delta} \infeg \Norm{E}{2,-\delta} + \Norm{F}{2,-\delta}.
\end{equation}
\begin{eqnarray*}
E &=& g^{-1 \slash 4} \, D\phisun.\xi - \t{g}^{-1 \slash 4} \, D\phiwsun .\xi \\
&=& (g^{-1 \slash 4} - \t{g}^{-1 \slash 4}) \, D\phisun.\xi + \t{g}^{-1 \slash 4} \, ( D\phisun.\xi -D\phiwsun .\xi ).
\end{eqnarray*}
Using (\ref{h1}), Sobolev inequality and inclusion ($\delta \infeg 0$)
\begin{eqnarray*}
\Norm{E}{2,-\delta} &\infeg& \Norm{(g^{-1 \slash 4} - \t{g}^{-1 \slash 4}) \, D\phisun.\xi}{2,-\delta} + \Norm{\t{g}^{-1 \slash 4} \, ( D\phisun.\xi -D\phiwsun .\xi )}{2,-\delta}\\
%&\infeg& \Norm{g^{-1 \slash 4} - \t{g}^{-1 \slash 4}}{\infty,\delta}\Norm{D\phisun.\xi}{2,-2\delta} + \Norm{\t{g}^{-1 \slash 4}}{\infty,0}\Norm{D\phisun.\xi -D\phiwsun .\xi}{2,-\delta}\\
&\infeg& c \, \Norm{g - \t{g}}{\F}\Norm{D\phisun.\xi}{2,-\delta} + c \, \Norm{D\phisun.\xi -D\phiwsun .\xi}{2,-\delta}.
\end{eqnarray*}
From (\ref{expdphi1dphi0}),
\begin{eqnarray*}
D\phisun.\xi -D\phiwsun .\xi &=& \lbrack D \phibf (g,0)^{*} - D \phibf (\wg,0)^{*} \rbrack \, (N,0) + (\Pi\sqrt{g} - \t{\Pi} \sqrt{\t{g}}) N  +  X \nz (\pi - \t{\pi}) \\
&& + (\pi - \t{\pi}) \nz X + A X (\pi - \t{\pi}).
\end{eqnarray*}
Using (\ref{h1}), (\ref{19g}), Sobolev inequality and inclusion ($\delta \infeg 0$)
\begin{eqnarray*}
\Norm{(\pi - \t{\pi}) \nz X}{2,-\delta} + \Norm{X \nz (\pi - \t{\pi})}{2,-\delta} &\infeg& \Norm{\pi - \t{\pi}}{1,2,\delta}\Norm{\nz X}{1,2,-2\delta}\\
&& + \Norm{\nz (\pi - \t{\pi})}{2,\delta} \Norm{X}{\infty,-2\delta}\\
&\infeg& c \, \Norm{\pi - \t{\pi}}{1,2,\delta} \; \Norm{X}{2,2,-\delta}.
\end{eqnarray*}
\begin{eqnarray*}
\Norm{A X (\pi - \t{\pi})}{2,-\delta} &\infeg& \Norm{A (\pi - \t{\pi})}{2,\delta} \Norm{X}{\infty,-2\delta}\\
%&\infeg& c \, \Norm{A}{1,2,\delta} \Norm{\pi - \t{\pi}}{1,2,\delta} \; \Norm{X}{2,2,-\delta}\\
&\infeg& C \, \Norm{\pi - \t{\pi}}{1,2,\delta} \; \Norm{X}{2,2,-\delta}.
\end{eqnarray*}
\begin{eqnarray*}
\Pi\sqrt{g} - \t{\Pi} \sqrt{\t{g}}&\thicksim& \frac{1}{\sqrt{g}}  g^{-1} \pi^{2}  - \frac{1}{\sqrt{\t{g}}} \t{g}^{-1} \t{\pi}^{2}\\
&\thicksim& \frac{1}{\sqrt{g}}  (g^{-1} -  \t{g}^{-1} ) \pi^{2} + \frac{1}{\sqrt{g}} \t{g}^{-1} (\pi^{2} - \t{\pi}^{2}) + (\frac{1}{\sqrt{g}} - \frac{1}{\sqrt{\t{g}}}) \t{g}^{-1} \t{\pi}^{2}.
\end{eqnarray*}
\begin{eqnarray*}
\Norm{(\Pi\sqrt{g} - \t{\Pi} \sqrt{\t{g}})N}{2,-\delta}
% &\infeg& c \, \Norm{g^{-1} -  \t{g}^{-1}}{\infty,\delta} \big(\Norm{(\pi-\piz)^{2}}{2,\delta} \Norm{N}{\infty,-3\delta} + \Norm{\piz^{2}}{\infty,0} \Norm{N}{\infty,-2\delta} \big)\\
%&& + c \, \Norm{\pi-\t{\pi}}{1,2,\delta}(\Norm{\pi-\piz}{1,2,\delta} + \Norm{\t{\pi}-\piz}{1,2,\delta})\Norm{N}{\infty,-3\delta}\\
%&& + c \, \Norm{\pi-\t{\pi}}{2,\delta}\Norm{\piz^{2}}{\infty,0} \Norm{N}{\infty,-2\delta}\\
%&& + c \, \BNorm{\frac{1}{\sqrt{g}} - \frac{1}{\sqrt{\t{g}}}}{\infty,\delta} \big(\Norm{(\t{\pi}-\piz)^{2}}{2,\delta}\Norm{N}{\infty,-3\delta} + \Norm{\piz^{2}}{\infty,0}\Norm{N}{\infty,-2\delta} \big)\\
%&\infeg& c \, \Norm{g -  \t{g}}{2,2,\delta} \big(\Norm{\pi-\piz}{1,2,\delta}^{2} + 1 \big) \Norm{N}{\infty,-\delta}\\
%&&  + c \, \big(\Norm{\pi-\piz}{1,2,\delta} + \Norm{\t{\pi}-\piz}{1,2,\delta} + 1 \big)\Norm{\pi - \t{\pi}}{1,2,\delta}\Norm{N}{\infty,-\delta}\\
%&& + c \, \BNorm{\frac{1}{\sqrt{g}} - \frac{1}{\sqrt{\t{g}}}}{2,2,\delta} \big(\Norm{\t{\pi}-\piz}{1,2,\delta}^{2} + 1 \big)\Norm{N}{\infty,-\delta} \\
&\infeg& C \, \Norm{(g - \t{g}, \pi - \t{\pi})}{\F} \Norm{N}{2,2,-\delta}.
\end{eqnarray*}
Given (\ref{lip}),
\begin{equation*}
\Norm{D\phisun.\xi -D\phiwsun .\xi}{2,-\delta} \infeg C \, \Norm{(g - \t{g}, \pi - \t{\pi})}{\F} \Norm{\xi}{2,2,-\delta}
\end{equation*}
and taking (\ref{majdphiun}) into account,
\begin{equation}\label{majE}
\Norm{E}{2,-\delta} \infeg C \, \Norm{(g - \t{g}, \pi - \t{\pi})}{\F} \Norm{\xi}{2,2,-\delta}.
\end{equation}
\begin{eqnarray*}
F &=& g^{1 \slash 4} \, \n D\phisdeux.\xi - \t{g}^{1 \slash 4} \, \nw D\phiwsdeux.\xi \\
&=& g^{1 \slash 4} \, (\n - \nw) D\phisdeux.\xi + (g^{1 \slash 4} - \t{g}^{1 \slash 4}) \, \nw D\phiwsdeux.\xi +  g^{1 \slash 4} \, \nw (D\phisdeux.\xi - D\phiwsdeux.\xi).
\end{eqnarray*}
Using (\ref{hess}),(\ref{h1}), (\ref{19g}), Sobolev inequality and inclusion ($\delta \infeg 0$)
\begin{eqnarray*}
\Norm{F}{2,-\delta} &\infeg& c \, \Norm{\n - \nw}{1,2,\delta}\Norm{D\phisdeux.\xi}{1,2,-2\delta} + \Norm{g^{1 \slash 4} - \t{g}^{1 \slash 4}}{\infty,-\delta}\Norm{ \nw D\phiwsdeux.\xi}{2,-2\delta}\\
&&+ c \, \Norm{\nz (D\phisdeux.\xi - D\phiwsdeux.\xi)}{2,-\delta} + c \, \Norm{A (D\phisdeux.\xi - D\phiwsdeux.\xi)}{2,-\delta}\\
&\infeg& C \, \Norm{g - \t{g}}{\F} \Norm{D\phisdeux.\xi}{1,2,-2\delta} +  c \, \Norm{g -  \t{g}}{\F}\Norm{ \nw D\phiwsdeux.\xi}{2,-2\delta}\\
&&+ c \, \Norm{\nz (D\phisdeux.\xi - D\phiwsdeux.\xi)}{2,-\delta} + c \, \Norm{A (D\phisdeux.\xi - D\phiwsdeux.\xi)}{2,-\delta}.
\end{eqnarray*}
Considering (\ref{majdphideux}) and (\ref{majdphideux12}),
\begin{eqnarray}\label{majFbis}
\Norm{F}{2,-\delta} &\infeg& C \, \Norm{g - \t{g}}{\F} \Norm{\xi}{2,2,-\delta} + c \, \Norm{\nz (D\phisdeux.\xi - D\phiwsdeux.\xi)}{2,-\delta}\\
&& + c \, \Norm{A (D\phisdeux.\xi - D\phiwsdeux.\xi)}{2,-\delta}. \nonumber
\end{eqnarray}
\begin{eqnarray*}
D\phisdeux.\xi - D\phiwsdeux.\xi &\thicksim& (K - \t{K}) N + ( A - \t{A})X \\
&\thicksim& (\pi - \t{\pi}) N + ( \n - \nw)X .
\end{eqnarray*}
Using (\ref{hess}),(\ref{h1}), (\ref{19g}), Sobolev inequality and inclusion ($\delta \infeg 0$)
\begin{eqnarray*}
\Norm{\nz (D\phisdeux.\xi - D\phiwsdeux.\xi)}{2,-\delta} &\infeg& c \, \Norm{\nz (\pi - \t{\pi})}{2,\delta}\Norm{N}{\infty,-2\delta} +  c \, \Norm{\pi - \t{\pi}}{1,2,\delta}\Norm{\nz N}{1,2,-2\delta}\\
&& + \Norm{\nz(\n - \nw)}{2,\delta}\Norm{X}{\infty,-2\delta} + \Norm{\n - \nw}{1,2,\delta}\Norm{\nz X}{1,2,-2\delta}\\
&\infeg&  c \, \Norm{\pi - \t{\pi}}{1,2,\delta}\Norm{N}{2,2,-\delta} + c \, \Norm{\n - \nw}{1,2,\delta}\Norm{X}{2,2,-\delta} \\
&\infeg& C \, \Norm{(g - \t{g},\pi - \t{\pi})}{\F} \Norm{\xi}{2,2,-\delta}.
\end{eqnarray*}
In the same way,
\begin{eqnarray*}
\Norm{A (D\phisdeux.\xi - D\phiwsdeux.\xi)}{2,-\delta} &\infeg& c \, \Norm{A(\pi - \t{\pi})}{2,\delta}\Norm{N}{\infty,-2\delta} + \Norm{A(\n - \nw)}{2,\delta}\Norm{X}{\infty,-2\delta}\\
%&\infeg&  c \Norm{A}{1,2,\delta} \Norm{\pi - \t{\pi}}{1,2,\delta}\Norm{N}{2,2,-\delta} + c \, \Norm{A}{1,2,\delta} \Norm{\n - \nw}{1,2,\delta}\Norm{X}{2,2,-\delta} \\
&\infeg& C \, \Norm{(g - \t{g},\pi - \t{\pi})}{\F} \Norm{\xi}{2,2,-\delta}.
\end{eqnarray*}
We deduce from (\ref{majFbis})
\begin{equation}\label{majF}
\Norm{F}{2,-\delta} \infeg C \, \Norm{(g - \t{g}, \pi - \t{\pi})}{\F} \Norm{\xi}{2,2,-\delta}.
\end{equation}
and \corec{the estimate}  (\ref{lipcg}) arises from (\ref{majlipEF}), considering (\ref{majE}) and (\ref{majF}).$\cqfd$\\[.3cm]

We show in the following proposition that the estimate \corec{(\ref{4.22biscg}) of Corollary \ref{corestiadjL2}} is also verified by weak solutions $\xi$ only in $\leb{2}{-\delta} (\T)$. We say that $\xi \in \L$ is a weak solution of $D \phigpis \xi = (f_{1} , f_{2}) \;$ , with $(f_{1} , f_{2}) \in \leb{2}{-\delta}(\w{\S}) \times \sob{1}{2}{-\delta}(\S)$ when
$$\int_{\m} \langle \xi, D \phigpi .(h,p) \rangle_{\gz} = \int_{\m} \langle (f_{1} , f_{2}) , (h,p) \rangle_{\gz} \; , \; \forall (h,p) \in \G \times \K .$$
It suffices to test with $(h,p) \in \cinf_{c} (S \times \w{\S})$ since this place is dense in $\G \times \K$.

\begin{propperso}{}\label{solforte}
Let $\delta \in \rbrack -(n+1) \slash 2 \, , 0 \rbrack \setminus \lbrace -(n-1) \slash 2 \rbrace$ with $n=3$ , $(g,\pi) \in \G^{+} \times \K \, , $\\ $(f_{1} , f_{2}) \in \leb{2}{-\delta}(\w{\S}) \times \sob{1}{2}{-\delta}(\S)$. Let $\xi \in \L$ be a weak solution of $D \phigpis \xi = (f_{1} , f_{2})$.\\
Then $\xi \in \sob{2}{2}{-\delta}(\T)$ is a strong solution and satisfies \corec{(\ref{4.22biscg})}.
\end{propperso}
%\begin{comment}
\textbf{Proof}: In \cite{Bartnik2005}, Bartnik shows that $\xi \in W^{2,2}_{\text{loc}}$. \corec{We can find a cut-off function $\chi_{R}$ as in Definition \ref{cutoff} such that
\begin{itemize}
\item[\textbullet] $\chi_{R} \in \cinf_{c}(\Omega_{R})$.
\item[\textbullet] $\chi_{R} = 1$ on $\Omega_{R \slash 2}$.
\end{itemize}
In particular, $\chi_{R} \xi \in \sob{2}{2}{-\delta}(\T)$ and from Proposition \ref{prop3.3g} , we can write: 
\begin{eqnarray}\label{expchir}
\Norm{\chi_{R}\xi}{2,2,-\delta} %&\infeg&  c \, \big(\Norm{D \phigpis_{1}.(\chi_{R}\xi) }{2, -\delta} + \Norm{D \phigpis_{2}.(\chi_{R}\xi) }{1,2, -\delta} \big) + C \, \Norm{\chi_{R}\xi}{2,-2 \delta} \nonumber\\
&\infeg& c \, \big(\Norm{D \phigpis_{1}.(\chi_{R}\xi) }{2, -\delta} + \Norm{D \phigpis_{2}.(\chi_{R}\xi) }{1,2, -\delta} \big) + C \, \Norm{\xi}{2,- \delta} , \nonumber\\
\end{eqnarray}
using Sobolev inclusion ($\delta \infeg 0$) and considering that $\chi_{R}\xi$ converge to $\xi$ in $\leb{2}{-\delta}$. We have to show that $\chi_{R}\xi$ is uniformly bounded in $\sob{2}{2}{-\delta}$ , $\ie$ bounded independently of $R$. In order to do so, we adapt S. McCormick's method found in \cite{McCormick2015}. From (\ref{expdphi1}) and (\ref{dphis2}),
\begin{eqnarray*}
D \phigpis_{1}.(\chi_{R}\xi) &\simeq& \chi_{R}( D \phigpis_{1}.(\xi) ) + N \, \n^{2} \chi_{R} + dN \, \n \chi_{R} +  X \, \n \chi_{R} + \xi \, \sobd{1} \, \n \chi_{R}. \\
D \phigpis_{2}.(\chi_{R}\xi) &\simeq& \chi_{R}( D \phigpis_{2}.(\xi) ) + X \, \n \chi_{R}.
\end{eqnarray*}
As derivatives of $\chi_{R}$ are supported in $A_{R}:= \Omega_{R} \setminus \Omega_{R \slash 2}$ , we can use (\ref{h1}), (\ref{uinfd}), (\ref{u3d}), Sobolev inclusion ($\delta \infeg 0$) and Ehrling inequality to obtain
\begin{eqnarray}\label{expdphi1chi}
\Norm{D \phigpis_{1}.(\chi_{R}\xi) }{2, -\delta} &\infeg& c \, \big( \Norm{\chi_{R}( D \phigpis_{1}.(\xi) )}{2,-\delta} + \Norm{ N \, \n^{2} \chi_{R}}{2,-\delta} \nonumber\\
&&+ \Norm{dN \, \n \chi_{R}}{2,-\delta} + \Norm{X \, \n \chi_{R}}{2,-\delta} + \Norm{\xi \, \sobd{1} \, \n \chi_{R}}{2,-\delta} \big) \nonumber\\
%&\infeg& c \, \big( \Norm{ D \phigpis_{1}.(\xi)}{2,-\delta} + \Norm{ N}{2,-\delta; A_{R}} \nonumber \\
%&&+ \Norm{dN}{2,-\delta; A_{R}} + \Norm{X}{2,-\delta; A_{R}} + \Norm{\xi \, \sobd{1}}{2,-\delta; A_{R}} \big) \nonumber\\
&\infeg& c \, \Norm{ D \phigpis_{1}.(\xi)}{2,-\delta} + C \, \Norm{\xi}{2,-\delta; A_{R}} +  \epsilon \Norm{\nz^{2}\xi}{2,-\delta; A_{R}} ,\nonumber\\
\end{eqnarray}
where $C$ is a constant depending on $\gz, \lambda, \epsilon, \delta$ and $\Norm{(g,\pi)}{\F}$.\\
Likewise,
\begin{eqnarray}\label{expdphi2chi}
\Norm{D \phigpis_{2}.(\chi_{R}\xi) }{1,2, -\delta} &\infeg& c \, \big( \Norm{\chi_{R}( D \phigpis_{2}.(\xi) )}{1,2,-\delta} + \Norm{X \, \n \chi_{R}}{1,2,-\delta} \big) \nonumber\\
%&\infeg& c \, \big( \Norm{ D \phigpis_{2}.(\xi)}{1,2,-\delta} + \Norm{ X}{1,2,-\delta; A_{R}} \big) \nonumber\\
&\infeg& c \, \Norm{ D \phigpis_{2}.(\xi)}{1,2,-\delta} + C \, \Norm{X}{2,-\delta; A_{R}} +  \epsilon \Norm{\nz^{2}X}{2,-\delta; A_{R}}. \nonumber\\
\end{eqnarray}
Plugging (\ref{expdphi1chi}) and (\ref{expdphi2chi}) into (\ref{expchir}) yields
\begin{eqnarray}\label{expchirbis}
\Norm{\chi_{R}\xi}{2,2,-\delta} &\infeg& c \, \Norm{ D \phigpis_{1}.(\xi)}{2,-\delta} + \Norm{ D \phigpis_{2}.(\xi)}{1,2,-\delta} + C \, \Norm{\xi}{2,-\delta} \nonumber \\
&& +  \epsilon \Norm{\nz^{2}\xi}{2,-\delta; A_{R}}.
\end{eqnarray}

Equation (\ref{defO}) allows us to determine the following link between operators $O$ and $T$:
$$O = T - g \,\corec{\tr}_{g}T - (n-1)gN$$
and we deduce from (\ref{T-gtrT})
\begin{eqnarray}\label{est0}
O &=& D\phigpis_{1}.\xi \slash \sqrt{g} + (n-1)\tau \big(D\phigpis_{2}.\xi - \tdemi \gz\, \corec{\tr}_{\gz} D\phigpis_{2}.\xi \big) \nonumber\\
&& + \lebd \xi +  \sob{1}{2}{\delta} \nz X - (n-1)gN .\label{estO}
\end{eqnarray}
Since $O(N) = L (\n^{2}N)$ with $L$ a linear invertible operator, we also get
\begin{equation*}
\Norm{\n^{2}N}{2, -\delta} \infeg c \, \Norm{O}{2,-\delta} ,
\end{equation*}
leading to
\begin{eqnarray}\label{oplinO}
\Norm{\nz^{2}N}{2, -\delta} &\infeg& c \, \left( \Norm{O}{2,-\delta} + \Norm{A\, d N}{2, -\delta} \right).
\end{eqnarray}
From (\ref{oplinO}), (\ref{est0}) and using (\ref{h1}), (\ref{uinfd}), (\ref{u3d}), Sobolev inclusion ($\delta \infeg 0$) and Ehrling inequality, there exists a constant $C$ depending on $\gz, \lambda, \epsilon, \delta$ and $\Norm{(g,\pi)}{\F}$ such that
\begin{eqnarray}\label{estN2}
\Norm{\nz^{2}N}{2, -\delta} &\infeg& c \, \big( \Norm{D\phigpis_{1}.\xi}{2,-\delta} + \Norm{D\phigpis_{2}.\xi}{2,-\delta} \big) \nonumber\\
&& + \epsilon \; \Norm{\nz^{2}\xi}{2,-\delta} + C \, \Norm{\xi}{2,-\delta}.
\end{eqnarray}
Moreover, combining (\ref{b29ter}) and (\ref{estsx}) and using Sobolev inclusion ($\delta \infeg 0$) and Ehrling inequality, we get
\begin{eqnarray}\label{estX2}
\Norm{\nz^{2}X}{2,-\delta} \infeg c \, \Norm{D\phigpis_{2}.\xi}{1,2,-\delta} + \epsilon \Norm{\nz^{2}\xi}{2,-\delta} + C \, \Norm{\xi}{2,-\delta}.
\end{eqnarray}
Now combination of (\ref{estN2}) and (\ref{estX2}) gives
\begin{equation}\label{estxi2}
\Norm{\nz^{2}\xi}{2,-\delta} \infeg  c \, \big( \Norm{D\phigpis_{1}.\xi}{2,-\delta} + \Norm{D\phigpis_{2}.\xi}{1,2,-\delta} \big) + C \, \Norm{\xi}{2,-\delta}.
\end{equation}
Given that $\xi \in W^{2,2}_{- \delta}(A_{R})$ and that all the inequalities used to obtain (\ref{estxi2}) are valid in particular on every compact set, we have the following local estimate on $ A_{R}$
\begin{eqnarray}\label{estxi2loc}
\Norm{\nz^{2}\xi}{2,-\delta; A_{R}} &\infeg& c \, \big( \Norm{D\phigpis_{1}.\xi}{2,-\delta; A_{R}} + \Norm{D\phigpis_{2}.\xi}{1,2,-\delta; A_{R}} \big) + C \, \Norm{\xi}{2,-\delta; A_{R}} \nonumber \\
&\infeg& c \, \big( \Norm{D\phigpis_{1}.\xi}{2,-\delta} + \Norm{D\phigpis_{2}.\xi}{1,2,-\delta} \big) + C \, \Norm{\xi}{2,-\delta}.
\end{eqnarray}
Finally, equation (\ref{expchirbis}) leads to
\begin{eqnarray*}\label{expchirter}
\Norm{\chi_{R}\xi}{2,2,-\delta} &\infeg& c \, \big(\Norm{ D \phigpis_{1}.(\xi)}{2,-\delta} + \Norm{ D \phigpis_{2}.(\xi)}{1,2,-\delta} \big) + C \, \Norm{\xi}{2,-\delta}\\
 &\infeg& c \, \big(\Norm{f_{1}}{2, -\delta} + \Norm{f_{2}}{1,2, -\delta} \big) + C \, \Norm{\xi}{2, -\delta}.
\end{eqnarray*}
$\chi_{R} \xi$ is then a Cauchy sequence in $\sob{2}{2}{-\delta}$, and so converge in $\sob{2}{2}{-\delta}$. As $\chi_{R} \xi$ converge to $\xi$ in $\leb{2}{-\delta}$ , uniqueness of the limit implies that $\chi_{R} \xi$ converge to $\xi$ in $\sob{2}{2}{-\delta}(\T)$. So $\xi \in \sob{2}{2}{-\delta}(\T)$ and consequently $\xi$ verifies (\ref{4.22biscg}).}$\cqfd$
%\end{comment}

\section{The operator $\Uz$}
Here we introduce an operator $\Uz$ inspired by the formula (\ref{b29}). It will allow us to control the $W^{2,2}_{\delta}$-norm of $X$ with the $L^2_\delta$-norms of $\Sz$ and $\Uz$ , \corec{in other words} with the $W^{1,2}_{\delta}$-norm of $\Sz$. The key estimate will arise from a succession of lemmas.\\
Let $\Uz$ be the operator defined on 1-forms by
%$$\Uz(X) = \nz^{2}X - Riem_{cc-1}\gz \,X \, ,$$
%où $\;Riem_{cc-1}\gz \;$ désigne le tenseur de Riemann de l'espace hyperbolique (à courbure constant -1).\\
%En coordonnées,
\begin{equation}\label{defUz}
\Uz_{kji}(X) = \nz^{2}_{kj}X_{i} - \gz_{jk} X_{i} + \gz_{ik} X_{j}.
\end{equation}
This readily implies
\begin{equation}\label{minUz}
\Norm{\nz^{2} X}{2,-\delta} - c\, \Norm{X}{2,-\delta} \infeg \Norm{\Uz(X)}{2,-\delta}.
\end{equation}
The next four lemmas are established on an asymptotically hyperbolic manifold $(M, \gz)$ , with $\gz = \rho^{-2} \hz$ and $ \norm{d\rho}{\hz}^{2} = 1 + o(1)$ near the boundary \corec{at infinity}.
\begin{lemperso}{}\label{lemihp0}
Let $(M, \gz)$ be an asymptotically hyperbolic manifold and $X \in \cinf(T^{*}\m)$ \\ compactly supported \corec{on $\m$}. $\forall \delta \in \R$ ,
\begin{eqnarray}
\int_{\m} \nz X(\frac{d\rho}{\rho},X) \rho^{2\delta} \, d\mu(\gz) &=& \int_{\m} \rho^{2\delta} \big(\tfrac{n-1}{2} - \delta +o(1)\big)\norm{X}{\gz}^{2} \, d\mu(\gz) \nonumber\\
&&+ \demi \int_{\d\m} \norm{X}{\gz}^{2} \langle \frac{d\rho}{\rho} , \eta \rangle_{\gz} \, \rho^{2\delta} \, d\sigma(\gz) \label{ihp0}\\
\int_{\m} \nz X(X,\frac{d\rho}{\rho}) \rho^{2\delta} \, d\mu(\gz) &=& -\int_{\m}  \rho^{2\delta} \, \Big(\corec{ \corec{\div }}X \langle X, \frac{d\rho}{\rho} \rangle + (2 \delta +1)\langle X, \frac{d\rho}{\rho} \rangle^{2} \Big) \, d\mu(\gz) \nonumber\\
&& +\int_{\m}  \rho^{2\delta}  \big(1 + o(1)\big)\norm{X}{\gz}^{2} \, d\mu(\gz) \nonumber\\
&&+ \int_{\d\m} \rho^{2\delta} \, \langle X, \frac{d\rho}{\rho} \rangle \langle X , \eta \rangle_{\gz} \,  d\sigma(\gz) \, , \label{ihp-1}
\end{eqnarray}
\end{lemperso}
\textbf{Proof}: For (\ref{ihp0}), we integrate by parts the term $\nz_{i}(\norm{X}{\gz}^{2} \, \nz^{i}\rho \, \rho^{2\delta-1})$ and the result follows from the Divergence theorem , along with the definition (\ref{defUz}) of $\Uz$ and (\ref{hessrhoah}). For \corec{(\ref{ihp-1})}, we integrate by parts the term $\nz_{i}\big( X^{i} \, \langle X, d\rho \rangle \, \rho^{2\delta -1}\big)$ and the result follows on from the Divergence theorem , along with the definition (\ref{defUz}). $\cqfd$ \\[.3cm]

%Proofs of the next two lemmas use Lemma \ref{lemihp0} in addition:
\begin{lemperso}{}\label{lemihp1}
Let $(M, \gz)$ be an asymptotically hyperbolic manifold and $X \in \cinf(T^{*}\m)$ compactly supported \corec{on $\m$}. $\forall \delta \in \R$ ,
\begin{eqnarray}\label{ihp1}
\int_{\m} \Uz_{kji}(X) \gz^{kj} X^{i} \rho^{2\delta} \, d\mu(\gz) &=& - \int_{\m}  \rho^{2\delta} \, \norm{\nz X}{\gz}^{2}  \, d\mu(\gz)\nonumber\\
&& + \int_{\m}  \rho^{2\delta} \, \lbrack 2 \delta^{2} - \delta (n-1) -( n-1 )+o(1) \rbrack \norm{X}{\gz}^{2}  \, d\mu(\gz)  \nonumber\\
&&+ \int_{\d\m}  \rho^{2\delta} \, \big( \nz X(\eta,X)-\delta \norm{X}{\gz}^{2}\langle \frac{d\rho}{\rho} , \eta \rangle_{\gz} \big) \, d\sigma(\gz).
\end{eqnarray}
\end{lemperso}
\textbf{Proof}: From (\ref{defUz}), $\; \Uz_{kji}(X) \,\gz^{kj} X^{i} = \langle X, \laplaz X \rangle - (n-1)\norm{X}{\gz}^{2}$.\\
We integrate the term $\nz_{k}(\nz^{k}X_{i} X^{i}\rho^{2\delta})$ and the lemma stems from the Divergence theorem and Lemma \ref{lemihp0}. $\cqfd$

\begin{lemperso}{}\label{lemihp2}
Let $(M, \gz)$ be an asymptotically hyperbolic manifold and $N \in \cinf(\m)$\\ compactly supported \corec{on $\m$}. $\forall \delta \in \R$ ,
\begin{eqnarray}\label{ihp2}
&& \int_{\m} \rho^{2\delta} \,  \big(\Uz_{kji}(X) \nz^{j}X^{i} \frac{\nz^{k}\rho}{\rho} +\corec{ \corec{\div }}X\langle X, \frac{d\rho}{\rho} \rangle \big) \, d\mu(\gz) = \nonumber \\
&&\int_{\m}  \rho^{2\delta} \, \big(\tfrac{n-1}{2} - \delta + o(1)\big) \lbrack  \norm{\nz X}{\gz}^{2} - \norm{X}{\gz}^{2} \rbrack \, d\mu(\gz) - \int_{\m}  \rho^{2\delta} \, \lbrack (2 \delta +1)\langle X, \frac{d\rho}{\rho} \rangle^{2} -\norm{X}{\gz}^{2} \rbrack  \, d\mu(\gz) \nonumber\\
&&+ \demi \int_{\d\m}  \rho^{2\delta} \, \big( \norm{\nz X}{\gz}^{2}\langle \frac{d\rho}{\rho} , \eta \rangle_{\gz} - \norm{X}{\gz}^{2}\langle \frac{d\rho}{\rho} , \eta \rangle_{\gz} \big) \, d\sigma(\gz) + \int_{\d\m}  \rho^{2\delta} \, \langle X, \frac{d\rho}{\rho} \rangle \langle X, \eta \rangle \, d\sigma(\gz).\nonumber\\
&&
\end{eqnarray}
\end{lemperso}
\textbf{Proof}: From (\ref{defUz}),
$$\Uz_{kji}(X) \nz^{j}X^{i} \frac{\nz^{k}\rho}{\rho} = \nz^{2}_{kj}X_{i}\nz^{j}X^{i} \frac{\nz^{k}\rho}{\rho} - \nz X(\frac{d\rho}{\rho},X) + \nz X(X,\frac{d\rho}{\rho}) .$$
We integrate the term $\nz_{k}(\norm{\nz X}{\gz}^{2} \frac{\nz^{k}\rho}{\rho}\rho^{2\delta})$ and the lemma stems from the Divergence theorem and Lemma \ref{lemihp0}. $\cqfd$\\

Considering also the equality $2 \, \norm{\Sz (X)}{\gz}^{2} = \norm{\nz X}{\gz}^{2} + \nz_{j}X_{k}\nz^{k}X^{j} \, ,$ we get
\begin{lemperso}{}\label{lemihp3}
Let $(M, \gz)$ be an asymptotically hyperbolic manifold and $X \in \cinf(T^{*}\m)$ compactly supported \corec{on $\m$}. $\forall \delta \in \R$ ,
\begin{eqnarray}\label{ihp3}
&&\int_{\m} \rho^{2\delta} \big( \Uz_{kji}(X) \gz^{ik} X^{j} + 2 \, \norm{\Sz (X)}{\gz}^{2} - 2 \delta \,\corec{ \corec{\div }}X\langle X, \frac{d\rho}{\rho} \rangle \big) \, d\mu(\gz)= \nonumber \\
&& \int_{\m}  \rho^{2\delta} \, \lbrack n-1 -  2 \delta +o(1) \rbrack \norm{X}{\gz}^{2}  \, d\mu(\gz) + \int_{\m}  \rho^{2\delta} \, \norm{\nz X}{\gz}^{2}  \, d\mu(\gz) \\
&&+ \int_{\m}  \rho^{2\delta} \, 2 \delta (2 \delta + 1) \langle X, \frac{d\rho}{\rho} \rangle_{\gz}^{2} \, d\mu(\gz) + \int_{\d\m}  \rho^{2\delta} \, \big( \nz X(X,\eta)-2 \delta \, \langle X, \frac{d\rho}{\rho} \rangle_{\gz} \langle X, \eta \rangle_{\gz} \big) \, d\sigma(\gz) \, , \nonumber
\end{eqnarray}
\end{lemperso}
\textbf{Proof}: From (\ref{defUz}), $\; \Uz_{kji}(X) \,\gz^{ik} X^{j} = \nz^{2}_{kj}X^{k}X^{j} + (n-1)\norm{X}{\gz}^{2}$.\\
We integrate by parts the term $\nz_{k}(\nz_{j}X^{k}X^{j}\rho^{2\delta})$ and the result follows on from the Divergence theorem along with (\ref{ihp-1}) and the equality $2 \, \norm{\Sz (X)}{\gz}^{2} = \norm{\nz X}{\gz}^{2} + \nz_{j}X_{k}\nz^{k}X^{j}$. $\cqfd$\\[.2cm]
\begin{comment}
\begin{eqnarray}\label{ihp3ter}
\int_{\m} \rho^{2\delta} \nz^{2}_{kj}X^{k}X^{j}  \, d\mu(\gz) &=& \int_{\d\m} \nz X(X,\eta) \rho^{2\delta} \, d\sigma(\gz) \nonumber \\
&-&\int_{\m} \rho^{2\delta} \big( \nz_{j}X_{k}\nz^{k}X^{j} + 2 \delta \, \nz X(X,\frac{d\rho}{\rho}) \big) \, d\sigma(\gz) .
\end{eqnarray}
L'égalité
$$2 \, \norm{\Sz (X)}{\gz}^{2} = \norm{\nz X}{\gz}^{2} + \nz_{j}X_{k}\nz^{k}X^{j}$$
nous permet de remplacer le terme $\nz_{j}X_{k}\nz^{k}X^{j}$ dans (\ref{ihp3ter}) afin d'obtenir (\ref{ihp3}) compte tenu de (\ref{ihp-1}).$\cqfd$
\end{comment}
We can now prove the following proposition, crucial in the demonstration of the adjoint kernel triviality.
\begin{propperso}{}\label{x2infcUzinf}
Set $E_{R} := \m \setminus \Omega_{R}$ , for an open set $\Omega_{R}$.
For all $\epsilon >0$ , for all $\delta \in \rbrack -2 \, , -1 \lbrack\,$ , there exists $R_{\epsilon,\delta}>0 $ such that for all $R>R_{\epsilon,\delta}$ , there exists a constant $c>0$ such that
\begin{equation}\label{x2infcUzbord}
\forall X \in \cinf_{c}(E_{R}) \; \; , \; \;\Norm{X}{1,2,-\delta; E_{R}} \infeg c \, (\Norm{\Uz(X)}{2,-\delta; E_{R}} + \Norm{\Sz(X)}{2,-\delta; E_{R}}).
\end{equation}
\end{propperso}
\textbf{Proof}: The linear combination $(\ref{ihp2}) - \tdemi (\ref{ihp1}) + \tdemi (\ref{ihp3}) + (\ref{eqcord10}) - \tdemi(\ref{eqcord11}) $ yields
\begin{eqnarray}\label{combigagnante1}
&&\int_{\m} \rho^{2\delta} \,  \Uz_{kji}(X) \Big(\nz^{j}X^{i} \frac{\nz^{k}\rho}{\rho} - \tdemi \gz^{kj} X^{i} + \tdemi \gz^{ik} X^{j} \Big) \, d\mu(\gz)  \nonumber\\
&&+ \int_{\m} \rho^{2\delta} \, \Big( \norm{\Sz (X)}{\gz}^{2} - \Sz(X)(\frac{d\rho}{\rho},\frac{d\rho}{\rho}) \; \langle X,\frac{d\rho}{\rho} \rangle_{\gz} + (2-\delta)\corec{\div }X\langle X, \frac{d\rho}{\rho} \rangle \Big) \,d\mu(\gz) \nonumber\\
 &=& \int_{\m} \lbrace - \delta^{2} + (\tfrac{n-3}{2}) \delta + n+1 +o(1) \rbrace \norm{X}{\gz}^{2} \rho^{2\delta} \, d\mu(\gz) \nonumber \\
&& + \int_{\m} \rho^{2\delta} \, \left\lbrace 2 \delta^{2} - 2\delta -(\tfrac{n+3}{2}) +o(1) \right\rbrace \langle X,\frac{d\rho}{\rho} \rangle_{\gz}^{2} \, d\mu(\gz)  \nonumber \\
&&  + \int_{\m} \left\lbrack \tfrac{n+1}{2} -\delta+o(1) \right\rbrack \norm{\nz X}{\gz}^{2} \rho^{2\delta} \, d\mu(\gz) + \int_{\d\m} \rho^{2\delta} \, \tdemi (\delta - 1) \norm{ X}{\gz}^{2}\langle \frac{d\rho}{\rho} , \eta \rangle_{\gz} \, d\sigma(\gz) \nonumber\\
&& + \int_{\d\m} \Big( \tdemi \norm{\nz X}{\gz}^{2}\langle \frac{d\rho}{\rho} , \eta \rangle_{\gz} + \tdemi \norm{ X}{\gz}^{2}\langle \frac{d\rho}{\rho} , \eta \rangle_{\gz} + \tdemi \nz X(X,\frac{d\rho}{\rho}) - \tdemi \nz X(\frac{d\rho}{\rho},X) \Big)  \, d\sigma(\gz)\nonumber \\
&& + \int_{\d\m} \rho^{2\delta} \,\Big((2- \delta)\langle X, \frac{d\rho}{\rho} \rangle_{\gz} \langle X, \eta \rangle_{\gz} -\tdemi \langle X, \frac{d\rho}{\rho} \rangle_{\gz}^{2} \langle \frac{d\rho}{\rho}, \eta \rangle_{\gz} \Big) \, d\sigma(\gz) .
\end{eqnarray}
Application on $E_{R}$:
\corec{With the same notations as in Proposition \ref{n2inftzinf}, since $X \in \cinf_{c}(E_{R})$, boundary terms in (\ref{7}) will only concern $\d \Omega_{R}$}. $X_{n}$ (resp. $X_{T}$) being the component of $X$ normal (resp. tangential) to $\corec{\d_{\infty}\m}$, $X_{n} := \langle X , \eta \rangle_{\gz} \; \; \; \text{and} \; \; \; \norm{X}{\gz}^{2} = X_{n}^{2}+ X_{T}^{2}$.
\begin{eqnarray}\label{combigagnante2}
&&\int_{E_{R}} \rho^{2\delta} \,  \Uz_{kji}(X) \Big(\nz^{j}X^{i} \frac{\nz^{k}\rho}{\rho} - \tdemi \gz^{kj} X^{i} + \tdemi \gz^{ik} X^{j} \Big) \, d\mu(\gz)  \nonumber\\
&&+ \int_{E_{R}} \rho^{2\delta} \, \big( \norm{\Sz (X)}{\gz}^{2} - \Sz(X)(\eta_{R},\eta_{R}) \; X_{n} + (2- \delta)\corec{\div }X X_{n} \big) \,d\mu(\gz) \nonumber\\
 &=& \int_{E_{R}} \rho^{2\delta} \, \lbrace - \delta^{2} + (\tfrac{n-3}{2}) \delta + n+1 +o(1) \rbrace X_{T}^{2} \, d\mu(\gz)\nonumber \\
&& + \int_{E_{R}} \rho^{2\delta} \, \left\{ \left\lbrack \tfrac{n+1}{2} -\delta +o(1) \right\rbrack \norm{\nz X}{\gz}^{2} + \left\lbrack \delta^{2} + (\tfrac{n-7}{2}) \delta + \tfrac{n-1}{2} +o(1) \right\rbrack X_{n}^{2} \right\} \, d\mu(\gz) \nonumber\\
&&  + \int_{\d E_{R}} \rho^{2\delta} \, \tdemi(\delta - 1) \big(X_{T}^{2} + o(1) \big)  \, d\sigma(\gz) + \int_{\d E_{R}} \rho^{2\delta} \,\tdemi(2- \delta)\big(X_{n}^{2} + o(1) \big) \, d\sigma(\gz)\nonumber \\
&& + \int_{\d E_{R}} \Big( \tdemi \norm{\nz X}{\gz}^{2}\norm{\eta}{\gz}^{2} + \tdemi \norm{ X}{\gz}^{2} + \tdemi \nz X(X,\eta) - \tdemi \nz X(\eta,X) \Big)  \, d\sigma(\gz). \nonumber
\end{eqnarray}
Given the following equalities
\begin{eqnarray*}\label{inegbord}
0 \infeg \tfrac{1}{4}\norm{\nz_{i} X_{j}\eta^{i}- X_{j}}{\gz}^{2} &=& \tfrac{1}{4}\norm{\nz X}{\gz}^{2}\norm{\eta}{\gz}^{2} + \tfrac{1}{4}\norm{X}{\gz}^{2} - \tdemi \nz X(\eta,X)\\
0 \infeg \tfrac{1}{4}\norm{\nz_{i} X_{j}\eta^{j}+ X_{i}}{\gz}^{2} &=& \tfrac{1}{4}\norm{\nz X}{\gz}^{2}\norm{\eta}{\gz}^{2} + \tfrac{1}{4}\norm{X}{\gz}^{2} + \tdemi \nz X(X,\eta)
\end{eqnarray*}
we get that for all $\epsilon >0$ , there exists $R_{\epsilon}>0 $ and $c_{\epsilon} >>1$ such that $ \forall R>R_{\epsilon}$ , 
\begin{eqnarray}\label{combigagnante3}
c_{\epsilon} \, \int_{E_{R}} \rho^{2\delta} \, \big( \norm{\Uz(X)}{\gz}^{2} + \norm{\Sz (X)}{\gz}^{2} \big) \, d\mu(\gz)
 %&\supeg& \int_{E_{R}} \rho^{2\delta} \,\lbrack \tfrac{n+1}{2} -\delta - \epsilon \rbrack \norm{\nz X}{\gz}^{2} \, d\mu(\gz) + \int_{E_{R}} \rho^{2\delta} \, \lbrace - \delta^{2} + (\tfrac{n-3}{2}) \delta + n+1 - \epsilon \rbrace X_{T}^{2} \, d\mu(\gz) \nonumber\\
%&& + \int_{E_{R}} \rho^{2\delta} \, \lbrace \delta^{2} + (\tfrac{n-7}{2}) \delta + \tfrac{n-1}{2} - \epsilon \rbrace X_{n}^{2} \, d\mu(\gz) \nonumber \\
%&&+ \demi \int_{\d E_{R}} \rho^{2\delta} \, \Big( (\delta - 1- \epsilon) X_{T}^{2} + (2- \delta - \epsilon) X_{n}^{2} \Big) \, d\sigma(\gz)\nonumber \\
%&& + \tfrac{1}{4} \int_{\d E_{R}} \Big( \norm{\nz_{i} X_{j}\eta^{i}- X_{j}}{\gz}^{2} + \norm{\nz_{i} X_{j}\eta^{j}+ X_{i}}{\gz}^{2} \Big)  \, d\sigma(\gz) \nonumber \\
 &\supeg& \int_{E_{R}} \rho^{2\delta} \, \lbrace - \delta^{2} + (\tfrac{n-3}{2}) \delta + n+1 - \epsilon \rbrace X_{T}^{2} \, d\mu(\gz) \nonumber\\
&& \int_{E_{R}} \rho^{2\delta} \,\lbrack \tfrac{n+1}{2} -\delta - \epsilon \rbrack \norm{\nz X}{\gz}^{2} \, d\mu(\gz) \nonumber \\
&& + \int_{E_{R}} \rho^{2\delta} \, \lbrace \delta^{2} + (\tfrac{n-7}{2}) \delta + \tfrac{n-1}{2} - \epsilon \rbrace X_{n}^{2} \, d\mu(\gz) \nonumber \\
&&+ \demi \int_{\d E_{R}} \rho^{2\delta} \, (2- \delta - \epsilon) X_{n}^{2} \, d\sigma(\gz)\nonumber\\
&&+ \demi \int_{\d E_{R}} \rho^{2\delta} \, (\delta - 1- \epsilon) X_{T}^{2} \, d\sigma(\gz).
\end{eqnarray}
\begin{itemize}
\item[\textbullet]The $X_{T}^{2}$ term
%, nous avons un trinôme du second degré en $\delta$ dont le discriminant est
%$$\Delta = \Big(\frac{n-3}{2}\Big)^{2} + 4(n+1) = \Big(\frac{n+5}{2}\Big)^{2}>0$$
%et dont les racines sont
%$$\delta_{\pm} = \frac{(n-3) \pm (n+5)}{4} = \Big \lbrace -2 \, ; \frac{n+1}{2} \Big \rbrace.$$
%Ce terme 
is non negative if $\delta \in \rbrack -2 \, ; \frac{n+1}{2} \lbrack$.\\
\item[\textbullet]The $X_{n}^{2}$ term
%, nous avons un trinôme du second degré en $\delta$ dont le discriminant est
%$$\Delta = \Big(\frac{n-7}{2}\Big)^{2} -2(n-1) = \Big(\frac{n^{2}-22n + 57}{4}\Big).$$
%Ce terme 
is non negative $\forall n \supeg 3$ , $ \forall \delta \in \rbrack \frac{n-1}{2} \, ; \frac{n+1}{2} \lbrack$.\\
\item[\textbullet]The $\norm{\nz X}{\gz}^{2}$ term is non negative if $\delta < (n+1) \slash 2$.\\
\item[\textbullet]The boundary term is non negative if $\delta \in \rbrack 1 , 2 \lbrack$.
\end{itemize}
Hence, $\forall \delta \in \rbrack 1 \, ; 2 \lbrack$, $\forall \epsilon >0 , \exists R_{\epsilon}>0 $ and $c_{\epsilon} >>1$ such that $ \forall R>R_{\epsilon}$ , 
\begin{eqnarray*}
&&c_{\epsilon} \, \int_{E_{R}} \rho^{2\delta} \, \big( \norm{\Uz(X)}{\gz}^{2} + \norm{\Sz (X)}{\gz}^{2} \big) \, d\mu(\gz) \supeg \int_{E_{R}} \rho^{2\delta} \,\lbrack \tfrac{n+1}{2} -\delta - \epsilon \rbrack \norm{\nz X}{\gz}^{2} \, d\mu(\gz) \nonumber\\
 &&+ \int_{E_{R}} \rho^{2\delta} \, \Big(\lbrace - \delta^{2} + (\tfrac{n-3}{2}) \delta + n+1 - \epsilon \rbrace X_{T}^{2} + \lbrace \delta^{2} + (\tfrac{n-7}{2}) \delta + \tfrac{n-1}{2} - \epsilon \rbrace X_{n}^{2} \, d\mu(\gz) \nonumber \Big) \, d\mu(\gz). \cqfd
\end{eqnarray*}
%C'est à dire
%$$\Norm{\Uz(X)}{2,-\delta; E_{R}} + \Norm{\Sz(X)}{2,-\delta; E_{R}} \supeg c \, \Norm{X}{1,2,-\delta; E_{R}}.$$ 

\section{The adjoint kernel triviality}

In this section we show that the kernel of $D \phigpis$ is trivial. We will need the following lemmas and propositions in the proof of Theorem \ref{trivialkercg}.

\begin{lemperso}{}\label{lemkertrivcg}
Let $\delta \in \R$ and $\xi \in \sob{2}{2}{-\delta}(\T)$ be a solution of $D \phiogpis \xi = 0$. Then $\xi=(N,X)$ satisfies \corec{a system of the form}
\begin{equation}\label{eqkertrivcg}
\begin{cases}
\Tz (N) = b_{0} \xi + b_{1} \nz \xi\\
\Sz(X) =  - \gz \tau N + b_{2} \xi
\end{cases}
,
\end{equation}

with $b_{0} \in \lebd$ \corec{and} $b_{1} , b_{2} \in \sobd{1}$.\\
\end{lemperso}
\textbf{Proof}: From (\ref{T-gtrT}) and (\ref{dphis2}) , $D \phiogpis \xi = 0$ leading to
\begin{eqnarray}
&&T- (\corec{\tr}_{g}T)g = \lebd \xi +  \sob{1}{2}{\delta} \nz X \label{T-gtrTbis}.\\
&&\Sz(X) =  - \gz \tau N + \sob{1}{2}{\delta} \xi. \nonumber
\end{eqnarray}
Taking the trace of (\ref{T-gtrTbis}), 
\begin{eqnarray}
&&\Tz(N) = \lebd \xi +  \sob{1}{2}{\delta} \nz \xi .\label{T-gtrTter}\\
&&\Sz(X) =  - \gz \tau N + \sob{1}{2}{\delta} \xi. \nonumber \cqfd
\end{eqnarray}

\begin{thmperso}{Triviality of ker \corec{$D \phiogpis$}.\\}\label{trivialkercg}
Let $\Omega \subset \m$ be a connected open set such that $E_{R} \subset \Omega$. We fix $(g, \pi) \in \F$. Set \corec{$\delta \in \rbrack -2\, , -1 \lbrack$} and $n=3$. Suppose $\xi \in \leb{2}{-\delta}(\m)$ verifies $\; D \phiogpis \xi = 0\;\mbox{ on } \Omega.$\\ Then $\xi \equiv 0$ on $\Omega$.
\end{thmperso}
\textbf{Proof}: From Proposition \ref{solforte} , $\xi \in \sob{2}{2}{-\delta}(\m)$. According to the Lemma \ref{lemkertrivcg} and since $\xi$ is a solution of $D \phiogs \xi = 0$ , $\xi$ satisfies
\begin{equation*}
\begin{cases}
\Tz(N) = b_{0} \xi + b_{1} \nz \xi\\
\Sz(X) =  - \gz \tau N + b_{2} \xi
\end{cases}
,
\end{equation*}
with $b_{0} \in \lebd$ \corec{and} $b_{1} , b_{2} \in \sobd{1}$.\\

From (\ref{integriem}),
\begin{equation}\label{defUzbis}
\Uz_{kji}(X) = \nz^{2}_{kj}X_{i} - Riem \, \gz_{ijkl} + \lebd X.
\end{equation}
Regarding (\ref{b29})
\begin{eqnarray}
\Uz(X) &=& c_{1} \nz \Sz (X) + \lebd X \nonumber\\
&=& -c_{1} \tau \gz \nz N + \lebd \xi + \sob{1}{2}{\delta} \nz \xi .\label{Uzbis}
\end{eqnarray}
We must show that a solution $\xi$ of (\ref{eqkertrivcg}) such that $\xi = o(\rho^{\delta})$ (from (\ref{tauxdec})) vanishes. Before pursuing the proof of the theorem, let us recall Proposition  3.9  of \cite{Bartnik2005} :
\begin{propperso}{}\label{prop 3.9cg}
In dimension $n=3$, set $\delta\infeg 0$ and $(g,\pi) \in \F$. Let $\Omega$ be a connected subset of $\m$. Let $\xi$ satisfy $D \phigpis \xi = 0$ \corec{on} $\Omega.$
If in addition $\xi\equiv 0$ on an open set $U \subset \Omega$ , then $\xi \equiv 0$ on $\Omega$.
\end{propperso}
Given the previous proposition, it remains to show $\xi$ vanishes near infinity.
\corec{As for (\ref{minUz}), we have}
\begin{equation}\label{minUzbord}
\Norm{\nz^{2} X}{2,-\delta;E_{R}} - c\, \Norm{X}{2,-\delta;E_{R}} \infeg \Norm{\Uz(X)}{2,-\delta;E_{R}}.
\end{equation}
Combine to (\ref{x2infcUzbord}), together with (\ref{Uzbis}), Sobolev inequality and (\ref{h1}), we obtain
\begin{eqnarray}\label{X2infcUzbordbis}
\Norm{X}{2,2,-\delta; E_{R}} &\infeg& c \, (\Norm{\Uz(X)}{2,-\delta; E_{R}} + \Norm{\Sz(X)}{2,-\delta; E_{R}}) \nonumber\\
&\infeg& c \, \Norm{N}{1,2,-\delta; E_{R}} + C \, (\Norm{\xi}{\infty,-2\delta;E_{R}} + \Norm{\nz \xi}{3,-2\delta;E_{R}}).
\end{eqnarray}
Using Proposition \ref{n2inftzinf} , there exists $\epsilon_{0} <<1$ such that $(\ref{n2infctzbord}) + \epsilon_{0} (\ref{X2infcUzbordbis})$ give
\begin{eqnarray*}%\label{xi2infcTzbord}
\Norm{\xi}{2,2,-\delta; E_{R}} &\infeg& c \, \Norm{\Tz}{2,-\delta; E_{R}} + C \, (\Norm{\xi}{\infty,-2\delta;E_{R}} + \Norm{\nz \xi}{3,-2\delta;E_{R}}).
\end{eqnarray*}
Considering Lemma \ref{lemkertrivcg} along with Sobolev inequality and (\ref{49bis}) with $\delta \infeg 0$,
\begin{eqnarray*}
\Norm{\xi}{2,2,-\delta;E_{R}} &\infeg& C \, \Norm{\xi}{2,2,-2\delta;E_{R}}\\
&\infeg& C \, e^{4R \delta} \, \Norm{\xi}{2,2,-\delta;E_{R}}.
\end{eqnarray*}
We end up with $ \forall \, \delta \in \rbrack -(n+1) \slash 2 \, , -1 \rbrack$,
\begin{equation*}
\Norm{\xi}{2,2,-\delta; E_{R}} \infeg C \, e^{4R \delta} \, \Norm{\xi}{2,2,-\delta;E_{R}}.
\end{equation*}
For $R$ large enough, $\Norm{\xi}{2,2,-\delta;E_{R}} = 0$ and thanks to Sobolev inequality, $\xi$ vanishes on $E_{R}$, for $R>>1$. From Prop.\ref{prop 3.9cg}, $\xi \equiv 0$ on $\Omega$, because $\Omega$ is connected by assumption. This ends the proof of Theorem \ref{trivialkercg}.$\cqfd$\\[.5cm]
%Le corollaire suivant est l'annalogue du Corollaire 3.11 de \cite{Bartnik2005} en asymptotiquement hyperbolique.
\begin{corperso}{}
Set $\delta \in \rbrack -(n+1) \slash 2 \, , 0 \rbrack \setminus \lbrace -(n-1) \slash 2 \rbrace$, with $n=3$. There exists a constant $C>0$ depending on $\Norm{(g,\pi)}{\F}$ such that for $\xi \in \sob{2}{2}{-\delta}(\T)$,
\begin{equation}
\Norm{\xi}{2,2,-\delta} \infeg C \, \Norm{P^{*} \xi}{2, -\delta}.
\end{equation}
\end{corperso}
\textbf{Proof}: In order to show that the kernel of $P^{*}$ is finite dimensional, we apply Riesz theorem showing every bounded subset of ker $P^{*}$ is $\Norm{\;}{2,2,-\delta}-$ compact. Let $\lbrace \xi_{k}\rbrace$ be a sequence of ker $P^{*}$ such that $\Norm{\xi_{k}}{2,2,-\delta} = 1$. Rellich theorem tells us we can extract from $\lbrace \xi_{k}\rbrace$ a sub-sequence, also noted $\lbrace \xi_{k}\rbrace$, converging in $\sob{1}{2}{-2\delta}$ to a limit $\bar{\xi}$. Hence, $\lbrace \xi_{k}\rbrace$ is a Cauchy sequence in $\sob{1}{2}{-2\delta}$. From (\ref{35cg}), considering that $\lbrace \xi_{k}\rbrace \in \text{ker} P^{*}$, $\lbrace \xi_{k}\rbrace$ is a Cauchy sequence in $\sob{2}{2}{-\delta}$, and so converges to $\bar{\xi}$ in $\sob{2}{2}{-\delta}$, from the limit uniqueness. This ends the proof of the finite dimension of ker $P^{*}$.
ker $P^{*}$ is thus a closed subspace of the Hilbert vector space $\sob{2}{2}{-\delta}$. Being a finite dimensional closed subspace of a normed vector space , it splits and if we set $W$ to be the closed complement of ker $P^{*}$,
$$\sob{2}{2}{-\delta} = \text{ker} P^{*} \oplus W .$$
From the same argument as in the proof of Theorem \ref{estellip} , there exists a constant $C>0$ depending on $\Norm{(g,\pi)}{\F}$ such that for all $\xi \in W$ ,
\begin{equation}
\Norm{\xi}{2,2,-\delta} \infeg C \, \Norm{P^{*} \xi }{2, -\delta}
\end{equation}
We conclude thanks to the triviality of ker $P^{*}$ from Theorem \ref{trivialkercg}.$\cqfd$\\

\section{The submanifold structure}

\begin{lemperso}{}\label{lemiso}
Let $X,Y$ be two Banach spaces and $T$ a linear operator with closed range.
\begin{eqnarray*}
T: X &\rightarrow &Y\\
T^{*}: Y^{*}& \rightarrow &X^{*}
\end{eqnarray*}
then $({\Coker} T)^{*} \simeq {\ker}\, T^{*}$,
where ${\Coker} T = \faktor{Y}{\Im T}$
\end{lemperso}
\textbf{Proof}: We define
\begin{eqnarray*}
 \psi \colon {\ker}\; T^{*}&\rightarrow &(\Coker T)^{*} = \mathcal{L} \left(\faktor{Y}{\Im T}, \R \right)\\
\rho & \mapsto &(\lambda \colon y + TX \mapsto \rho(y))
\end{eqnarray*}
The map $\lambda$ is well defined because $\forall x \in X \, , \rho(Tx) = T^{*}(\rho)(x) = 0$.\\
The map $\psi$ is invertible and
\corec{\begin{eqnarray*}
\psi^{-1} \colon (\Coker T)^{*}&\rightarrow &\text{ker}\, T^{*}\\
\lambda & \mapsto & \rho \hspace{1cm}\text{ where } \rho(y) :=\lambda(y + TX).
\end{eqnarray*}
Note that $\rho \in \text{ker}\, T^{*}$} because $T^{*}(\rho)(x) = \rho(Tx) = \lambda(\bar{0}) = 0. \cqfd$\\
Remark: The closed range of $T$ imply ${\Coker} T $ is a Banach space.

\begin{thmperso}{}\label{mainthmmanifold}
Let $\phibf : \F \rightarrow \L^{*}$ be the constraint operator
in dimension $n=3$.\\
For every $\epsilon \in \L^{*}$ , for all \corec{$\delta  \in \rbrack - 2 \, , -1 \lbrack$} , the set of solutions of the constraint equations
$$\contrainte (\epsilon) := \lbrace (g,\pi) \in \F : \phigpi = \epsilon \rbrace$$
is a submanifold of $\F$. In particular,  the space of solutions of the vacuum constraint equations  $\contrainte = \contrainte (0)$ has a Hilbert submanifold structure.
\end{thmperso}
 \corec{In order to prove the Theorem 
\ref{mainthmmanifold}, we will  use the implicit function theorem, so we have to show:}
\begin{list}{$\bullet$}
\item $\ker\; D \phigpi$ splits.
\item \item $D \phigpi$ is surjective.
\end{list}
\medskip
Given that the kernel of $D \phigpi$ is finite dimensional, we show that ker $D \phigpi$ is closed and hence splits. $D \phigpi$ being a bounded operator, its kernel is closed by continuity.\\
\corec{The triviality of ker $D \phigpis$, established in Theorem \ref{trivialkercg}, leads to}
\begin{equation*}
\left( {\ker} D \phigpis \right)^{\perp} = \L^{*}.
\end{equation*}
Using the classic relation
\begin{equation*}
\left( {\ker} D \phigpis \right)^{\perp} = \overline{{\Im} D \phigpi} ,
\end{equation*}
we get
$$\overline{{\Im} D \phigpi} = \L^{*}.$$
Thus, in order to have the surjectivity of $D \phiog$ , it suffices to prove $D \phigpi$ has closed range. To do so, we consider particular variations $(h,p)$ of $(g,\pi)$ of the form 
\begin{equation}\label{defhp}
\begin{cases}
h_{ij} = 2 \, y \, g_{ij}\\
p^{ij} = \big(2 \S(Y)^{ij} - g^{ij} \,\corec{\tr}_{g}\S(Y) -(n-1)(n-2) \tau y \, g^{ij} \big) \sqrt{g}
\end{cases}
,
\end{equation}
determined from fields $(y,Y^{i})$. We define the operator
\begin{equation}
F(y, Y^{i}) = \lbrack F_{0}(y, Y^{i}),F_{i}(y, Y^{i})\rbrack = \lbrack D \phiogpi (h,p),  D \phiigpi (h,p) \rbrack.
\end{equation}
From (\ref{def0}) and (\ref{defi}),
\begin{equation*}
\left\lbrace
\begin{array}{rcl}
 F_{0}(y, Y^{i})& =& 2(n-1) \sqrt{g} \, \lbrack -\Delta y + n y  \rbrack + (4-n) \, \phiogpi \, y + 2(n-2)\tau\,{ \corec{\div }}Y \sqrt{g}\\
&& + \leb{2}{\delta} \, \lbrack y + Y \rbrack + \sob{1}{2}{\delta} \, \nz Y.\\[.2cm]
F_{i}(y, Y^{i}) &=& -2 \sqrt{g} \, \lbrack -\Delta Y^{i} + (n-1) Y^{i}  \rbrack + 2 \, \phiigpi \, y + \sob{1}{2}{\delta} \, \nz y + \leb{2}{\delta} \, \lbrack y + Y \rbrack.
\end{array}
\right.
\end{equation*}

\begin{defperso}{Operator asymptotic to $\laplaz$.\\}\label{defopasympt}
We say an operator $P$ of the form
\begin{equation*}
Pu = a^{ij}(x)\nz^{2}_{ij}u + b^{i}(x)\d_{i}u + c(x)u
\end{equation*}
is asymptotic to $\laplaz$ with a decaying rate $\tau$ if there exists $n<q<\infty$ , $\tau \infeg0$ and two positive constants $C_{1} , \lambda$ such that
\begin{eqnarray*}
&&\lambda \norm{\xi}{\gz}^{2} \infeg a^{ij}(x) \xi_{i}\xi_{j} \infeg \lambda^{-1} \norm{\xi}{\gz}^{2} , \forall x \in \m \; , \; \xi \in T\m. \\
&&\Norm{a^{ij} - \gz^{ij}}{1,q,\tau} + \Norm{b^{i}}{q,\tau} + \Norm{c}{q\slash 2,\tau} \infeg C_{1}.
\end{eqnarray*}
\end{defperso}
\begin{propperso}{}\label{deltaasymptdelta0}
Let $g \in \G^{+}$ with $\delta \infeg 0$. Then $ \Delta$ is asymptotic to $\laplaz$ with a decaying rate $\delta$.
\end{propperso}
\textbf{Proof}:\begin{eqnarray}\label{expressionDelta}
\Delta = g^{ij} \n^{2}_{ij} &=& g^{ij} \nz^{2}_{ij} + g^{ij} (\n_{i} - \nz_{i}) \nz_{j} \nonumber\\
&=& g^{ij} \nz^{2}_{ij} - g^{ij} \ade{i}{k}{j} \, \nz_{k}.
\end{eqnarray}
The metrics $g$ and $\gz$ being equivalent, equation (\ref{holdcontcg}) directly gives
\begin{equation*}
\lambda \norm{\xi}{\gz}^{2} \infeg g^{ij}(x) \xi_{i}\xi_{j} \infeg \lambda^{-1} \norm{\xi}{\gz}^{2} , \forall x \in \m \; , \; \xi \in T\m.
\end{equation*}
Setting $$b^{k} = g^{ij} \ade{i}{k}{j},$$ then $b \in \sob{1}{2}{\delta}$ from (\ref{A}).
Given the Sobolev inequality, there exists a constant $C_{1} >0 $ such that
\begin{equation*}
\Norm{g^{ij} - \gz^{ij}}{1,6,\delta} + \Norm{b^{k}}{6,\delta} \infeg c \, (\Norm{g^{ij} - \gz^{ij}}{2,2,\delta} + \Norm{b^{k}}{1,2,\delta}) \infeg C_{1}. \cqfd
\end{equation*}
We will justify later on the terminology "asymptotic" used in this definition.\\
The operator $\A = - \Delta + n$ , acting on functions, will be of great interest.
\begin{propperso}{}\label{propelliptique}
Let $g\in \Gplus$ with $\delta\infeg 0$ and $\A = - \Delta + n$. Let $s \in \R $. There exists a constant $C = C(n,p,q,s,\delta,C_{1},\lambda)$ such that if $u \in \leb{2}{s}$ and $\A u \in \leb{2}{s}$ , then $u \in \sob{2}{2}{s}$ \corec{and}
\begin{equation*}
\Norm{u}{2,2,s} \infeg C \, \left(\Norm{\A u }{2,s} + \Norm{u}{2,s} \right).
\end{equation*}
\end{propperso}
\textbf{Proof}: By elliptic regularity, $u \in \sob{2}{2}{loc}$ and the estimate arises from interior estimates (see \cite{GT} for example) and scaling.$\cqfd$
%Le théorème suivant est l'analogue du Théorème 1.10 de %\cite{Bartnik1986} en asymptotiquement hyperbolique:
\begin{thmperso}{}\label{thcA}
Let $g\in \Gplus$ with $\delta\infeg 0$ so that $ \Delta$ is asymptotic to $\laplaz$.\\
Set $\A = - \Delta + n$ , with $n=3$.\\
Let $|s|< (n+1) \slash 2$. Then $\A: \sob{2}{2}{s}(\m) \rightarrow \leb{2}{s}(\T^{*} \otimes \Lambda^{3} T^{*}\m)$ is bounded.\\
Moreover, it satisfies the following elliptic estimate
\begin{equation}\label{sbes}
\Norm{u}{2,2,s} \infeg C \,\left(\Norm{\A u}{2,s} + \Norm{u}{2,s;\Omega_{R}} \right).
\end{equation}
In particular, $\A$ is a semi-Fredholm operator , $\ie$ $\A$ has finite dimensional kernel and closed range.
\end{thmperso}
\textbf{Proof}: We define the following operator norm:
$$||\Delta - \mathring\Delta||_{op}= \underset{\m}{\sup} \lbrace \Norm{(\Delta - \mathring\Delta)u}{2,s} : u\in \sob{2}{2}{s} , \Norm{u}{2,2,s} = 1 \rbrace$$
and $||\bullet||_{op,R}$ denotes the same norm restricted to functions supported in $E_{R} = \m \setminus \Omega_{R}$.\\
If $supp(u) \subset E_{R}$ , then from the expression (\ref{expressionDelta}) of $\Delta \, $,
\begin{eqnarray*}
\Norm{(\Delta - \mathring\Delta)u}{2,s} &\infeg& \Norm{(g^{ij} - \gz^{ij})\nz^{2}_{ij}u}{2,s} + \Norm{b^{k}\nz_{k}u}{2,s}\\
&\infeg& \underset{E_{R}}{\sup} \lbrace g^{ij} - \gz^{ij} \rbrace \Norm{\nz^{2}_{ij}u}{2,s} + \Norm{b^{k}\nz_{k}u}{2,s}\\
&\infeg& c \, \Norm{g - \gz}{\infty,0;E_{R}} \Norm{\nz^{2} u}{2,s} + \Norm{b \, \nz u}{2,s}.
\end{eqnarray*}
Using (\ref{h1}), Sobolev inequality and inclusion ($\delta \infeg 0$),
\begin{eqnarray*}
\Norm{(\Delta - \mathring\Delta)u}{2,s} 
%&\infeg& c \, \Norm{g - \gz}{\infty,\delta;E_{R}} \Norm{\nz^{2}u}{2,s} + \Norm{b}{6,\delta; E_{R}}\Norm{\nz u}{3,s}\\
&\infeg& c \, \big( \Norm{g - \gz}{2,2,\delta;E_{R}} + \Norm{b}{1,2,\delta;E_{R}} \big) \Norm{u}{2,2,s}.\\
\end{eqnarray*}
Recalling that  $\Norm{g - \gz}{2,2,\delta} + \Norm{b}{1,2,\delta}$ is bounded because $g\in\Gplus$,\\
 %et compte tenu que $\Norm{u}{2,2,s} = 1$ ,
\begin{equation}\label{opasympt}
||\Delta - \mathring\Delta||_{op,R} = o(1) \; \; \corec{\text{when}} \; \; R \rightarrow + \infty .
\end{equation}
This justifies {\it a posteriori} the terminology used in Definition \ref{defopasympt}.\\

Let $\chi_{R}$ be a cut-off function as in Definition \ref{cutoff}
\begin{eqnarray*}
\chi_{R} &=& \begin{cases}
1 \; \text{on} \; \Omega_{R \diagup 2} \\
0 \; \text{on} \; \m \setminus \Omega_{R}
\end{cases}
\end{eqnarray*}
Then we can decompose $u = u_{0} + u_{\infty}$ , with $u_{\infty}=(1-\chi_{R})u$.\\
We look into the operator $\Az = - \mathring\Delta + n $ acting on functions.
Using Corollary 3.13 of \cite{Andersson1993} with $\lambda = n \,$ , we obtain that $\forall |s|< (n+1) \slash 2 $ ,
$\Az : \sob{2}{2}{s} \rightarrow \leb{2}{s}$ is a Fredholm operator and an isomorphism. So there exists a positive constant $C= C(n,s)$ such that
\begin{equation}\label{1.22}
\Norm{u}{2,2,s} \infeg C \,\Norm{\Az u}{2,s}.
\end{equation}
Applying (\ref{1.22}) to $u_{\infty}$ , 
\begin{eqnarray}\label{1.28}
\Norm{u_{\infty}}{2,2,s} &\infeg& C \,\Norm{\Az u_{\infty}}{2,s} \nonumber\\
&\infeg&C \,\Norm{\A u_{\infty}}{2,s} + ||\Delta - \mathring\Delta||_{op,R}\Norm{u_{\infty}}{2,2,s}.
\end{eqnarray}
Yet $\A u_{\infty} = \A u - \A u_{0} = \A u \mathbf{- \chi_{R}\A u} \mathbf{+ \chi_{R}\A u} -  \A u_{0}$.\\
Thus,
\begin{eqnarray*}
\Norm{\A u_{\infty}}{2,s} &\infeg& \Norm{\A u}{2,s} + \Norm{\chi_{R}\A u }{2,s} +\Norm{\chi_{R}\A u -  \A u_{0}}{2,s}\\
&\infeg&C \,\Norm{\A u}{2,s} + \Norm{\chi_{R}\A u -  \A u_{0}}{2,s;\Omega_{R}}.
\end{eqnarray*}
From the expression (\ref{expressionDelta}) of $\Delta$,
\begin{eqnarray*}
\chi_{R}\A u -  \A u_{0} &=& -u \A \chi_{R} + n \chi_{R}u + 2 g^{ij}\d_{i}u \d_{j}\chi_{R}\\
&=& 2 g^{ij}\d_{i}u \, \d_{j}\chi_{R} + (g^{ij}\nz^{2}_{ij}\chi_{R} + b^{i}\d_{i}\chi_{R})u \, ,
\end{eqnarray*}
leading to
\begin{equation*}
\Norm{\chi_{R}\A u -  \A u_{0} }{2,s;\Omega_{R}} \infeg c \, \Norm{u}{1,2,s;\Omega_{R}}.
\end{equation*}
Finally
\begin{equation*}
\Norm{\A u_{\infty}}{2,s} \infeg C \,\left(\Norm{\A u}{2,s} + \Norm{u}{1,2,s;\Omega_{R}}\right).
\end{equation*}
Replacing (\ref{1.28}) and considering (\ref{opasympt}), we obtain for $R$ large enough
\begin{eqnarray}\label{1.28bis}
\Norm{u_{\infty}}{2,\delta} \infeg \Norm{u_{\infty}}{2,2,s} \infeg C \,\left(\Norm{\A u}{2,s} + \Norm{u}{1,2,s;\Omega_{R}} \right).
\end{eqnarray}
Using (\ref{1.28bis}) and the fact that on $\Omega_{R} , \norm{u_{0}}{\gz} \infeg \norm{u}{\gz}$
\begin{eqnarray*}
\Norm{u}{2,s} &\infeg& \Norm{u_{\infty}}{2,s} +  \Norm{u_{0}}{2,s}\\
&\infeg& C \,\left(\Norm{\A u}{2,s} + \Norm{u}{1,2,s;\Omega_{R}} \right) + \Norm{u}{2,s;\Omega_{R}}.
\end{eqnarray*}
Thanks to Ehrling inequality (\ref{Ehrling}),
\begin{eqnarray*}
\Norm{u}{2,s} &\infeg& C \,\left(\Norm{\A u}{2,s} + \Norm{u}{2,s;\Omega_{R}} \right) + \epsilon \Norm{u}{2,2,s;\Omega_{R}}\\
&\infeg& C \,\left(\Norm{\A u}{2,s} + \Norm{u}{2,s;\Omega_{R}} \right) + \epsilon \Norm{u}{2,2,s}.
\end{eqnarray*}
and we conclude with Proposition \ref{propelliptique}. $\cqfd$

Let $\Bz = - \mathring\Delta + n-1$ be an operator acting on $1$-forms.
\begin{thmperso}{}\label{thcB}
Let $\delta\infeg 0$ and $g\in\Gplus$.
Setting $B = - \Delta + n-1$ and $|s|< (n+1) \slash 2$.\\
Then $B: \sob{2}{2}{s}(T^{*}\m) \rightarrow \leb{2}{s}(\T^{*} \otimes \Lambda^{3} T^{*}\m)$ is bounded. Furthermore, it satisfies
\begin{equation}\label{sbeB}
\Norm{Y}{2,2,s} \infeg C \,\left(\Norm{BY}{2,s} + \Norm{Y}{2,s;\Omega_{R}} \right).
\end{equation}
In particular, $B$ is a semi-Fredholm operator , $\ie$ $B$ has finite dimensional kernel and closed range.
\end{thmperso}
\textbf{Proof}: From Proposition E of \cite{Lee2006}, the indicial radius of $\Bz$ is $(n+1) \slash 2$ and by Theorem C of \cite{Lee2006}, $\forall |s|< (n+1) \slash 2 $ , $\Bz : \sob{2}{2}{s}\rightarrow \leb{2}{s}$ is a Fredholm operator. By Corollary 3.13 of \cite{Andersson1993} , $\Bz$ is an isomorphism for the same span of $s$. So there exists a positive constant $C= C(n,s)$ such that
\begin{equation}\label{1.22B}
\Norm{Y}{2,2,s} \infeg C \,\Norm{\Bz Y}{2,s}.
\end{equation}
The proof is nearly identical to the one of Theorem \ref{thcA} with (\ref{1.22B}) replacing (\ref{1.22})$.\cqfd$
\begin{thmperso}{}\label{thc}
Let $\delta  \in \rbrack -(n+1) \slash 2 , 0 \rbrack$ with $n=3$. Then the operator\\ $F: \sobd{2}(\T) \rightarrow \leb{2}{\delta}(\T^{*} \otimes \Lambda^{3} T^{*}\m) :=\L^{*}$ is bounded. Furthermore, it verifies
\begin{equation}\label{sbedelta}
\Norm{(y,Y)}{2,2,\delta} \infeg C \, \left(\Norm{F(y,Y)}{2,\delta} + \Norm{(y,Y)}{2,0} + \Norm{(y,Y)}{2,\delta;\Omega_{R}}\right).
\end{equation}
In particular, $F$ is a semi Fredholm operator, $\ie$ $F$ has finite dimensional kernel and closed range.
\end{thmperso}
\textbf{Proof}:
Starting from the definition of $F$, the Triangle inequality together with (\ref{16g}) and the Sobolev inclusion (with $\delta \infeg 0$) directly yield
$$\Norm{F(y,Y)}{2, \delta}\infeg C \, \Norm{(y,Y)}{2,2,\delta} ,$$
where $C$ is a constant depending on $\gz$ and $\Norm{g}{\F}$.\\
Hence $F$ is a bounded (continuous) operator. Plugging the expression of $F_{0}(y, Y^{i})$ in (\ref{sbes}) and using Hölder inequality (\ref{h1}) , (\ref{uinfd}) ,Ehrling inequality (\ref{Ehrling}) along with Sobolev inclusion (with $\delta \infeg 0$) and $\phiogpi\in \lebd$,
\begin{eqnarray}\label{fredA}
\Norm{y}{2,2,\delta} &\infeg& C \, \left(\Norm{-\Delta y + ny}{2,\delta} + \Norm{y}{2,\delta;\Omega_{R}}\right) \nonumber\\
%&\infeg& C \, \left(\Norm{F_{0}(y, Y)}{2,\delta} + \Norm{ \phiogpi y}{2,\delta} + \Norm{(y,Y)}{2,0} + \Norm{\corec{\div }Y}{2,\delta} + \epsilon \Norm{Y}{2,2,0} + \Norm{y}{2,\delta;\Omega_{R}}\right) \nonumber\\
&\infeg& C \, \left(\Norm{F_{0}(y, Y)}{2,\delta} + \Norm{(y,Y)}{2,0} + \Norm{Y}{2,2,\delta} + \Norm{y}{2,\delta;\Omega_{R}}\right).
\end{eqnarray}
Plugging the expression of $F_{i}(y, Y^{i})$ in (\ref{sbeB}) and using Hölder inequality (\ref{h1}) , (\ref{uinfd}) ,Ehrling inequality (\ref{Ehrling}) along with Sobolev inclusion (with $\delta \infeg 0$) and $\phiigpi\in \lebd$,
\begin{eqnarray}\label{fredB}
\Norm{Y}{2,2,\delta} &\infeg& C \, \left(\Norm{-\Delta Y + (n-1)Y}{2,\delta} + \Norm{Y}{2,\delta;\Omega_{R}}\right) \nonumber\\
%&\infeg& C \, \left(\Norm{F_{i}(y, Y)}{2,\delta} + \Norm{ \phiigpi y}{2,\delta} + \Norm{(y,Y)}{2,0} + \epsilon \Norm{Y}{2,2,0} + \Norm{Y}{2,\delta;\Omega_{R}}\right) \nonumber\\
%&\infeg& C \, \left(\Norm{F_{i}(y, Y)}{2,\delta} + \Norm{(y,Y)}{2,0} + \epsilon \Norm{Y}{2,2,\delta} + \Norm{Y}{2,\delta;\Omega_{R}}\right) \nonumber\\
&\infeg& C \, \left(\Norm{F_{i}(y, Y)}{2,\delta} + \Norm{(y,Y)}{2,0} + \Norm{Y}{2,\delta;\Omega_{R}}\right).
\end{eqnarray}
Finally, combination of (\ref{fredA}) and (\ref{fredB}) gives (\ref{sbedelta}). For all $\delta  \in \rbrack -(n+1) \slash 2 , 0 \rbrack$ , the estimate (\ref{sbedelta}) verified by $F$ is analoguous to the one of Theorem \ref{estellip} and by a similar proof, we show $F$ is semi-Fredholm, $\ie$ $F$ has finite dimensional kernel and closed range.$\cqfd$\\[.3cm]
Now $F$ and its adjoint $F^{*}$ have similar structure ($F$ is formally self-adjoint)
$$F^{*}: \leb{2}{-\delta}(\T)\rightarrow \sob{-2}{2}{-\delta}(\T^{*} \otimes \Lambda^{3}T^{*}\m).$$
Let $\w{F^{*}}$ be the restriction of $F^{*}$ defined as follows
$$\w{F^{*}}: \sob{2}{2}{-\delta}(\T) \rightarrow \leb{2}{-\delta}(\T^{*} \otimes \Lambda^{3}T^{*}\m).$$
We can apply Theorem \ref{thc} to $\w{F^{*}}$: 

\begin{thmperso}{}\label{thcftilde}
Let $\delta  \in \rbrack -(n+1) \slash 2 , 0 \rbrack$ with $n=3$. Then the operator\\ $\w{F^{*}}: \sob{2}{2}{-\delta}(\T) \rightarrow \leb{2}{-\delta}(\T^{*} \otimes \Lambda^{3} T^{*}\m)$ is bounded. Furthermore, it satisfies
\begin{equation}\label{sbeftilde}
\Norm{(y,Y)}{2,2,-\delta} \infeg C \, \left(\Norm{\w{F^{*}}(y,Y)}{2,-\delta} + \Norm{(y,Y)}{2,-2 \delta} + \Norm{(y,Y)}{2,-\delta;\Omega_{R}}\right).
\end{equation}
In particular,  $\w{F^{*}}$ is a semi-Fredholm operator , $\ie$ $\w{F^{*}}$ has finite dimensional kernel and closed range.
\end{thmperso}
\textbf{Proof}:
From the definition of $\w{F^{*}}$, the Triangle inequality together with (\ref{16g}) and the Sobolev inclusion with $\delta \infeg 0$ directly yield
$$\Norm{\w{F^{*}}(y,Y)}{2, -\delta}\infeg C \, \Norm{(y,Y)}{2,2,-\delta}.$$
 Plugging the expression of $\w{F^{*}_{0}}(y, Y^{i})$ (formally identical to $F_{0}(y, Y^{i}))$ in (\ref{sbes}) and using Hölder inequality (\ref{h1}) , (\ref{uinfd}) ,Ehrling inequality (\ref{Ehrling}) along with Sobolev inclusion (with $\delta \infeg 0$) and $\phiogpi\in \lebd$,
\begin{eqnarray}\label{fredAtilde}
\Norm{(y,Y)}{2,2,-\delta} &\infeg& C \, \left(\Norm{-\Delta y + ny}{2,-\delta} + \Norm{N}{2,-\delta;\Omega_{R}}\right) \nonumber\\
%&\infeg& C \, \left(\Norm{\w{F^{*}_{0}}(y, Y)}{2,-\delta} + \Norm{ \phiogpi y}{2,-\delta} + \Norm{(y,Y)}{2,-2\delta} + \Norm{\corec{\div }Y}{2,-\delta} + \epsilon \Norm{Y}{2,2,-2\delta} + \Norm{y}{2,\delta;\Omega_{R}}\right) \nonumber\\
&\infeg& C \, \left(\Norm{\w{F^{*}_{0}}(y, Y)}{2,\delta} + \Norm{(y,Y)}{2,-2\delta} + \Norm{Y}{2,2,-\delta} + \Norm{y}{2,-\delta;\Omega_{R}}\right).
\end{eqnarray}
Plugging the expression of $\w{F^{*}_{i}}(y, Y^{i})$ (formally identical to $F_{i}(y, Y^{i}))$ in (\ref{sbeB}) and using Hölder inequality (\ref{h1}) , (\ref{uinfd}) ,Ehrling inequality (\ref{Ehrling}) along with the Sobolev inclusion (with $\delta \infeg 0$) and $\phiigpi\in \lebd$,
\begin{eqnarray}\label{fredBtilde}
\Norm{Y}{2,2,-\delta} &\infeg& C \, \left(\Norm{-\Delta Y + (n-1)Y}{2,-\delta} + \Norm{Y}{2,-\delta;\Omega_{R}}\right) \nonumber\\
%&\infeg& C \, \left(\Norm{\w{F^{*}_{i}}(y, Y)}{2,-\delta} + \Norm{ \phiigpi y}{2,-\delta} + \Norm{(y,Y)}{2,-2\delta} + \epsilon \Norm{Y}{2,2,-2\delta} + \Norm{Y}{2,-\delta;\Omega_{R}}\right) \nonumber\\
%&\infeg& C \, \left(\Norm{\w{F^{*}_{i}}(y, Y)}{2,-\delta} + \Norm{(y,Y)}{2,-2\delta} + \epsilon \Norm{Y}{2,2,-\delta} + \Norm{Y}{2,-\delta;\Omega_{R}}\right) \nonumber\\
&\infeg& C \, \left(\Norm{\w{F^{*}_{i}}(y, Y)}{2,-\delta} + \Norm{(y,Y)}{2,-2\delta} + \Norm{Y}{2,-\delta;\Omega_{R}}\right).
\end{eqnarray}
Finally, combination of (\ref{fredAtilde}) and (\ref{fredBtilde}) gives (\ref{sbeftilde}).\\
Similarly to $F$, for all $\delta \in \rbrack -(n+1) \slash 2 , 0 \rbrack$ , $\w{F^{*}}$ is a semi-Fredholm operator. $\cqfd$\\

\corec{We are now in possession of all the tools necessary to finish the proof 
of Theorem \ref{mainthmmanifold}.
}\\
 By elliptic regularity, $\text{ker} F^{*} = \text{ker} \w{F^{*}}$ is also finite dimensional.
%De plus , Coker $F = \faktor{Y}{ImF}$ est également de dimension finie.
If we apply Lemma \ref{lemiso} to $F$ , we get $$(\text{Coker} F)^{*} \simeq \text{ker} F^{*}.$$ So $(\text{Coker} F)^{*}$ is finite dimensional so $(\text{Coker} F)^{*} \simeq \text{Coker} F$. Thus we have the isomorphism
$$\text{Coker} F = \quotient{\L^{*}}{\text{Im} F} \simeq \text{ker} F^{*}.$$
The operator $F$ satisfies
$$\text{Im} F \subset \text{Im} D \phigpi \subset \L^{*}.$$
%et on peut écrire $$\L^{*} = \text{Im} F \:\oplus \text{ker} F^{*}$$ 
Let $\pi$ be the canonic projection:
$$\pi : \L^{*} \rightarrow \quotient{\L^{*}}{\text{Im} F}.$$
$\pi (\text{Im} D \phi)$ is closed, \corec{being the} subspace of a finite dimensional vector space. $\text{Im} (D \phi)$ is closed, being the inverse image of a closed set by a continuous map.
\corec{This ends the proof of the manifold structure of $\contrainte$, as a smooth submanifold of $\F$. In fact, all level sets of $\phigpi$ are smooth submanifolds of $F$. $\cqfd$}

%bibliographie bibtex
\bibliography{biblioarticle}

\begin{thebibliography}{10}

\bibitem{Andersson1993}
L.~Andersson.
\newblock Elliptic systems on manifolds with asymptotically negative curvature.
\newblock {\em Indiana University Mathematics Journal}, 42(4):1359--1387, 1993.

\bibitem{Andersson1996}
L.~Andersson and P.~Chrusciel.
\newblock Solutions of the constraint equations in general relativity
  satisfying "hyperboloïdal boundary conditions".
\newblock {\em Dissertationes Mathematicae}, 355:1--100, 1996.

\bibitem{Bartnik1986}
R.~Bartnik.
\newblock The mass of an asymptotically flat manifold.
\newblock {\em Comm. Pure Appl. Math.}, 34:661--693, 1986.

\bibitem{Bartnik2005}
R.~Bartnik.
\newblock Phase space for the {E}instein equations.
\newblock {\em Comm; Anal. Geom}, 13(5):845--885, 2005.

\bibitem{Chrusciel2004}
P.~Chrusciel and E.~Delay.
\newblock On mapping properties of the general relativistic constraints
  operator in weighted functions space, with applications.
\newblock {\em arXiv:gr-qc/0301073}, 2004.

\bibitem{ChruscielDelayExotic}
P.~Chrusciel and E.~Delay.
\newblock Exotic hyperbolic gluings.
\newblock {\em Journal of Differential Geometry, [hal-01233383]}, to appear.

\bibitem{ChruscielHerzlich2003}
P.~T. Chrusciel and M.~Herzlich.
\newblock The mass of asymptotically hyperbolic {R}iemannian manifold.
\newblock {\em Pac. J. Mah.}, 212:231--264, 2003.

\bibitem{DGS2013}
M.~Dahl, R.~Gicquaud, and A.~Sakovich.
\newblock Asymptotically hyperbolic manifolds with small mass.
\newblock {\em Comm. Math. Phys.}, 325(2):757--801, 2014.

\bibitem{Fischer1979}
A.~E. Fischer and J.~E. Marsden.
\newblock Topics in the dynamics general relativity.
\newblock {\em Italian physical society}, pages 322--395, 1979.

\bibitem{JFthese}
J.~Fougeirol.
\newblock {\em Etude des équations de contraintes sur une variété
  Asymptotiquement hyperbolique}.
\newblock PhD thesis, University of Avignon, in preparation.

\bibitem{GT}
D.~Gilbarg and N.~S. Trudinger.
\newblock {\em Elliptic Partial Differential Equations of Second Order}.
\newblock Springer, 2001.

\bibitem{Herzlich2005}
M.~Herzlich.
\newblock Mass formulae for asumptotically hyperbolic metrics, in {T}he
  {A}d{S}/{CFT} correspondence ({E}instein metrics and conformal geometry).
\newblock {\em IRMA Lect. in Math. Theor. Physics.}, 8:103--121, 2005.
\newblock Eur. Math. Soc. , Zurich.

\bibitem{Hille1957}
E.~Hille and R.~Phillips.
\newblock Functional analysis and semi-groups.
\newblock {\em Colloquim Series, Am. Math. Soc}, volume 31, 1957.

\bibitem{Lee2006}
J.~M. Lee.
\newblock {F}redholm operators and {E}instein metrics on conformally compact
  manifolds.
\newblock {\em Mem. Amer. Math. Soc}, 183(864), July 2006.

\bibitem{Mazzeo1991}
R.~Mazzeo.
\newblock Unique continuation at infinity and embedded eigenvalues for
  asymptotic hyperbolic manifolds.
\newblock {\em American Journal of Mathematics}, 113:25--45, 1991.

\bibitem{McCormick2014}
S.~McCormick.
\newblock The phase space for the {E}instein-{Y}ang-{M}ills equations and the
  first law of black hole thermodynamics.
\newblock {\em Adv. Theor . Math. Phys}, 18(4):799--825, 2014.

\bibitem{McCormick2015}
S.~McCormick.
\newblock The {H}ilbert manifold of asymptotically flat metric extensions.
\newblock {\em arXiv:1512.02331v1}, 2015.

\bibitem{RaiSaraykar2016}
J.~H. Rai and R.~Saraykar.
\newblock {H}ilbert manifold structure of the set of solutions of constraint
  equations for coupled {E}instein and scalar fields.
\newblock {\em arXiv:1605.08858v1}, 2016.

\bibitem{SergiuKlainerman2013}
I.~R. Sergiu~Klainerman and J.~Szeftel.
\newblock Overview of the proof of the bounded {L}$^2$ curvature conjecture.
\newblock {\em arXiv:1204.1772v2}, 2013.

\end{thebibliography}
\bibliographystyle{abbrv}
%\addcontentsline{toc}{chapter}{Bibliographie}

\end{document}